\documentclass[12pt]{amsart}
\usepackage{amsmath,amsfonts,latexsym,graphicx,amssymb,url}
\usepackage{hyperref}
\usepackage{amsmath,txfonts,pifont,bbding,pxfonts}
\usepackage[active]{srcltx}
\usepackage{wasysym,pstricks,pdfsync}

\newcommand{\dist}{\text{dist}} 

\newcommand{\R}{{\mathbb R}} 

\newcommand{\manuscriptdate}{\today}
\renewcommand{\(}{\left(}
\renewcommand{\)}{\right)}
\theoremstyle{plain}
\newtheorem{theorem}{Theorem}[section]
\newtheorem{corollary}[theorem]{Corollary}
\newtheorem{lemma}[theorem]{Lemma}
\newtheorem{proposition}[theorem]{Proposition}
\newtheorem{definition}[theorem]{Definition}
\newtheorem{remark}[theorem]{Remark}

\begin{document}

\title[Partial regularity for fully nonlinear equations ~~~\manuscriptdate]{Partial $C^{2,\alpha}$ and $W^{2, VMO}$ regularity for nondifferentiable fully nonlinear elliptic equations and Isaacs equations}
\author[C. E. Guti\'errez and Qingbo Huang\\~~~\tiny{\manuscriptdate}]{Cristian E. Guti\'errez and Qingbo Huang\\\manuscriptdate}
\address{Department of Mathematics\\Temple University\\Philadelphia, PA 19122}
\email{gutierre@temple.edu}
\address{Department of Mathematics \& Statistics, Wright State University, Dayton, OH 45435}
\email{qingbo.huang@wright.edu}
\thanks{\manuscriptdate. \\The first author was partially supported by NSF grant DMS--1600578.}

\maketitle

\begin{abstract}

We investigate the partial regularity of viscosity solutions to fully nonlinear elliptic equations of the form $F(D^2u, x) = f(x)$ in previously unexplored settings:

(i) When $F$ is non-differentiable with respect to $D^2u$.
(ii) When either $F$ or $f$ lacks Dini continuity in $x$.

In the first scenario, we introduce a tangent operator property for $F$ such that, for a fixed point $x_0$ and a given matrix $M_0$, $F(M, x_0)$ is tangential to a convex or concave operator at $M_0$. We establish partial $C^{2,\alpha}$ regularity for equations exhibiting this property.

In the second scenario, we establish new partial $W^{2, \text{VMO}}$ regularity provided that the oscillation of $F$ in $x$ under the $L^n$-norm is of order $o(\log^{-1}(1/r))$.

We address the Hausdorff dimensions of the singular sets of solutions for equations of the form $F(D^2u)=f(x)$ with data $f$ in $W^{1,n}$.

Additionally, we demonstrate that the tangent operator property holds for various classes of equations, including nondifferentiable twisted elliptic equations, infimum of $C^1$ equations, semiconvex equations, and piecewise differentiable equations.

Finally, we tackle the issue of smooth enveloping and obtain partial regularity of solutions for three subclasses of Isaacs equations.
\end{abstract}

\maketitle

\tableofcontents

\setcounter{equation}{0}
\section{Introduction and main results}

This paper addresses the partial regularity of viscosity solutions to  
fully nonlinear elliptic equations of form
\begin{equation}\label{eq:main equation}
F(D^2u, x) = f(x)
\end{equation}
where $F:\mathcal S\times \Omega\longrightarrow \R$ is continuous but not necessarily differentiable, with $\Omega \subset \R^n$ a bounded domain and $\mathcal S$ is the space of real $n\times n$ symmetric matrices, satisfying the following uniform ellipticity condition with constants 
$0<\lambda \le \Lambda$:
\begin{equation}\label{eq:ellipticity condition}
\lambda \|N\|\leq F(M+N,x)-F(M,x)\leq \Lambda \|N\|
\end{equation}
for each $M\in \mathcal S$, $x\in \Omega$, and for all $N\in \mathcal S$ with $N\geq 0$. Here $\|N\|=\sup_{|x|=1}|Nx|$. 

It is well known, as shown in \cite{NV}, that 
for general non-convex fully nonlinear elliptic equations like \eqref{eq:main equation}, one cannot expect to obtain solutions with regularity higher than $C^{1, \varepsilon}$.

Recent research has also examined the partial regularity of viscosity solutions for differentiable, fully nonlinear elliptic equations.   
In \cite{ASS}, the authors initially established partial $C^{2,\alpha}$ regularity, along with an estimate for the Hausdorff dimension of the singular set, for viscosity solutions of the equation $F(D^2u)=0$ when $F$ belongs to the class $C^1$.  
This result was extended in \cite{dPT} 
to $F(D^2u)=f(x)$ under the assumptions that $F\in C^1$ and $f$ is 
Lipschitz continuous.  
One component to show the partial 
$C^{2,\alpha}$ regularity is the regularity of flat viscosity solutions 
which itself is also very interesting.   
For $F(D^2u, x)=f(x)$ where $F$  
is $C^1$ in $M$ and both $F$ and $f$ are H\"older continuous in $x$, 
the regularity of flat viscosity solutions asserts that if the solution $u$ 
is a small perturbation of a quadratic polynomial solution, then $u$ is
$C^{2,\alpha}$;
this has been established in \cite{Sav}, \cite{ASS}, and \cite{dPT}.   
In case that $F$ is $C^1$ in $M$ and both $F$ and $f$ are Dini 
continuous in $x$, a $C^2$ regularity of flat viscosity solutions can be developed, see \cite{daSP}.         
Despite these significant advancements and techniques, 
the partial regularity of solutions and the regularity of flat solutions for \eqref{eq:main equation} 
remain undetermined in two cases when $F$  
is nondifferentiable and either $F$ or $f$ is not Dini continuous in $x$.

Important nondifferentiable equations include Isaacs equations 
arising from stochastic differential games (see \cite{Nis}).  
Also the natural setting for viscosity solutions of \eqref{eq:main equation} is 
that $F$ is elliptic and both $F$ and $f$ are continuous in $x$.     
These motivate us to investigate the partial regularity issue of
\eqref{eq:main equation} in the general setting. 

Our first purpose in this paper is to address partial $C^{2,\alpha}$-regularity for nondifferentiable equations having the tangent operator 
property defined as follows.  
\begin{definition}\label{def:tangent property operator}
Let $\omega:[0,\, \infty)\to [0,\infty)$ be a nondecreasing function such that 
$\omega(t)>0$ for $t>0$ and $\lim_{t\to 0^+}\omega(t)=\omega(0)=0$.
Given  
$F:\mathcal S\longrightarrow \R$ uniformly elliptic with ellipticity 
constants $\lambda$, $\Lambda$ and a set $\mathcal E\subset \mathcal S$, 
we say that $F$ has tangent operator property in $\mathcal E$ with modulus $\omega$, denoted by $F\in \mathbb T^\omega(\mathcal E, \lambda, \Lambda)$ 
(or $F\in \mathbb T^\omega(\mathcal E)$ when $\lambda$ and $\Lambda$ 
are clear from the context), 
if for any $M_0\in \mathcal E$, 
there exists a function $T_{M_0}:\mathcal S\longrightarrow \R$ convex or concave 
uniformly elliptic with ellipticity constants $\lambda$, $\Lambda$ 
such that 
$$
|F(M)-T_{M_0}(M)| \le \omega(\|M-M_0\|) \|M-M_0\|   
$$
for all $M\in \mathcal S$ with $\|M-M_0\| \le 1$.   
\end{definition}

Note that in this definition $T_{M_0}$ is not necessarily unique and 
the convexity of $T_{M_0}$ may vary as $M_0$ changes.  Geometrically, 
$F$ is tangent to $T_{M_0}$ at $M_0$.   
It is easy to verify that if $F$ is of class $C^1$ and $DF$ is uniformly 
continuous in $\mathcal S$, then $F\in \mathbb T^\omega(\mathcal S)$ for some
$\omega$.   Like $C^1$ operators, operators having the tangent 
operator property are in general neither convex nor concave.  
It will be shown in Sections \ref{sec:nondiff fully nonlinear satisfying (H2)} and \ref{sec:Isaacs equations} below that a variety of classes of nondifferentiable 
operators satisfy the tangent operator property.  

In the definition of $\mathbb T^\omega(\mathcal E)$, the convexity requirement
on $T_{M_0}$ can be weakened to semiconvexity.
This is because the tangency between two operators is a transitive 
relation and by Theorem \ref{thm:Theorem 1.4} below a semiconvex operator is tangent  
to some convex operator at any given point.     

Let $\mathcal B_k$ denote the ball in $\mathcal S$ centered at $0$ with radius $k$.  
To study partial $C^{2,\alpha}$ regularity, we will need the following 
assumption saying that $F(M,x)$ satisfies a local tangent operator property 
for $M\in \mathcal S$ uniformly in $x$; that is:

\medskip
\begin{enumerate}
  \item[\hypertarget{H1}{\bf(H1)}]
  There exists $\{\omega_k\}_{k\ge 1}$ such that
$F(\cdot, x)\in \mathbb T^{\omega_k}(\overline{\mathcal B}_k, \lambda, \Lambda)$ 
for every $k\ge 1$ and $x\in \Omega$, 
where $\omega_k(t)$ is nondecreasing in $t\ge 0$,  
$\omega_k(t)>0$ for $t>0$, $\lim_{t\to 0^+}\omega_k(t)=\omega_k(0)=0$,  
and $\omega_k\le \omega_{k+1}$ for $k\ge 1$.    
More precisely,  for each $M_0\in \overline{\mathcal B}_k$ and $x_0\in \Omega$, 
there exists a  function $T_{M_0, x_0}$ convex or concave 
and uniformly elliptic with ellipticity constants $\lambda$, $\Lambda$ 
such that 
$$
|F(M, x_0)-T_{M_0, x_0}(M)| \le \omega_k(\|M-M_0\|) \|M-M_0\|   
$$
for all $M\in \mathcal S$ with $\|M-M_0\| \le 1$. 
\end{enumerate}
\medskip

We are now ready to state our main result on partial $C^{2, \alpha}$ regularity. 
 
\begin{theorem}\label{thm:Theorem 1.1} 
Let $u$ be a viscosity solution of \eqref{eq:main equation} in $\Omega$.  
Assume that $F$ satisfies \eqref{eq:ellipticity condition} and (H1).  
Suppose that $f\in C^\alpha(\Omega)$ and $F$ satisfies 
\begin{equation}\label{eq:holder in x}
|F(M, x_2)-F(M, x_1)| \le A(1+\|M\|)\,|x_2-x_1|^\alpha  
\qquad\forall\ x_1, x_2\in \Omega,\ M\in \mathcal S, 
\end{equation}
for some $A>0$ and $0<\alpha<\alpha_0$ with $0<\alpha_0<1$ 
depending only on $n$, $\lambda$, and $\Lambda$.   
Then there exists open set $\Omega_1\subset \Omega$ such that 
$u\in C^{2, \alpha}(\Omega_1)$ with 
$|\Omega\setminus\Omega_1|=0$ where $|\cdot |$ 
denotes Lebesgue measure.      
\end{theorem}

Our second purpose in the paper is to establish partial 
$W^{2,\text{VMO}}$ regularity for \eqref{eq:main equation} with $F$ non-Dini 
continuous in $x$ but satisfying an appropriate oscillation condition in $x$. 
To state the next result, we consider the following global tangent operator assumption \hyperlink{H2}{\bf(H2)\rm}   
and introduce some notation:  

\medskip
{\noindent
\begin{enumerate}
\item[\hypertarget{H2}{\bf(H2)\rm}]  $F(\cdot, x)\in \mathbb T^{\omega}(\mathcal S, \lambda, \Lambda)$ 
for every $x\in \Omega$, 
where $\omega(t)$ is nondecreasing in $t\ge 0$,  
$\omega(t)>0$ for $t>0$, and $\lim_{t\to 0^+}\omega(t)=\omega(0)=0$. 
\end{enumerate}    
}
\medskip

Throughout this paper, $B_r(x_0)$ denotes the ball centered 
at $x_0$ with radius $r$, $B_r=B_r(0)$, and    
$f_E=\fint_{E}f\,dx$ denotes the average of $f$ over $E$.   
$W^{2,\text{VMO}}(\Omega)$ denotes the space of functions $u$ such that
$u$, $Du$, $D^2u$ are in $\text{VMO}(\Omega)$, that is, they have vanishing mean oscillation in $\Omega$.  
A constant $C=C(n, \lambda, \Lambda)$ depending only on 
$n$, $\lambda$, and $\Lambda$ is called a universal constant.   
Set 
\begin{equation}\label{eq:definition betaF}
\beta_F(x, x_0)=\sup_{M\in \mathcal S}  
\frac{|F(M, x)-F(M, x_0)|}{\|M\|+1} 
\end{equation}
and for simplicity $\beta_F(x)=\beta_F(x, 0)$.  
For a function $f$ on a set $H$, as usual, we denote $\text{osc}_H f=\sup_{x,y\in H}|f(x)-f(y)|$.
  
Next is our primary partial $W^{2,\text{VMO}}$-regularity result, which represents a new type of partial regularity even for continuously differentiable equations. Its proof introduces a novel combinantion of observations and techniques.

\begin{theorem}\label{thm:Theorem 1.2} 
Let $u$ be a viscosity solution of \eqref{eq:main equation} in $\Omega$.  
Assume that $F$ satisfies \eqref{eq:ellipticity condition} and (H2).  
Let $\iota(r)$ be a nonnegative function on $(0,\, 1]$ 
with $\lim_{r\to 0^+}\iota(r)=0$.  
Suppose that $f\in C(\Omega)$, $F\in C(\mathcal S\times \Omega)$, and  
\begin{equation}\label{eq:Ln norm of betaF}
\(\fint_{B_r(x_0)}[\beta_F(x, x_0)]^n\,dx\)^{1/n}
\le \frac{\iota(r)}{\log(2/r)}  
\end{equation}
for any $B_r(x_0)\subset \Omega$ with $r\le 1$.       
Then there exists an open set $\Omega_2\subset \Omega$ such that 
$u\in W^{2,\text{VMO}}_{loc}(\Omega_2)$ with 
$|\Omega\setminus \Omega_2|=0$.      
\end{theorem}

Notice that from Theorem \ref{thm:Theorem 1.2}, any viscosity solution of $F(D^2u)=f(x)$ 
with $F\in \mathbb T^\omega(\mathcal S)$ and $f\in C(\Omega)$ has partial 
$W^{2,\text{VMO}}$-regularity.   

Our next result, Theorem \ref{thm:Theorem 1.3}, concerns the
Hausdorff dimension of singular sets for solutions to $F(D^2u)=f(x)$ 
under assumption $f\in W^{1,n}(\Omega)$.   
The mechanism to obtain Hausdorff dimension estimates in Theorem \ref{thm:Theorem 1.3} 
is analogous to that in \cite{ASS}.   
Besides $W^{2,\text{VMO}}$ regularity for flat solutions,  
another key component for our argument are the 
$W^{3,\delta}$-estimates for $W^{1,n}$ data $f$, Theorem \ref{Theorem 4.4}, 
which also refines the corresponding $W^{3,\delta}$ estimates in 
\cite{ASS} for $f=0$ and in \cite{dPT} for $f$ Lipschitz.


\begin{theorem}\label{thm:Theorem 1.3} 
Let $u$ be a viscosity solution of $F(D^2u)=f(x)$ in $\Omega$.  
Assume that $f\in C(\Omega)\cap W^{1,n}(\Omega)$ and $F$ satisfies \eqref{eq:ellipticity condition} and (H2).       
Then there exists open set $\Omega_2\subset \Omega$ such that 
$u\in W^{2,\text{VMO}}_{loc}(\Omega_2)$ with       
$\mathcal H_{dim}(\Omega\setminus \Omega_2)\le n-\delta$, 
where $\mathcal H_{dim}(\Omega\setminus \Omega_2)$ is the Hausdorff dimension 
of $\Omega\setminus \Omega_2$ and $\delta$ depends only on $n$, 
$\lambda$, and $\Lambda$.  
\end{theorem}

Our third objective in the paper is to examine diverse classes of non-differentiable equations and Isaac’s equations that satisfy the tangent operator property. In this context, the classes of non-differentiable equations satisfying one of the following sets of conditions \hypertarget{A1}{\bf(A1)\rm}--\hypertarget{A4}{\bf(A4)\rm} are discussed.


\medskip
{\noindent
\begin{enumerate}
\item[\hypertarget{A1}{\bf(A1)\rm}] (Convex$+C^1$ equations) 
$F(M, x)=F_\cup(M, x)+F_C(M, x)$ on $\mathcal S\times \Omega$, 
where $F_\cup(M, x)$ is uniformly elliptic with constants 
$0<\lambda \le \Lambda$ and convex in $M$, $F_C(M, x)$ is 
degenerate elliptic with constants $0$, $\Lambda$ and of $C^1$ 
in $M$, and $D_MF_C$ has a modulus of continuity $\tau$ 
such that 
$$
|D_MF_C(M_2, x) - D_MF_C(M_1, x)| \le \tau(\|M_2-M_1\|)
$$
for any $M_1, M_2\in \mathcal S$, $x\in\Omega$.  
\end{enumerate}     
}
\medskip

The equations satisfying \hypertarget{A1}{\bf(A1)\rm} include non-differentiable twisted elliptic 
equations of weakly concave $F_C$.  

\medskip
{\noindent
\begin{enumerate}
\item[\hypertarget{A2}{\bf(A2)\rm}] ($\inf$ of $C^1$ equations)
$F(M, x)=\inf_{\gamma\in \Gamma} F_\gamma(M, x)$ 
on $\mathcal S\times \Omega$, 
where $\{F_\gamma\}_{\gamma\in \Gamma}$ is a family of fully nonlinear 
elliptic operators with ellipticity constants $0<\lambda \le \Lambda$,  
$\{F_\gamma(0, x)\}_{\gamma\in \Gamma}$ is uniformly bounded from
below, each $F_\gamma(M, x)$ is of $C^1$ in $M$, and 
$\{D_MF_\gamma\}$
has a uniform modulus of continuity $\tau$ such that 
$$
|D_MF_\gamma(M_2, x) - D_MF_\gamma(M_1, x)| \le \tau(\|M_2-M_1\|)
$$
for any $M_1, M_2\in \mathcal S$, $x\in\Omega$, $\gamma\in\Gamma$.   
\end{enumerate}        
}
\medskip

The equations satisfying \hypertarget{A2}{\bf(A2)\rm} include certain Isaacs equations.   

\medskip
{\noindent
\begin{enumerate}
\item[\hypertarget{A3}{\bf(A3)\rm}] (Semiconvex equations) 
$F(M, x)$ satisfies \eqref{eq:ellipticity condition} and is semiconvex in $M$ in the sense that
there exists a constant $K>0$ such that 
$F(M, x)+K|M|^2$ is convex on $\mathcal S$ for every $x\in \Omega$
with $|M|=(\sum m_{ij}^2)^{1/2}$ being the Euclidean norm of $M=[m_{ij}]$.  
\end{enumerate}     
}
\medskip

\medskip
{\noindent
\begin{enumerate}
\item[\hypertarget{A4}{\bf(A4)\rm}] (Piecewise differentiable equations) 
$F(M, x)$ satisfies \eqref{eq:ellipticity condition}, and there exist (bounded or unbounded) closed surfaces $\Xi_i , 1\le i \le N$ of class $C^1$ in the sense of (S) in Subsection \ref{subsect:piecewise-diff-equations} below such that $d(\Xi_i, \Xi_j)>2d_0$ for $i\neq j$, $\mathcal S\setminus \cup_{i=1}^N \Xi_i$ is partitioned into disjoint domains $\mathcal D_i$,  $1\le i \le N+1$, and on each $\mathcal D_i$, $F$ is $C^1$ in $M$ and $D_MF$ has a modulus of continuity $\tau$ so that 
$$ |D_MF(M_2, x)-D_MF(M_1, x)|\le \tau(\|M_2-M_1\|)$$
for $M_1, M_2\in \mathcal D_i, x\in \Omega$.   
\end{enumerate}
}
\medskip

\begin{theorem}\label{thm:Theorem 1.4} 
If $F$ satisfies one of \hyperlink{A1}{\bf(A1)\rm}-- \hyperlink{A4}{\bf(A4)\rm}, then 
$F(\cdot, x)\in \mathbb T^{\omega}(\mathcal S, \lambda', \Lambda')$ 
for every $x\in \Omega$, 
where constants $\lambda'$, $\Lambda'$ and modulus $\omega$ 
depending only on $n$ and the structure data in \hyperlink{A1}{\bf(A1)\rm}-- \hyperlink{A4}{\bf(A4)\rm}, and
$\omega(t)$ is nondecreasing in $t\ge 0$,  
$\omega(t)>0$ for $t>0$, and $\lim_{t\to 0^+}\omega(t)=\omega(0)=0$. 
\end{theorem}

According to Theorem \ref{thm:Theorem 1.4}, Theorems \ref{thm:Theorem 1.1}-- \ref{thm:Theorem 1.3} hold for \eqref{eq:main equation} 
satisfying one of the conditions \hyperlink{A1}{\bf(A1)\rm}-- \hyperlink{A4}{\bf(A4)\rm}.  The proof of Theorem \ref{thm:Theorem 1.4} for the cases of
\hyperlink{A3}{\bf(A3)\rm} and \hyperlink{A4}{\bf(A4)\rm} is more involved than it appears.    

To place our results in perspective, we enumerate some recent work on everywhere regularity and a priori estimates for several special non-convex equations.
In recent years, regularity theories of second derivatives for viscosity solutions 
to various 
non convex fully nonlinear elliptic equations has been developed.  
In fact, the $C^{2,\alpha}$ regularity theory 
was established in 
\cite{CC2} 
for the equation $\min\{F^\cap(M,x), F^\cup(M,x)\}=f(x)$ with $F^\cap$ concave and $F^\cup$ convex in $M$.
%
%
For equations exhibiting approximate convexity 
at $M=\infty$, various regularity theories were developed under  
suitable topologies of approximation.  
In these instances, the $W^{2,p}$ regularity theory of viscosity solutions was addressed in \cite{Kr1}, \cite{Kr2}, and \cite{PT}. The $W^{2,p, \mu}$ and $C^{2,\alpha}$ regularity theories were obtained in \cite{H2}. Additionally, the $W^{2,\text{BMO}}$ regularity theory of viscosity solutions was established in \cite{H2}, and the $C^{2,\alpha}$ regularity theory for quasiconcave equations is established in \cite{Gof}.


Apriori estimates for the second derivatives of solutions to several non-convex fully nonlinear elliptic equations have also been proved.     
The work \cite{CY} 
establishes $C^{2,\alpha}$ estimates for $F(D^2u)=0$ when the set 
$\{M: F(M)=0\} \cap \{M: tr(M)=t\}$ is a strictly convex variety in 
$\mathcal S$ for every $t$. 
In \cite{SW},  $C^{2,\alpha}$ estimates for the twisted Monge-Amp\`ere 
equation $\log\det D^2_xu-\log\det (-D^2_yu)=0$ 
were established.   
In \cite{Col}, $C^{2,\alpha}$ estimates for general $C^2$ twisted elliptic equations, such as $F(D^2u)=0$, are presented. The operator $F$ is twisted if it can be expressed as $F=F_\cup + F_\cap$, where $F_\cup$ is convex and $F_\cap$ is weakly concave. \cite{Col} also noted that several fully nonlinear equations arising from differential geometry fall under the category of twisted equations.   
In \cite{DK}, weighted $W^{2,p}$ estimates for equation \eqref{eq:main equation} were established under suitable approximate convexity conditions.
It is important to note that these a priori estimates for fully nonlinear equations cannot directly yield corresponding regularity without additional work.

We now turn to the Isaacs equations of form 
\begin{equation}\label{eq:Isaacs equation}
\inf_{\gamma\in\mathcal G}\sup_{\beta\in\mathcal B}
[\text{tr\,}(A_{\beta\gamma}(x)D^2u)-f_{\beta\gamma}(x)]=0,
\end{equation}
where $A_{\beta\gamma}(x)\in {\mathcal A}_{\lambda, \Lambda}$ for 
$\beta\in \mathcal B$, $\gamma\in \mathcal G$, $x\in \Omega$, and   
${\mathcal A}_{\lambda, \Lambda}$ denotes all symmetric matrices 
whose eigenvalues belong to $[\lambda,\,\Lambda]$. 

Progress regarding regularity for Isaacs equations has been limited, with no partial regularity results previously established. However, recent investigations have addressed Hessian estimates and regularity for several specific Isaacs equations.   
A $C^{2,\alpha}$ regularity theory for   
the 3-operator equation $\min\{L_1u, \max\{L_2u, L_3u\}\}=0$ was 
obtained in \cite{CC2}, where 
$L_i$ are linear operators.  
In \cite{Pim}, a $W^{2,p}$ regularity of viscosity 
solutions to the equation 
\[
\sup_{\alpha\in\mathcal A}\inf_{\beta\in\mathcal B}
[-\text{tr\,}(A_{\alpha,\beta}(x)D^2u)-b_{\alpha,\beta}(x)Du]=f
\]
was considered under the assumption that $A_{\alpha,\beta}(x)$ are close to 
$\overline A_{\beta}(x)$, i.e., when this equation is close to a Bellman equation.  
In \cite{Kov},  $W^{2,p}$ estimates were established 
for $C^2$ solutions of the Isaacs equation 
$\max_{y\in Y}\min_{z\in Z}\text{tr}[a(y,z)D^2u]=0$ when 
the control matrix $a(y,z)=a^1(y)+a^2(z)$ is diagonal, 
and $a^1(y)=\text{diag}[a^1_{ii}]$ with $a^1_{kk}=y$ and constant 
$a^1_{ii}$ ($i\neq k$).
In \cite[Sect. 6.3]{H2}, $W^{2,\text{BMO}}$ regularity was proved for the Isaacs 
equations of the form 
$\inf_{N\in\mathcal S}\sup_{A\in\mathcal A_{\lambda, \Lambda}}
[\text{tr}(AD^2u+(B-A)N)+Z(A,N)]=0$ under the assumption  
$\lim_{N\to\infty}D_NZ=0$. 

Our objective is to examine the problem of smooth enveloping and to analyze specific classes of
Isaac’s equations are classified as twisted elliptic, half-$C^1$, or piecewise differentiable.  
We first introduce two structural conditions that yield
the local tangent operator property \hyperlink{H1}{\bf(H1)\rm}.

\medskip
\begin{enumerate}
\item[\hypertarget{I1}{\bf(I1)\rm}]
$A_{\beta\gamma}(x)\in\mathcal A_{\lambda,\Lambda}$ for every
$\beta\in\mathcal B$, $\gamma\in\mathcal G$, and $x\in\Omega$, and there
exist a constant $C_0>0$ and a modulus $h$ such that
$|f_{\beta\gamma}(x)|\le C_0$ and, for every $\gamma\in\mathcal G$ and
$x\in\Omega$, the set
\[
\mathcal E_{\gamma,x}
=\{(A_{\beta\gamma}(x),-f_{\beta\gamma}(x)):\beta\in\mathcal B\}
\]
satisfies \hyperlink{LSE}{\bf(LSE)\rm} with modulus $h$ as defined in
Subsection \ref{subsec:locally-smooth-enveloping} below.

\item[\hypertarget{I2}{\bf(I2)\rm}]
$A_{\beta\gamma}=A^I_\beta+A^{II}_\gamma$ and
$f_{\beta\gamma}=f^I_\beta+f^{II}_\gamma$, where
$A^I_\beta(x)\in\mathcal A_{0,\Lambda}$ and
$A^{II}_\gamma(x)\in\mathcal A_{\lambda,\Lambda}$ for every
$\beta\in\mathcal B$, $\gamma\in\mathcal G$, and $x\in\Omega$.
In addition, there exist a constant $C_0>0$ and a modulus $h$ such that
\[
|f^I_\beta(x)|+|f^{II}_\gamma(x)|\le C_0
\]
and, for every $x\in\Omega$, the set
\[
\mathcal E_x=\{(A^I_\beta(x),-f^I_\beta(x)):\beta\in\mathcal B\}
\]
satisfies \hyperlink{LSE}{\bf(LSE)\rm} with modulus $h$.
\end{enumerate}
\medskip

\begin{theorem}\label{thm:Theorem 1.5}
If \hyperlink{I1}{\bf(I1)\rm} or \hyperlink{I2}{\bf(I2)\rm} holds, then
the Isaacs operator
\begin{equation}\label{eq:Isaacs operator}
F(M,x)=\inf_{\gamma\in\mathcal G}\sup_{\beta\in\mathcal B}
\bigl[\text{tr}(A_{\beta\gamma}(x)M)-f_{\beta\gamma}(x)\bigr]
\end{equation}
satisfies \hyperlink{H1}{\bf(H1)\rm} with ellipticity parameters
$\lambda'$, $\Lambda'$ and moduli $\omega_k$, $k\ge1$, depending only on
$n$ and the structural data in \hyperlink{I1}{\bf(I1)\rm} or
\hyperlink{I2}{\bf(I2)\rm}.
\end{theorem}

We next introduce three structural conditions that yield the global tangent
operator property \hyperlink{H2}{\bf(H2)\rm}.

\medskip
\begin{enumerate}
\item[\hypertarget{I3}{\bf(I3)\rm}]
$A_{\beta\gamma}(x)\in\mathcal A_{\lambda,\Lambda}$ and
$|f_{\beta\gamma}(x)|\le C_0$ for every $\beta\in\mathcal B$,
$\gamma\in\mathcal G$, and $x\in\Omega$, where $C_0$ is a constant.
In addition, for every $\gamma\in\mathcal G$ and $x\in\Omega$, the set
\[
\mathcal E_{\gamma,x}
=\{(A_{\beta\gamma}(x),-f_{\beta\gamma}(x)):\beta\in\mathcal B\}
\]
is closed and satisfies \hyperlink{GSE}{\bf(GSE)\rm} with modulus $\tau$
as defined in Subsection \ref{subsec:global-smooth-enveloping} below.

\item[\hypertarget{I4}{\bf(I4)\rm}]
$A_{\beta\gamma}=A^I_\beta+A^{II}_\gamma$ and
$f_{\beta\gamma}=f^I_\beta+f^{II}_\gamma$, where
$A^I_\beta(x)\in\mathcal A_{0,\Lambda}$,
$A^{II}_\gamma(x)\in\mathcal A_{\lambda,\Lambda}$, and
\[
|f^I_\beta(x)|+|f^{II}_\gamma(x)|\le C_0
\]
for every $\beta\in\mathcal B$, $\gamma\in\mathcal G$, and $x\in\Omega$,
where $C_0$ is a constant.  In addition, for every $x\in\Omega$, the set
\[
\mathcal E_x=\{(A^I_\beta(x),-f^I_\beta(x)):\beta\in\mathcal B\}
\]
is closed and satisfies \hyperlink{GSE}{\bf(GSE)\rm} with modulus $\tau$.

\item[\hypertarget{I5}{\bf(I5)\rm}]
$\mathcal B$ and $\mathcal G$ are finite sets, $A_{\beta\gamma}$ and
$f_{\beta\gamma}$ are independent of $x$, and
$A_{\beta\gamma}\in\mathcal A_{\lambda,\Lambda}$ for
$\beta\in\mathcal B$ and $\gamma\in\mathcal G$.
Moreover, the hyperplanes
\[
\{\Pi(\beta,\gamma,\beta',\gamma'):
(\beta,\gamma)\ne(\beta',\gamma')\}
\]
are parallel in $\mathcal S$, where $\Pi(\beta,\gamma,\beta',\gamma')$
is given by
\[
\text{tr}(A_{\beta\gamma}M)-f_{\beta\gamma}
=\text{tr}(A_{\beta'\gamma'}M)-f_{\beta'\gamma'}.
\]
\end{enumerate}
\medskip

\begin{theorem}\label{thm:Theorem 1.6}
If one of \hyperlink{I3}{\bf(I3)\rm}--\hyperlink{I5}{\bf(I5)\rm}
holds, then the Isaacs operator $F$ in \eqref{eq:Isaacs operator} satisfies
\[
F(\cdot,x)\in\mathbb T^\omega(\mathcal S,\lambda',\Lambda')
\quad\text{for every }x\in\Omega,
\]
where $\lambda'$, $\Lambda'$, and $\omega$ depend only on $n$ and the
structural data in \hyperlink{I3}{\bf(I3)\rm}--\hyperlink{I5}{\bf(I5)\rm},
and $\lim_{t\to0^+}\omega(t)=0$.
\end{theorem}

The organization of the paper is as follows.  
Section \ref{sec:C2alpha partial regularity} contains the partial $C^{2,\alpha}$ regularity and 
a Liouville type property for nondifferentiable equations.   
The partial $W^{2,\text{VMO}}$ regularity is developed in Section 3.   
$W^{3,\delta}$ estimates with $W^{1,n}$ data and 
Hausdorff dimensions of singular sets of viscosity solutions 
are the contents of Section 4.  
Theorem \ref{thm:Theorem 1.4} is proved in Section \ref{sec:nondiff fully nonlinear satisfying (H2)}.  
Finally, Section \ref{sec:Isaacs equations} is devoted to the issue of smooth enveloping
and partial regularity for Isaacs equations.

\setcounter{equation}{0}

\section{Partial $C^{2,\alpha}$ regularity}\label{sec:C2alpha partial regularity}

\subsection{Viscosity solutions and an approximation lemma}
We recall several notions and results for viscosity solutions.    

Let $f\in C(\Omega)$.  The function $u\in C(\Omega)$ is a viscosity subsolution 
(supersolution) to \eqref{eq:main equation} in $\Omega$ if for each $x_0\in\Omega$ and 
$\varphi\in C^2(\Omega)$ such that $u-\varphi$ has a local maximum (minimum) 
at $x_0$ we have 
$$
F(D^2\varphi(x_0), x_0) \ge (\le) f(x_0); 
$$    
$u$ is a viscosity (or C-viscosity) solution if it is both a viscosity subsolution and 
a viscosity supersolution.  

An important result is that if $u$ is a viscosity solution of \eqref{eq:main equation} in 
$\Omega$, then $u$ is twice differentiable (in the sense of having 
second order Taylor expansions) a.e. in $\Omega$ and its pointwise derivatives 
satisfy \eqref{eq:main equation} a.e. in $\Omega$.  See \cite{CCKS} for more details.  

The Pucci extremal operators are defined as follows:
\begin{align*}
\mathcal P^+(M)&=\mathcal P^+_{\lambda, \Lambda}(M)=
   \sup_{A\in {\mathcal A}_{\lambda, \Lambda}} \text{tr\,}(AM), \\
\mathcal P^-(M)&=\mathcal P^-_{\lambda, \Lambda}(M)=
   \inf_{A\in {\mathcal A}_{\lambda, \Lambda}} \text{tr\,}(AM).
\end{align*}

We adopt the following notations
\begin{align*}
\underline S(\lambda,\Lambda,f)&=\{u\in C(\Omega): 
\mathcal P^+_{\lambda, \Lambda}(D^2u)\ge f  \text{ in the viscosity sense in } \Omega\}, 
\\
\overline S(\lambda,\Lambda,f)&=\{u\in C(\Omega): 
\mathcal P^-_{\lambda, \Lambda}(D^2u)\le f  \text{ in the viscosity sense in } \Omega\}, 
\\
S(\lambda,\Lambda,f)&=\underline S(\lambda,\Lambda,f)\cap 
\overline S(\lambda,\Lambda,f). 
\end{align*}

If $u$ is a viscosity solution to $F(D^2u, x)=f$ in $\Omega$, it follows from the  definitions
that 
$$
u\in S(\lambda/n, \Lambda, f(x)-F(0,x)). 
$$
Consequently, the $W^{2,\delta}$ estimates 
([CC1, Proposition 7.4]) yield that if $u$ is a viscosity solution to \eqref{eq:main equation} 
in $B_1$ then 
\begin{equation}\label{eq:Ldelta norm of Hessian}
\|D^2u\|_{L^\delta(B_{1/2})}
\le C(\|u\|_{L^\infty(B_1)} + \|f(x)-F(0,x)\|_{L^n(B_1)}),
\end{equation}
where $0<\delta<1$ and $C$ are universal.

For convex equations $F(D^2u)=0$, the celebrated Evans-Krylov Theorem asserts 
\begin{equation}\label{eq:Evans Krylov theorem}
\| u \|_{C^{2, \alpha_o}(\overline B_{1/2})} \le C(\|u\|_{L^\infty(B_1)} + |F(0)|), 
\end{equation}
where $0<\alpha_0<1$ and $C$ are universal constants.  
We mention that Campanato's type inequalities that are equivalent to \eqref{eq:Evans Krylov theorem} were considered 
in \cite{CH}.  

We shall prove the following approximation lemma.  

\begin{lemma}\label{lm:Lemma 2.1}
Let $u$ be a viscosity solution to $F(D^2u,x)=f(x)$ in $B_1$.    
Suppose that $F$ satisfies \eqref{eq:ellipticity condition}, $F(0,0)=f(0)=0$, and
$$
|F(M,0)-T(M)| \le \omega(\|M\|) \|M\| \qquad \text{for } M\in \overline{\mathcal B}_1, 
$$
where $T$ is a convex or concave uniformly elliptic operator with ellipticity
constants $\lambda$, $\Lambda$, and $\omega$ is nondecreasing with
$\lim_{t\to 0^+}\omega(t)=\omega(0)=0$.  

Let $\alpha_0$ be the exponent in \eqref{eq:Evans Krylov theorem} and $0<\alpha<\alpha_0$. 
Then there exist constants 
$\eta_0=\eta_0(n, \lambda, \Lambda, \alpha, \omega)>0$, and  
$0<\sigma=\sigma(n, \lambda, \Lambda, \alpha)$, 
$\mu=\mu(n, \lambda, \Lambda, \alpha)<1$, such that if 
$$
\|u\|_{L^\infty(B_1)}\le \varepsilon, \quad
\(\fint_{B_1}|f|^n\,dx\)^{1/n}\le \sigma\varepsilon, \quad
\(\fint_{B_1}[\beta_F(x)]^n\,dx\)^{1/n}\le \sigma\varepsilon, \quad 
\varepsilon\le \eta_0,   
$$
then 
\begin{equation}\label{eq:Ldelta norm small}
\(\fint_{B_{1/2}}|D^2u|^\delta\,dx\)^{1/\delta}\le C_0\varepsilon
\end{equation}
and there exists a second degree polynomial $P$ satisfying 
\begin{align}\label{eq:approx by polynomials}
F(D^2P, 0)&=0,
\\
\|u-P\|_{L^\infty(B_\mu)}&\le \mu^{2+\alpha}\varepsilon,\notag
\\
|P(0)|+|DP(0)|+\|D^2P\|&\le C_0\varepsilon,\notag
\end{align}
where $\delta$ is as in \eqref{eq:Ldelta norm of Hessian} and $C_0$ depends only on $n$, 
$\lambda$, and $\Lambda$.  
\end{lemma}

\begin{proof}
Note that $u\in S(\lambda/n, \Lambda, f(x)-F(0,x))$.  By the definition of 
$\beta_F$ and $F(0,0)=0$, $|F(0,x)|\le \beta_F(x)$ for $x\in B_1$.  It follows 
from \eqref{eq:Ldelta norm of Hessian} that 
$$
\|D^2u\|_{L^\delta(B_{1/2})}
\le C(\|u\|_{L^\infty(B_1)} + \|f\|_{L^n(B_1)} + \| \beta_F\|_{L^n(B_1)})
\le C\varepsilon. 
$$
So \eqref{eq:Ldelta norm small} is proved.  

In addition, the H\"older estimate ([CC1, Proposition 4.10]) yields that 
\begin{equation}\label{eq:holder estimate in terms of betaF}
\|u\|_{C^{\alpha_1}(\overline B_{1/2})}
\le C(\|u\|_{L^\infty(B_1)} + \|f\|_{L^n(B_1)} + \| \beta_F\|_{L^n(B_1)})
\le C\varepsilon, 
\end{equation}
where $0<\alpha_1<1$ and $C$ are universal constants.  

To prove \eqref{eq:approx by polynomials}, we use a perturbation argument.  Let $v$ be the 
classical solution to the problem
$$
\left\{
\aligned
 &T(D^2v)=0\\ 
 &v=u
\endaligned
\right.
\quad
\aligned 
&\text {in }B_{1/2},\\
&\text {on }\partial B_{1/2}.
\endaligned
$$
From \eqref{eq:Evans Krylov theorem}, $T(0)=0$, and the maximum principle, one obtains that
\begin{equation}\label{eq:consequence max principle}
\| v \|_{C^{2, \alpha_0}(\overline B_{1/4})} \le C\|v\|_{L^\infty(B_{1/2})}
\le C\|u\|_{L^\infty(B_{1/2})}.
\end{equation}
Since $u$ is $C^{\alpha_1}$ on $\partial B_{1/2}$ and  
$v\in S(\lambda/n, \Lambda, 0)$, applying the global H\"older estimate 
from [CC1, Proposition 4.13] in $B_{1/2}$ yields 
\begin{equation}\label{eq:estimate holder seminorm}
[ v ]_{C^{\alpha_1/2}(\overline B_{1/2})} 
\le C\|u\|_{C^{\alpha_1}(\partial B_{1/2})}.
\end{equation}

To estimate $u-v$, note that $v\in C^{2, \alpha_0}(B_{1/2})$.  Therefore 
$u-v\in S(\lambda/n, \Lambda, g)$ in $B_{1/2}$ with $g(x)=f(x)-F(D^2v,x)$.  
For $0<\theta<1/4$, the ABP estimate leads to 
\begin{equation}\label{eq:ABP consequence}
\|u-v\|_{L^\infty(B_{1/2-\theta})}
\le \|u-v\|_{L^\infty(\partial B_{1/2-\theta})} + C(\|f\|_{L^n(B_{1/2-\theta})} + 
\|F(D^2v, x)\|_{L^n(B_{1/2-\theta})}).   
\end{equation}
Let $y\in \partial B_{1/2-\theta}$.  Recall that $u-v=0$ on $\partial B_{1/2}$. 
By \eqref{eq:estimate holder seminorm} and \eqref{eq:holder estimate in terms of betaF}, we obtain   
\begin{equation}\label{eq:u-v close}
|(u-v)(y)| \le [u-v]_{C^{\alpha_1/2}(\overline B_{1/2})} \theta^{\alpha_1/2}
\le C\|u\|_{C^{\alpha_1}(\overline B_{1/2})} \theta^{\alpha_1/2}
\le C\varepsilon \theta^{\alpha_1/2}. 
\end{equation}

Next estimate $\|F(D^2v, x)\|_{L^n(B_{1/2-\theta})}$ in \eqref{eq:ABP consequence}.  For each 
$x\in \overline B_{1/2-\theta}$, we claim that 
\begin{equation}\label{eq:bound hessian of v}
|D^2v(x)| \le C_1\theta^{\alpha_1/2-2} \varepsilon,
\end{equation}
where $C_1$ is universal.   
To prove this, let $w(z)=v(x+\theta z)/ \theta^2$ and pick $z_1\in \partial B_1$.   We 
have $T\big(D^2(w-w(z_1))\big)=0$ in $B_1$ and so \eqref{eq:Evans Krylov theorem} yields 
$$
\|D^2w\|_{L^\infty(B_{1/2})}=
\|D^2(w-w(z_1))\|_{L^\infty(B_{1/2})}
\le C\|w-w(z_1)\|_{L^\infty(B_1)}.
$$
Therefore
$$
|D^2v(x)|=|D^2w(0)| 
\le C\theta^{\alpha_1/2-2} [v]_{C^{\alpha_1/2}(\overline B_{1/2})}. 
$$ 
So \eqref{eq:bound hessian of v} follows from \eqref{eq:estimate holder seminorm} and \eqref{eq:holder estimate in terms of betaF}.  

The constant $\theta$ will be determined later.  
By considering small $\eta_0$, one can assume that 
$\eta_0\le C_1^{-1}\theta^{2-\alpha_1/2}$.    
In view of  \eqref{eq:bound hessian of v}, 
$|D^2v(x)| \le  C_1\theta^{\alpha_1/2-2} \eta_0 \le 1$ for 
$x\in \overline B_{1/2-\theta}$.  It gives rise to     
$$
\aligned
|F(D^2v(x), x)|
&\le
|F(D^2v(x), x)-F(D^2v(x), 0)| + |F(D^2v(x), 0)-T(D^2v(x))|  \\
& \le
\beta_F(x)(\|D^2v(x)\|+1)+\omega(\|D^2v(x)\|) \|D^2v(x)\|  \\  
& \le
2\beta_F(x) + \omega(C_1\theta^{\alpha_1/2-2} \eta_0)
C_1\theta^{\alpha_1/2-2} \varepsilon. 
\endaligned
$$
Consequently
$$
\|F(D^2v, x)\|_{L^n(B_{1/2-\theta})}
\le 2\|\beta_F(x)\|_{L^n(B_{1/2-\theta})} +C_1 \omega(C_1\theta^{\alpha_1/2-2} \eta_0)
\theta^{\alpha_1/2-2} \varepsilon. 
$$
Inserting this estimate in \eqref{eq:ABP consequence} and using \eqref{eq:u-v close}, one obtains 
\begin{align}\label{eq:Linfty estimate of u-v}
\|u-v\|_{L^\infty(B_{1/4})} 
&\le 
C\theta^{\alpha_1/2} \varepsilon + C\|f\|_{L^n(B_1)} + C\|\beta_F(x)\|_{L^n(B_1)}
+C\omega(C_1\theta^{\alpha_1/2-2} \eta_0)\theta^{\alpha_1/2-2} \varepsilon       \notag  \\
&\le
C(\theta^{\alpha_1/2} + \sigma + 
\omega(C_1\theta^{\alpha_1/2-2} \eta_0)\theta^{\alpha_1/2-2})  \varepsilon. 
\end{align}

Let $\Phi$ be the second Taylor polynomial of $v$ at $0$ and $0<\rho<1/4$.  
Clearly, $\sup_{B_\rho} |v-\Phi| \le C[D^2v]_{C^{\alpha_0}(\overline B_{1/4})} 
\rho^{2+\alpha_0}$.  Applying \eqref{eq:Linfty estimate of u-v} and \eqref{eq:consequence max principle} yields 
$$
\aligned
\|u-\Phi\|_{L^\infty(B_\rho)}
&\le 
\|u-v\|_{L^\infty(B_\rho)} + \|v-\Phi\|_{L^\infty(B_\rho)}    \\
&\le
C(\theta^{\alpha_1/2} + \sigma 
+ \omega(C_1\theta^{\alpha_1/2-2} \eta_0)\theta^{\alpha_1/2-2} 
+ \rho^{2+\alpha_0})  \varepsilon. 
\endaligned
$$ 
In addition, it follows from \eqref{eq:consequence max principle} that  
$$
|\Phi(0)|+|D\Phi(0)|+\|D^2\Phi\|\le C\varepsilon. 
$$

We next claim that there is $t_0$ such that $|t_0| \le |F(D^2\Phi, 0)|/ \lambda$ 
and $F(D^2\Phi+t_0I_n, 0)=0$ with $I_n$ as the $n\times n$ identity matrix. 
Indeed, from the ellipticity if $t > |F(D^2\Phi, 0)|/ \lambda$, then 
$F(D^2\Phi+tI_n, 0)\ge F(D^2\Phi, 0) + \lambda t>0$.  And  
if $t < -|F(D^2\Phi, 0)|/ \lambda$, then 
$F(D^2\Phi+tI_n, 0)\le F(D^2\Phi, 0) + \lambda t<0$.  The claim follows by continuity.  
From the assumption on $F$, we get 
$ |F(D^2\Phi, 0)| \le \omega(\|D^2\Phi\|) \|D^2\Phi\| $.  

Set $P(x)=\Phi(x)+\frac12 t_0|x|^2$.  Obviously, $F(D^2P, 0)=0$ and 
\begin{equation}\label{eq:estimate polynomial at 0}
|P(0)|+|DP(0)|+\|D^2P\|\le C\varepsilon 
+ C\omega(C_1\theta^{\alpha_1/2-2} \eta_0)\varepsilon. 
\end{equation}
It is easy to see that 
\begin{equation}\label{eq:Linfty estimate of u-P}
\|u-P\|_{L^\infty(B_\rho)}  
\le C(\theta^{\alpha_1/2} + \sigma 
+ \omega(C_1\theta^{\alpha_1/2-2} \eta_0)\theta^{\alpha_1/2-2} 
+ \rho^{2+\alpha_0}) \varepsilon
+ C\omega(C_1\theta^{\alpha_1/2-2} \eta_0)\rho^2\varepsilon.   
\end{equation}
The constant $C$ in \eqref{eq:estimate polynomial at 0} and \eqref{eq:Linfty estimate of u-P} depends only on $n$, $\lambda$, 
and $\Lambda$.  
To finish the proof, first pick $\mu=\mu(n, \lambda, \Lambda, \alpha)<1/4$ 
such that $C\mu^{\alpha_0-\alpha} \le 1/5$,  
and set $\rho=\mu$.  With this value of $\rho$, we have 
$C\rho^{2+\alpha_0}\le \mu^{2+\alpha}/5$.  
Second, pick $\theta$, $\sigma>0$ depending only on $n$, $\lambda$, 
$\Lambda$, $\alpha$ such that $C\theta^{\alpha_1/2}\le \mu^{2+\alpha}/5$ 
and $C\sigma\le \mu^{2+\alpha}/5$.  
Finally, with this fixed value of $\theta$, pick $\eta_0$ sufficiently small 
such that $C\omega(C_1\theta^{\alpha_1/2-2} \eta_0)\theta^{\alpha_1/2-2} 
\le \mu^{2+\alpha}/5$.  Therefore, $P$ satisfies \eqref{eq:approx by polynomials}.  
The proof of the lemma is complete.        
\qed
\end{proof}

\subsection{Partial $C^{2,\alpha}$ regularity}
We are ready to derive the following form of pointwise $C^{2,\alpha}$ estimates.  

\begin{lemma}\label{lm:pointwise C2alpha at 0}
$(\text{Pointwise } C^{2,\alpha} \text{ at } 0)$   
Let $u$ be a viscosity solution to $F(D^2u,x)=f(x)$ in $B_1$.    
Suppose that $F$ satisfies \eqref{eq:ellipticity condition}, $F(0,0)=f(0)=0$, and
$F(\cdot, 0)\in \mathbb T^{\omega}(\overline{\mathcal B}_1, \lambda, \Lambda)$,  
where $\omega(t)$ is nondecreasing in $t\ge 0$, $\omega(t)>0$ for $t>0$, and
$\lim_{t\to 0^+}\omega(t)=\omega(0)=0$.  

Let $\alpha_0$ be the exponent in \eqref{eq:Evans Krylov theorem} and $0<\alpha<\alpha_0$. 
Then there exist constants 
$\eta_1=\eta_1(n, \lambda, \Lambda, \alpha, \omega)>0$ and  
$0<\sigma_1=\sigma_1(n, \lambda, \Lambda, \alpha)<1$ such that if 
$\varepsilon\le \eta_1$,  $\|u\|_{L^\infty(B_1)}\le \varepsilon$,  
$$
\(\fint_{B_r}|f|^n\,dx\)^{1/n}\le \sigma_1\varepsilon r^\alpha, \quad\text{and }
\(\fint_{B_r}[\beta_F(x)]^n\,dx\)^{1/n}\le \sigma_1\varepsilon r^\alpha 
\qquad \forall\ 0<r\le 1,    
$$
then there exists a second degree polynomial $P$ satisfying 
\begin{align}\label{eq:Linfty estimate u-P solution}
F(D^2P, 0)&=0,\notag
\\
\|u-P\|_{L^\infty(B_r)}&\le Cr^{2+\alpha}\varepsilon   \qquad\forall\ 0<r\le 1,
\\
|P(0)|+|DP(0)|+\|D^2P\|&\le C\varepsilon,\notag
\end{align}
with $C$ depending only on $n$, $\lambda$, $\Lambda$, and $\alpha$.    
\end{lemma}

\begin{proof}
Choose $\sigma_1=\sigma/2$ and $\eta_1=\min\{\eta_0, \frac{1-\mu^\alpha}{C_0}\}$ 
with $\sigma$, $\eta_0$, $\mu$, and $C_0$ being the constants given in Lemma \ref{lm:Lemma 2.1}.  
We will first construct by induction a sequence of quadratic polynomials  
$$
P_k(x)=a_k + b_k\cdot x + \frac12x^tc_kx 
$$ 
with the following properties     
\begin{align}
F(D^2P_k, 0)&=0,  	\label{eq:Pk solution}
\\
\|u-P_k\|_{L^\infty(B_{\mu^{k+1}})}&\le \mu^{(k+1)(2+\alpha)}\varepsilon,  
\label{eq:u-Pk estimate}
\\
\mu^{-k(2+\alpha)}|a_k - a_{k-1}| + \mu^{-k(1+\alpha)}|b_k - b_{k-1}|
+ \mu^{-k\alpha}\|c_k - c_{k-1}\| &\le C_0\varepsilon,  	
\label{eq:estimate coefficients}			
\\
\|D^2P_k\|&\le (1+\mu^\alpha+\cdots+\mu^{k\alpha})\, C_0\varepsilon,  
\label{eq:estimate D2Pk} 
\end{align}
for $k\ge 0$.  Here $P_{-1}\equiv 0$. 

For case $k=0$, apply Lemma \ref{lm:Lemma 2.1} to $u$ and let $P_0$ be equal to the polynomial 
obtained in Lemma \ref{lm:Lemma 2.1}. 
It follows from \eqref{eq:approx by polynomials} that $P_0$ satisfies \eqref{eq:Pk solution} -- \eqref{eq:estimate D2Pk} with 
$k=0$.  

Suppose $P_k(x)$ is constructed from the induction hypothesis.  We want to construct 
$P_{k+1}$ satisfying \eqref{eq:Pk solution} -- \eqref{eq:estimate D2Pk}.  Let  
$$
\aligned
&u_{k+1}(z)=\frac{(u-P_k)(\mu^{k+1}z)}{\mu^{(k+1)(2+\alpha)}},  \\
&F_{k+1}(M,z)=\mu^{-(k+1)\alpha} F(\mu^{(k+1)\alpha} M+D^2P_k,\, \mu^{k+1}z),  \\
&f_{k+1}(z)=\mu^{-(k+1)\alpha} f(\mu^{k+1}z).  
\endaligned
$$  
Obviously, 
$
F(\mu^{(k+1)\alpha} D^2u_{k+1}(z)+D^2P_k,\, \mu^{k+1}z)=f(\mu^{k+1}z)
$
for $|z|<1$.  
Therefore, $u_{k+1}$ is a viscosity solution to 
$$
F_{k+1}(D^2u_{k+1}(z), z)=f_{k+1}(z)  \qquad \text{in }B_1. 
$$  
Below we verify that $u_{k+1}$, $F_{k+1}$, and $f_{k+1}$ satisfy the assumptions 
of Lemma \ref{lm:Lemma 2.1}, so it can be applied to construct $P_{k+1}$. 
Clearly, $F_{k+1}$ has ellipticity constants $\lambda$ and $\Lambda$. 
From \eqref{eq:Pk solution}, $F_{k+1}(0, 0) =\mu^{-(k+1)\alpha} F(D^2P_k, 0)=0=f_{k+1}(0)$.   
Since $F(\cdot, 0)\in \mathbb T^{\omega}(\overline{\mathcal B}_1, \lambda, \Lambda)$  
and $\|D^2P_k\| \le \frac{C_0\eta_1}{1-\mu^\alpha} \le 1$ by \eqref{eq:estimate D2Pk},    
there exists $T_{D^2P_k, 0}(M)$ such that 
$$
|F(M,0)-T_{D^2P_k, 0}(M)| \le \omega(\|M-D^2P_k\|) \|M-D^2P_k\| 
$$
for all $M\in \mathcal S$ with $\|M-D^2P_k\| \le 1$,  
where $T_{D^2P_k, 0}$ is a convex or concave uniformly elliptic operator 
with ellipticity constants $\lambda$, $\Lambda$. 
Set 
$$
T_{k+1}(M)=\mu^{-(k+1)\alpha} T_{D^2P_k, 0}(\mu^{(k+1)\alpha} M+D^2P_k).
$$
For $M\in \overline{\mathcal B}_1$, we have 
$$
\aligned
|F_{k+1}&(M, 0) - T_{k+1}(M)|    \\
&=\mu^{-(k+1)\alpha} |F(\mu^{(k+1)\alpha} M+D^2P_k,\, 0)-
  T_{D^2P_k, 0}(\mu^{(k+1)\alpha} M+D^2P_k)|   			\\
&\le \omega(\mu^{(k+1)\alpha} \|M\|) \|M\| \le \omega(\|M\|) \|M\|.
\endaligned
$$
Moreover, $T_{k+1}$ has ellipticity constants $\lambda$, $\Lambda$ 
and is convex or concave.  By the induction assumption 
$$
\|u_{k+1}\|_{L^\infty(B_1)}
=\frac{\|u-P_k\|_{L^\infty(B_{\mu^{k+1}})}}{\mu^{(k+1)(2+\alpha)}}
\le \varepsilon.
$$    
The assumptions on $f$ give  
$$
\(\fint_{B_1}|f_{k+1}(z)|^n\,dz\)^{1/n}
=\mu^{-(k+1)\alpha}\(\fint_{B_{\mu^{k+1}}}|f(x)|^n\,dx\)^{1/n}
\le \sigma_1\varepsilon.  
$$
Note that 
$$
\aligned
\beta_{F_{k+1}}(z) 
&=\mu^{-(k+1)\alpha}\sup_{M\in\mathcal S} 
\frac{ |F(\mu^{(k+1)\alpha} M+D^2P_k,\, \mu^{k+1}z)-
F(\mu^{(k+1)\alpha} M+D^2P_k,\, 0)| } {\|M\|+1}  		\\
&\le \mu^{-(k+1)\alpha}\sup_{M\in\mathcal S}  
\frac{\|\mu^{(k+1)\alpha} M+D^2P_k\|+1}{\|M\|+1} \beta_F(\mu^{k+1}z)   \\
&\le \mu^{-(k+1)\alpha}\big (1+(1+\mu^\alpha+\cdots+\mu^{k\alpha})\,C_0\varepsilon\big)
\beta_F(\mu^{k+1}z)  
\endaligned
$$
from \eqref{eq:definition betaF} and \eqref{eq:estimate D2Pk}.  Consequently, 
$$
\(\fint_{B_1}[\beta_{F_{k+1}}(z)]^n\,dz\)^{1/n}
\le 2\mu^{-(k+1)\alpha} 
\(\fint_{B_1}[\beta_F(\mu^{k+1}z)]^n\,dz\)^{1/n}
\le 2\sigma_1\varepsilon
$$
from the selection of $\eta_1$ and the assumption on $\beta_F$.  
We apply Lemma \ref{lm:Lemma 2.1} with the operator $F_{k+1}$ and $u_{k+1}$ to obtain 
a quadratic polynomial 
$$
\Phi(x)=q_0+q_1\cdot x+\frac12x^tq_2x
$$
satisfying 
\begin{align}\label{eq:estimates for uk+1}
F_{k+1}(D^2\Phi, 0)&=0,\notag
\\
\|u_{k+1}-\Phi\|_{L^\infty(B_\mu)}&\le \mu^{2+\alpha}\varepsilon,
\\
|\Phi(0)|+|D\Phi(0)|+\|D^2\Phi\|&\le C_0\varepsilon.  \notag
\end{align}
From the definition of $u_{k+1}$ and \eqref{eq:estimates for uk+1} 
$$
\left\| \frac{(u-P_k)(\mu^{k+1}z)}{\mu^{(k+1)(2+\alpha)}} - \Phi(z) 
\right\|_{L^\infty(B_\mu)}\le \mu^{2+\alpha}\varepsilon
$$
implying
$$
\| u(x)-P_k(x)-\mu^{(k+1)(2+\alpha)}\Phi(x/\mu^{k+1}) \|_{L^\infty(B_{\mu^{k+2}})}
\le \mu^{(k+2)(2+\alpha)}\varepsilon. 
$$
Now we let 
$$
P_{k+1}(x)=P_k(x)+\mu^{(k+1)(2+\alpha)}\Phi(x/\mu^{k+1}), 
$$
and verify that $P_{k+1}$ is the desired polynomial satisfying 
\eqref{eq:Pk solution} -- \eqref{eq:estimate D2Pk} with $k+1$ in place of $k$.  To verify \eqref{eq:Pk solution}, 
we have $D^2P_{k+1}=D^2P_k+\mu^{(k+1)\alpha}D^2\Phi$.   Since 
$F_{k+1}(D^2\Phi, 0)=\mu^{-(k+1)\alpha}F(D^2P_{k+1}, 0)$,  it follows that 
$F(D^2P_{k+1}, 0)=0$.  Clearly, \eqref{eq:u-Pk estimate} follows from the definition of 
$P_{k+1}$.  Again from the definition of $P_{k+1}$, 
$$
a_{k+1}=a_k+\mu^{(k+1)(2+\alpha)}q_0, \quad
b_{k+1}=b_k+\mu^{(k+1)(1+\alpha)}q_1, \quad
c_{k+1}=c_k+\mu^{(k+1)\alpha}q_2. \quad
$$  
So \eqref{eq:estimate coefficients} and \eqref{eq:estimate D2Pk} for $k+1$ follow from  \eqref{eq:estimates for uk+1}. 

To finish the proof, let us construct the desired polynomial $P$ satisfying \eqref{eq:Linfty estimate u-P solution}.  
From \eqref{eq:estimate coefficients}, $\{a_k\}$, $\{b_k\}$, and $\{c_k\}$ are Cauchy sequences.  Let
$$
a=\lim_{k\to \infty} a_k, \quad\ b=\lim_{k\to \infty} b_k, \quad  c=\lim_{k\to \infty} c_k, 
$$
and 
$$
P(x)=a + b\cdot x + \frac12x^tcx. 
$$
Since $D^2P_k\longrightarrow D^2P$ and $F$ is continuous, we get $F(D^2P, 0)=0$.  
It is easy to see from \eqref{eq:estimate coefficients} that 
$$
|P(0)| \le \sum_{k=0}^\infty|a_k-a_{k-1}| \le \sum_{k=0}^\infty \mu^{k(2+\alpha)}
C_0\varepsilon = \frac{C_0\varepsilon}{1-\mu^{2+\alpha}}.  
$$
Similarly
$$
|DP(0)| \le  \frac{C_0\varepsilon}{1-\mu^{1+\alpha}},  \qquad
\|D^2P\| \le  \frac{C_0\varepsilon}{1-\mu^\alpha}.
$$ 
Also for $k\ge 1$,  \eqref{eq:u-Pk estimate} and \eqref{eq:estimate coefficients} give rise to 
$$
\aligned
\|u-P\|_{L^\infty(B_{\mu^k})}
&\le \|u-P_{k-1}\|_{L^\infty(B_{\mu^k})} + \|P_{k-1}-P\|_{L^\infty(B_{\mu^k})}    \\
&\le \mu^{k(2+\alpha)}\varepsilon + |a_{k-1}-a| + \mu^k|b_{k-1}-b|  
   + \frac12\mu^{2k}\|c_{k-1}-c\|		\\
&\le \mu^{k(2+\alpha)}\varepsilon +  \( \sum_{j=k}^\infty \mu^{j(2+\alpha)}
   + \mu^k \sum_{j=k}^\infty \mu^{j(1+\alpha)} + 
    \frac12\mu^{2k}\sum_{j=k}^\infty \mu^{j\alpha}  \)\,C_0\varepsilon 	\\ 
&\le \mu^{k(2+\alpha)}\varepsilon
\bigl[1+C_0\bigl(\frac{1}{1-\mu^{2+\alpha}}
   + \frac{1}{1-\mu^{1+\alpha}} + \frac{1}{2(1-\mu^\alpha)}\bigr)\bigr].
\endaligned 
$$ 
When $k=0$, the same inequality holds because 
$\|u-P\|_{L^\infty(B_1)} \le \|u\|_{L^\infty(B_1)}  + \|P\|_{L^\infty(B_1)} $ 
and $\|P\|_{L^\infty(B_1)}$ can be bounded using the estimates for $P(0)$, 
$DP(0)$, and $D^2P$.  Then given $\mu^{k+1}<r \le \mu^k$ ($k\ge 0$), we obtain 
$$
\|u-P\|_{L^\infty(B_r)} \le \|u-P\|_{L^\infty(B_{\mu^k})}
\le C\mu^{k(2+\alpha)}\varepsilon \le \frac{C}{\mu^{2+\alpha}}r^{2+\alpha}\varepsilon.
$$  
So we have proved Lemma \ref{lm:pointwise C2alpha at 0}. 
\end{proof}


For the convenience of reader, we include a proof of the following auxiliary lemma.  

\begin{lemma}\label{lm:2.3}  
Let $u\in C(B_2)$ and suppose that there exist constants $K>0$ and $0<\alpha<1$ 
such that for any $x_0\in \overline B_1$ there exists a polynomial $P_{x_0}$ of 
degree at most two satisfying 
\begin{equation}\label{eq:approx by polynomials bis}
|u(x)-P_{x_0}(x)| \le K|x-x_0|^{2+\alpha}, \qquad\text{for }|x-x_0|<1.  
\end{equation}
Then $u\in C^{2,\alpha}(\overline B_1)$ and there exists a dimensional constant $C>0$ such that 
\[
[D^2 u]_{C^\alpha(\bar B_1)}\leq C\,K.
\] 
Moreover for $x_0\in B_1$,  $P_{x_0}(x)$ is the second Taylor Polynomial of 
$u$ at $x_0$. 
\end{lemma}

\begin{proof}
Let $J$ be a standard mollifier and 
$J_\epsilon(x)=\epsilon^{-n}J(x/\epsilon)$.  Consider $u_\epsilon=J_\epsilon*u$.  
We want to show that $\{u_\epsilon\}$ is a Cauchy family in $C^2(\overline B_1)$
as $\epsilon \longrightarrow 0$.  Note that for $x\in \overline B_1$, $1\le i$, $j\le n$,  
$0<\epsilon<1$,  
$$
D_{ij}u_\epsilon(x)
=\frac1{\epsilon^2}\int_{\R^n}u(y)\epsilon^{-n}D_{ij}J\left((x-y)/\epsilon\right)\,dy
=\frac1{\epsilon^2}\int_{\R^n}u(x-\epsilon z)D_{ij}J(z)\,dz.  
$$
Writing $u(y)=P_x(y)+E_x(y)$ yields 
\begin{align}\label{eq:derivatives of uepsilon}
D_{ij}u_\epsilon(x)
&=\frac1{\epsilon^2}\int_{\R^n}P_x(x-\epsilon z)D_{ij}J(z)\,dz + 
     \frac1{\epsilon^2}\int_{\R^n}E_x(x-\epsilon z)D_{ij}J(z)\,dz		\\
&=D_{ij}P_x + \frac1{\epsilon^2}\int_{\R^n}E_x(x-\epsilon z)D_{ij}J(z)\,dz. \notag
\end{align} 
It follows from \eqref{eq:approx by polynomials bis} that 
\begin{equation}\label{eq:estimate Dij uepsilon}
|D_{ij}u_\epsilon(x)-D_{ij}P_x| \le C_nK\epsilon^\alpha. 
\end{equation}  
Therefore,  $\{D_{ij}u_\epsilon\}$ is a Cauchy family in $C(\overline B_1)$ as 
$\epsilon \longrightarrow 0$.  Since $u_\epsilon \longrightarrow u$ in $C(\overline B_1)$
and by the interpolation inequalities, $u_\epsilon \longrightarrow v$ in $C^2(\overline B_1)$. 
Thus, $u=v\in C^2(\overline B_1)$.  Furthermore, $D^2u(x)=D^2P_x$.     

To show that $D_{ij}u$ are $C^\alpha$, let $x$, $x_0\in \overline B_1$, $0<|x-x_0|<1/2$, 
and $\epsilon=|x-x_0|$.  Analogous to \eqref{eq:derivatives of uepsilon}, one obtains 
$$
D_{ij}u_\epsilon(x)
=D_{ij}P_{x_0} + \frac1{\epsilon^2}\int_{\R^n}E_{x_0}(x-\epsilon z)D_{ij}J(z)\,dz.
$$  
This equality and \eqref{eq:derivatives of uepsilon} with $x_0$ in place of $x$ imply 
\begin{align}\label{eq:holder Dijuepsilon}
|D_{ij}u_\epsilon(x) - D_{ij}u_\epsilon(x_0)|
&\le \frac1{\epsilon^2}\int_{\R^n} |E_{x_0}(x-\epsilon z) - E_{x_0}(x_0-\epsilon z)|
|D_{ij}J(z)|\,dz 	\notag	\\
&\le C_nK\epsilon^\alpha = C_nK|x-x_0|^\alpha. 
\end{align}
From \eqref{eq:estimate Dij uepsilon} and \eqref{eq:holder Dijuepsilon}, we get 
\begin{align}\label{eq:2:24}
|D_{ij}u&(x) - D_{ij}u(x_0)| \notag	\\
&\le |D_{ij}u(x) - D_{ij}u_\epsilon(x)| + |D_{ij}u_\epsilon(x) - D_{ij}u_\epsilon(x_0)|
   + |D_{ij}u_\epsilon(x_0) - D_{ij}u(x_0)|	\notag	\\
&\le C_nK|x-x_0|^\alpha.
\end{align}
For the case when $1/2 \le |x-x_0| < 1$, let $\overline x=(x+x_0)/2$.  Then we apply 
\eqref{eq:2:24} to $|D_{ij}u(x) - D_{ij}u(\overline{x})|$ and 
$|D_{ij}u(\overline x) - D_{ij}u(x_0)|$ 
to get $|D_{ij}u(x) - D_{ij}u(x_0)|\le C_nK|x-x_0|^\alpha$.  The lemma now follows.

\end{proof}

The following result concerns the $C^{2, \alpha}$ regularity of flat viscosity solutions 
to nondifferentiable fully nonlinear elliptic equations.  

\begin{theorem}\label{thm:2.4}
$(C^{2,\alpha} \text{ in } \overline B_1)$   
Let $u$ be a viscosity solution to $F(D^2u,x)=f(x)$ in $B_2$.    
Suppose that $F$ satisfies \eqref{eq:ellipticity condition}, $F(0,0)=f(0)=0$, and
$F(\cdot, x)\in \mathbb T^{\omega}(\overline{\mathcal B}_2, \lambda, \Lambda)$ 
for every $x\in \overline B_1$,  
where $\omega(t)$ is nondecreasing in $t\ge 0$, $\omega(t)>0$ for $t>0$, and
$\lim_{t\to 0^+}\omega(t)=\omega(0)=0$.  

Let $\alpha_o$ be as in \eqref{eq:Evans Krylov theorem} and $0<\alpha<\alpha_o$. 
Then there exist constants 
$\eta_2=\eta_2(n, \lambda, \Lambda, \alpha, \omega)>0$ and  
$0<\sigma_2=\sigma_2(n, \lambda, \Lambda, \alpha)<1$ such that if 
\begin{align*}
&\varepsilon\le \eta_2, \quad  
\|u\|_{L^\infty(B_2)}\le \varepsilon, \quad 
[f]_{C^\alpha(\overline B_2)}\le \sigma_2\varepsilon,      \\
&\beta_F(x_1, x_2)\le \sigma_2\varepsilon |x_1-x_2|^\alpha 
\qquad \forall\ x_1, x_2\in \overline B_2, 
\end{align*}   
then $u\in C^{2,\alpha}(\overline B_1)$ and 
$$
\|u\|_{C^{2,\alpha}(\overline B_1)}\le C\varepsilon
$$
with $C$ depending only on $n$, $\lambda$, $\Lambda$, and $\alpha$.    
\end{theorem}

\begin{proof}
Let $x_0\in \overline B_1$.  We claim that there exists $t_0\in \mathbb R$ such that
$F(t_0I_n, x_0)=f(x_0)$ with 
\begin{equation}\label{2.25}
|t_0|\le \frac{|f(x_0) - F(0, x_0)|}\lambda \le \frac{2\sigma_2\varepsilon}\lambda.
 \end{equation}
Indeed, if $t > |f(x_0) - F(0, x_0)|/ \lambda$, then 
$$
F(tI_n, x_0)=F(tI_n, x_0) - F(0, x_0) + F(0, x_0) \ge \lambda t + F(0, x_0) > f(x_0)
$$
by the ellipticity.  
Similarly, if $t < -|f(x_0) - F(0, x_0)|/ \lambda$, then 
$F(tI_n, x_0) < f(x_0)$.  So the existence of $t_0$ follows from the continuity of $F$.  
Noting that $F(0, 0)=0$ and $f(0)=0$, we obtain that 
$|F(0, x_0)|\le \beta_F(x_0, 0)\le \sigma_2\varepsilon$ and 
$|f(x_0)|\le [f]_{C^\alpha(\overline B_2)}|x_0|^\alpha \le \sigma_2\varepsilon$. 
Thus, \eqref{2.25} is proved.  

Let $\widetilde u(y):= u(y+x_0)-\frac12t_0|y|^2$.  Set  
\begin{align*}
&\widetilde F(M,y)=F(M+t_0I_n,\, y+x_0) - f(x_0),  \\
&\widetilde f(y)=f(y+x_0)-f(x_0).  
\end{align*} 
Clearly, 
$
F(D^2\widetilde u(y)+t_0I_n,\, y+x_0)=f(y+x_0)
$
for $|y|<1$.  
Therefore, $\widetilde u$ satisfies in the viscosity sense 
$$
\widetilde F(D^2\widetilde u(y), y)=\widetilde f(y)  \qquad \text{in }B_1. 
$$  

Choose $\sigma_2=\min\{\sigma_1, \frac{\lambda}2\}$ and 
$\eta_2=\min\{1, \frac{\eta_1}2\}$ 
with $\sigma_1$ and $\eta_1$ being the constants given in Lemma \ref{lm:pointwise C2alpha at 0}.  

Below we verify that $\widetilde u$, $\widetilde F$, and $\widetilde f$ satisfy 
the assumptions of Lemma \ref{lm:pointwise C2alpha at 0}. 
Obviously, $\widetilde F$ has ellipticity constants $\lambda$ and $\Lambda$.   
Also, $\widetilde f(0)=0$ and 
$\widetilde F(0, 0) = F( t_0I_n, x_0) - f(x_0)=0$.     
Let $M_0\in \overline{\mathcal B}_1$ and $N_0=M_0+t_0I_n$.   
Then $\|N_0\| \le 1+2\sigma_2/ \lambda \le 2$ by \eqref{2.25}.      
Since $F(\cdot, x_0)\in \mathbb T^{\omega}(\overline{\mathcal B}_2, \lambda, \Lambda)$,      
there exists $T_{N_0, x_0}(M)$ such that 
\begin{equation}\label{2.26}
|F(M, x_0) - T_{N_0, x_0}(M)| \le \omega(\|M-N_0\|) \|M-N_0\|
\end{equation}
for all $M\in \mathcal S$ with $\|M-N_0\| \le 1$,  
where $T_{N_0, x_0}$ is a convex or concave uniformly elliptic operator 
with ellipticity constants $\lambda$, $\Lambda$. 
Set 
$$
\widetilde T_{M_0, 0}(M)=T_{N_0, x_0}(M+t_0I_n)-f(x_0).
$$
For $M-M_0\in \overline{\mathcal B}_1$, from \eqref{2.26}, we obtain 
\begin{align*}
|\widetilde F(M, 0) - \widetilde T_{M_0, 0}(M)|    
&=|F(M+t_0I_n,\, x_0)-T_{N_0, x_0}(M+t_0I_n)|   		\\
&\le \omega(\|M-M_0\|) \|M-M_0\|.
\end{align*}
Obviously, $\widetilde T_{M_0, 0}$ has ellipticity constants $\lambda$, $\Lambda$ 
and is convex or concave.  As a result,  
$\widetilde F(\cdot, 0)\in \mathbb T^{\omega}(\overline{\mathcal B}_1, \lambda, \Lambda)$.  

Furthermore, $2\varepsilon \le 2\eta_2 \le \eta_1$.  
Set $\widetilde\varepsilon=2\varepsilon$.  Note $|t_0|\le \varepsilon$ 
by \eqref{2.25} and the selection of $\sigma_2$.  Thus, we get 
$\|\widetilde u\|_{L^\infty(B_1)} \le \|u\|_{L^\infty(B_2)} + |t_0| 
\le \widetilde\varepsilon$.    
It is easy to see that   
$
\sup_{B_r}|\widetilde f| \le [f]_{C^\alpha(\overline B_r(x_0))}r^\alpha
\le \sigma_1 \widetilde\varepsilon r^\alpha  
$
for $0<r\le 1$.  
For $|y| \le r \le 1$, we have  
\begin{align*}
\beta_{\widetilde F}(y) 
&=\sup_{M\in\mathcal S}
\frac{ |F(M+t_0I_n,\, y+x_0)-F(M+t_0I_n,\, x_0)| } {\|M\|+1}  		\\
&\le \sup_{M\in\mathcal S}
\frac{\|M+t_0I_n\|+1}{\|M\|+1} \beta_F(y+x_0, x_0)   \\
&\le 2 \sigma_2\varepsilon |y|^\alpha
   \le \sigma_1 \widetilde\varepsilon r^\alpha
\end{align*}
from \thetag{1.4} and the assumption on $\beta_F$.  

We apply Lemma \ref{lm:pointwise C2alpha at 0} with $\widetilde{F}$,  $\widetilde f$, and $\widetilde u$ 
to obtain a quadratic polynomial $\Phi_{x_0}(y)$ satisfying 
\begin{align*}
&|\widetilde u(y)-\Phi_{x_0}(y)| \le C\widetilde\varepsilon |y|^{2+\alpha}
\qquad\text{for }|y| \le 1, 
\\
&|\Phi_{x_0}(0)|+|D\Phi_{x_0}(0)|+\|D^2\Phi_{x_0}\|\le C\widetilde\varepsilon,  
\end{align*}
with $C=C(n, \lambda, \Lambda, \alpha)$.  
Set $P_{x_0}(x)=\frac12t_0|x-x_0|^2+\Phi_{x_0}(x-x_0)$.  Then 
$$
|u(x)-P_{x_0}(x)| \le C\widetilde\varepsilon |x-x_0|^{2+\alpha}
\qquad\text{for }|x-x_0| \le 1.  
$$
By Lemma \ref{lm:2.3}, $u\in C^{2,\alpha}(\overline B_1)$ with 
$[D^2u]_{C^\alpha(\overline B_1)}\le C\varepsilon$.  
Since $P_{x_0}$ is the second Taylor polynomial of $u$ at $x_0\in B_1$,  
the estimates on $\Phi_{x_0}$ and $t_0$ yield 
$\|u\|_{C^2(\overline B_1)} \le C\varepsilon$. 
The proof of Theorem \ref{thm:2.4} is finished.  
\end{proof}

We are ready to prove Theorem \ref{thm:Theorem 1.1}.

\begin{proof}(of Theorem \ref{thm:Theorem 1.1}) 
Set $E=\{x\in\Omega: u \text{ is twice differentiable at }x \}$.  
Note that $|\Omega - E|=0$.  
Let $\alpha_0$ be the exponent in \eqref{eq:Evans Krylov theorem}. 

Given $x_0\in E$, we want to show that  
$u\in C^{2,\alpha}(\overline B_{r_{x_0}}(x_0))$ for some small $r_{x_0}$.   
By the definition of $E$, there exists a quadratic polynomial $P$ such 
that $u(x)-P(x)=o(|x-x_0|^2)$ for $x$ near $x_0$ and $F(D^2P, x_0)=f(x_0)$. 
Let $k_0 \ge \|D^2P\|+2$ and $B_{4r}(x_0)\subset \Omega$ with small $r$ 
to be determined later.  Consider $v(y)= (u-P)(ry+x_0)/r^2$.  Set  
\begin{align*}
&G(M,y)=F(M+D^2P,\, ry+x_0) - f(x_0),  \\
&g(y)=f(ry+x_0)-f(x_0).  
\end{align*}
It is easy to see that $v$ is a viscosity solution to   
$$
G(D^2v, y)=g(y)  \qquad \text{in }B_2. 
$$  

Below we verify that $v$, $G$, and $g$ satisfy 
the assumptions of Theorem \ref{thm:2.4}.    
Obviously, $G$ has ellipticity constants $\lambda$ and $\Lambda$   
and $G(0, 0) = g(0)=0$.     Let $M_1\in \overline{\mathcal B}_2$.  
Then $N_1=M_1+D^2P\in \overline{\mathcal B}_{k_0}$.     
It follows from (H1) that for $x\in\Omega$, there exists $T_{N_1, x}(M)$ 
convex or concave with ellipticity constants $\lambda$, $\Lambda$
such that 
\begin{equation}\label{2.27}
|F(M, x) - T_{N_1, x}(M)| \le \omega_{k_0}(\|M-N_1\|) \|M-N_1\|
\end{equation}
for all $M-N_1\in \overline{\mathcal B}_1$.  
Let $|y_1| \le 1$ and $x_1=ry_1+x_0$.  
Clearly, $\widetilde T_{M_1, y_1}(M):=T_{N_1, x_1}(M+D^2P)-f(x_0)$
has ellipticity constants $\lambda$, $\Lambda$ and is convex or concave.  
For $M-M_1\in \overline{\mathcal B}_1$, from \eqref{2.27}, we have 
\begin{align*}
|G(M, y_1) - \widetilde T_{M_1, y_1}(M)|    
&= |F(M+D^2P,\, x_1)-T_{N_1, x_1}(M+D^2P)|   	\\	
&\le \omega_{k_0}(\|M-M_1\|) \|M-M_1\|. 
\end{align*}
Therefore, 
$G(\cdot, y)\in \mathbb T^{\omega_{k_0}}(\overline{\mathcal B}_2, \lambda, \Lambda)$ 
for any $y\in \overline B_1$.      

Let $\eta_2$ and $\sigma_2$ be the constants in Theorem \ref{thm:2.4} with $\omega_{k_0}$ 
in place of $\omega$.  In the following estimates on $v$, $g$, $\beta_G$, 
we will take $r$ sufficiently small.     
Since $P$ is the second Taylor polynomial of $u$ at $x_0$, we obtain
$\|v\|_{L^\infty(B_2)} \le \eta_2$.   Also
$[g]_{C^\alpha(\overline B_2)} \le r^\alpha [f]_{C^\alpha(\overline B_{2r}(x_0))}
\le \sigma_2\eta_2$.  
For $z_1$, $z_2 \in \overline B_2$, \eqref{eq:definition betaF} and \eqref{eq:holder in x} yield 
\begin{align*}
\beta_G(z_1, z_2) 
&=\sup_{M\in\mathcal S}
\frac{ |F(M+D^2P,\, rz_1+x_0)-F(M+D^2P,\, rz_2+x_0)| } {\|M\|+1}  		\\
&\le \sup_{M\in\mathcal S}
\frac{A(1+\|M+D^2P\|)r^\alpha |z_1-z_2|^\alpha}{\|M\|+1}      \\
&\le Ak_0 r^\alpha |z_1-z_2|^\alpha
   \le \sigma_2\eta_2 |z_1-z_2|^\alpha. 
\end{align*}

Applying Theorem \ref{thm:2.4} with $v$, $G$, $g$, and $\omega_{k_0}$,  
we conclude that $v\in C^{2,\alpha}(\overline B_1)$ and hence
$u\in C^{2,\alpha}(\overline B_{r_{x_0}}(x_0))$ with $r_{x_0}=r$.  

Choose $\Omega_1=\bigcup_{x_0\in E} B_{r_{x_0}}(x_0)$.  
Theorem \ref{thm:Theorem 1.1} follows.             
\end{proof}


\subsection{A Liouville type property and everywhere $C^{2,\alpha}$ 
estimates} 
The elliptic equation $F(D^2u)=0$ or $F$ is said to satisfy the Liouville 
property if any entire viscosity solution $u\in C^{1,1}_{loc}(\mathbb R^n)$ 
to $F(D^2u)=0$ with bounded $D^2u$ in $\mathbb R^n$ must be a polynomial 
of degree at most 2.  

The work \cite{H1} established $C^{2,\alpha}$ estimates and regularity for
$C^{1,1}$ solutions of the equation $F(D^2u)=0$ if $F\in C^1$ and the 
Liouville property holds for $F$.  The following theorem extends this result 
to nondifferentiable equations. 

\begin{theorem}\label{Theorem 2.5}  
Let $u\in C^{1,1}(B_1)$ be a viscosity solution to $F(D^2u)=0$ in $B_1$.    
Suppose that $F$ satisfies \eqref{eq:ellipticity condition}, $F(0)=0$, and
$F\in \mathbb T^{\omega_k}(\overline{\mathcal B}_k, \lambda, \Lambda)$ for 
every $k\ge 1$,  
where $\omega_k(t)$ is nondecreasing in $t\ge 0$, $\omega_k(t)>0$ for $t>0$, 
$\lim_{t\to 0^+}\omega_k(t)=\omega_k(0)=0$, and 
$\omega_k \le \omega_{k+1}$ for $k\ge 1$.   
Let $\alpha_0$ be the exponent in \eqref{eq:Evans Krylov theorem}.   

If $F$ satisfies the Liouville property, then $u\in C^{2,\alpha}(\overline B_{1/2})$ 
for any $0<\alpha<\alpha_0$.  
Furthermore,  $[D^2u]_{C^\alpha(\overline B_{1/2})}\le C$  
with $C$ depending only on $n$, $\lambda$, $\Lambda$, $\alpha$, $F$, 
$\{\omega_k\}$, and $\|u\|_{C^{1,1}(B_1)}$.    
\end{theorem}

\begin{proof}  
Let $|D^2u|\le K$ in $B_1$.  
Since the Liouville property holds for $F$, we apply \cite[Lemma 1]{H1} to obtain
that $D^2u\in \text{VMO}(B_{1/2})$ and 
\begin{equation}\label{2.28}
\fint_{B_R(x_0)}|D^2u-(D^2u)_{x_0,R}|\,dx \le \psi(R)
\qquad \text{for }x_0\in \overline B_{1/2}, \ R\le \frac12, 
\end{equation}  
where $(D^2u)_{x_0,R}=(D^2u)_{B_R(x_0)}$ for simplicity, 
$\lim_{R\to 0^+}\psi(R)=0$, and $\psi$ depends only on $n$, $\lambda$, 
$\Lambda$, $F$, and $K$. 

Fix $x_0\in \overline B_{1/2}$.  For $r\le 1/2$ and $p>n$, the Sobolev inequality
and \eqref{2.28} yields  
\begin{align}\label{2.29}
&\sup_{B_r(x_0)}|u(x)-u(x_0)-Du(x_0)\cdot (x-x_0)-\frac12(x-x_0)^t(D^2u)_{x_0,r}(x-x_0)|   \\
&\qquad \le Cr^2\(\fint_{B_r(x_0)}|D^2u-(D^2u)_{x_0,r}|^p\,dx\)^{1/p}
  \le Cr^2[K^{p-1}\psi(r)]^{1/p}.\notag
\end{align}
Note that $F(D^2u(x))=0$ for a.e. $x\in B_1$.  With the ellipticity of $F$ and 
\eqref{2.28}, we proceed to get     
$$
|F((D^2u)_{x_0,r})| = \fint_{B_r(x_0)}|F((D^2u)_{x_0,r})-F(D^2u(x))|\,dx
\le \Lambda \psi(r).  
$$
The ellipticity of $F$ and intermediate value theorem imply that 
there exists $t_0\in \mathbb R$ such that 
$|t_0| \le |F((D^2u)_{x_0,r})|/ \lambda$ and $F((D^2u)_{x_0,r}+t_0I_n)=0$.  

Let $P(x)=u(x_0)+Du(x_0)\cdot (x-x_0)+\frac12(x-x_0)^t[(D^2u)_{x_0,r}+t_0I_n](x-x_0)$    
and $v(y)=(u-P)(\frac{1}2ry+x_0)/(\frac{r}2)^2$ with small $r$ to be determined later.   
It is easy to check that $v$ satisfies in the viscosity sense    
$$
\widetilde F(D^2v)=0  \qquad \text{in }B_2, 
$$  
where $\widetilde F(M)=F(M+D^2P)$. 

To apply Theorem \ref{thm:2.4} to gain $C^{2,\alpha}$ regularity of $v$, we 
consider properties of $D^2P$, $\widetilde F$, and $v$.  
The estimates of $t_0$ and $F((D^2u)_{x_0,r})$ lead to 
$\|D^2P\| \le K+\Lambda\psi(\frac12)/\lambda$.  Let $k_0 \ge \|D^2P\|+2$.  We claim that 
$\widetilde F\in \mathbb T^{\omega_{k_0}}(\overline{\mathcal B}_2, \lambda, \Lambda)$.     
Indeed, let $M_1\in \overline{\mathcal B}_2$.  
Then $N_1=M_1+D^2P\in \overline{\mathcal B}_{k_0}$.  
From the assumptions on $F$, there exists $T_{N_1}(M)$ 
convex or concave with ellipticity constants $\lambda$, $\Lambda$ such that 
\begin{equation}\label{2.30}
|F(M) - T_{N_1}(M)| \le \omega_{k_0}(\|M-N_1\|) \|M-N_1\|
\end{equation}
for all $M-N_1\in \overline{\mathcal B}_1$. 
Obviously, $\widetilde T_{M_1}(M):=T_{N_1}(M+D^2P)$
has ellipticity constants $\lambda$, $\Lambda$ and is convex or concave.  
Furthermore, \eqref{2.30} yields 
$$
|\widetilde F(M) - \widetilde T_{M_1}(M)|    
\le \omega_{k_0}(\|M-M_1\|) \|M-M_1\| 
$$
for any $M-M_1\in \overline{\mathcal B}_1$.  Thus, the claim is verified.  

It is also easy to see that $\widetilde F(0)=0$ and $\widetilde F$ 
has ellipticity constants $\lambda$ and $\Lambda$.   
 
Let $\eta_2$ be the constant in Theorem \ref{thm:2.4} with $\omega_{k_0}$ 
in place of $\omega$.  In light of \eqref{2.29} and the estimate of $t_0$, we obtain 
$$
\|v\|_{L^\infty(B_2)} =\sup_{x\in B_r(x_0)}\frac{4|(u-P)(x)|}{r^2}
\le C[K^{p-1}\psi(r)]^{1/p}+C\psi(r).
$$
For sufficiently small $r$, we have $\|v\|_{L^\infty(B_2)} \le \eta_2$.   

Applying Theorem \ref{thm:2.4} with $v$, $\widetilde F$, and $\omega_{k_0}$,  
we get that $[D^2v]_{C^{\alpha}(\overline B_1)}\le C\eta_2$.  Equivalently, 
$[D^2u]_{C^{\alpha}(\overline B_{r/2}(x_0))}\le C\eta_2r^{-\alpha}$.  
Therefore, $[D^2u]_{C^{\alpha}(\overline B_{1/2})}\le C(\eta_2+K)r^{-\alpha}$.  
Theorem \ref{Theorem 2.5} follows.             
\end{proof}

\setcounter{equation}{0}
\section{Partial $W^{2,\text{VMO}}$ regularity}

We first establish the following mean oscillation estimates at $0$ for 
the pointwise second derivatives of viscosity solutions.   

\begin{lemma}\label{Lemma 3.1} 
Let $u$ be a viscosity solution to $F(D^2u,x)=f(x)$ in $B_1$.    
Suppose that $F$ satisfies \eqref{eq:ellipticity condition}, $F(0,0)=f(0)=0$, and
$F(\cdot, 0)\in \mathbb T^{\omega}(\mathcal S, \lambda, \Lambda)$,  
where $\omega(t)$ is nondecreasing in $t\ge 0$, $\omega(t)>0$ for $t>0$, and
$\lim_{t\to 0^+}\omega(t)=\omega(0)=0$.  
Let $\alpha_0$ be the exponent in \eqref{eq:Evans Krylov theorem} and $0<\alpha<\alpha_0$.   

Then there exist constants 
$\eta_3=\eta_3(n, \lambda, \Lambda, \alpha, \omega)>0$, and  
$0<\sigma_3=\sigma_3(n, \lambda, \Lambda, \alpha)$, 
$\mu=\mu(n, \lambda, \Lambda, \alpha)<1$, such that if   
\begin{align*}
\mu^\alpha\varepsilon_k&\le \varepsilon_{k+1}\le \varepsilon_k\le \eta_3
\qquad \forall\ k\ge 0,
\\  
\|u\|_{L^\infty(B_1)}&\le \varepsilon_0, 
\\  
\(\fint_{B_{\mu^k}}|f|^n\,dx\)^{1/n}&\le \sigma_3\varepsilon_k, \quad
\(\fint_{B_{\mu^k}}[\beta_F(x)]^n\,dx\)^{1/n}\le 
\frac{\sigma_3\varepsilon_k}{k+1} \qquad \forall\ k\ge 0,    
\end{align*}   
then  
\begin{equation}\label{3.1}
\(\inf_{ M\in \mathcal S } \fint_{B_{\mu^k/2}}|D^2u-M|^\delta\,dx\)^{1/\delta}
\le C\varepsilon_k  \qquad \forall\ k\ge 0, 
\end{equation}
where $\delta$ is as in \eqref{eq:Ldelta norm of Hessian} and $C$ depends only on $n$, $\lambda$, and 
$\Lambda$.    
\end{lemma}

\begin{proof}
Let $\eta_0$, $\sigma$, $\mu$, and $C_0$ be the constants given in Lemma \ref{lm:Lemma 2.1}. 
Choose $\sigma_3=\sigma$ and $\eta_3=\min\{\eta_0, \frac1{C_0}\}$.  
To prove the lemma, it suffices to show that there is a sequence $\{P_k\}$ of 
quadratic polynomials such that      
\begin{align}
F(D^2P_k, 0)&=0,  	\label{3.2}
\\
\|u-P_k\|_{L^\infty(B_{\mu^{k+1}})}&\le \mu^{2(k+1)}\varepsilon_{k+1},  \label{3.3}  
\\
\|D^2P_k\|&\le C_0(\varepsilon_0+\varepsilon_1+\cdots+\varepsilon_k),  \label{3.4}
\\
\(\fint_{B_{\mu^k/2}}|D^2u-D^2P_{k-1}|^\delta\,dx\)^{1/\delta}
&\le C_0\varepsilon_k,  \label{3.5}
\end{align}
for $k\ge 0$.  Here $P_{-1}\equiv 0$. 

We proceed to construct $\{P_k\}$ by induction. 
For case $k=0$, apply Lemma \ref{lm:Lemma 2.1} to $u$ and let $P_0$ be equal to the polynomial 
obtained in Lemma \ref{lm:Lemma 2.1}. 
It follows from \eqref{eq:Ldelta norm small} and \eqref{eq:approx by polynomials} that \eqref{3.2}, \eqref{3.4}, and 
\eqref{3.5} with $k=0$ hold.  From \eqref{eq:approx by polynomials} and the assumptions, we have
$$
\|u-P_0\|_{L^\infty(B_\mu)}\le \mu^{2+\alpha}\varepsilon_0\le \mu^2\varepsilon_1. 
$$  

Suppose $P_k$ is constructed from the induction hypothesis.  We want to construct 
$P_{k+1}$ such that \eqref{3.2} -- \eqref{3.5} hold.  Let  
\begin{align*}
u_{k+1}(z)&=\frac{(u-P_k)(\mu^{k+1}z)}{\mu^{2(k+1)}},  \\
F_{k+1}(M,z)&=F( M+D^2P_k,\, \mu^{k+1}z),  \\
f_{k+1}(z)&=f(\mu^{k+1}z).  
\end{align*} 
Obviously, $u_{k+1}$ satisfies in the viscosity sense  
$$
F_{k+1}(D^2u_{k+1}(z), z)=f_{k+1}(z)  \qquad \text{in }B_1. 
$$  
Below we verify that the assumptions 
of Lemma \ref{lm:Lemma 2.1} hold for $u_{k+1}$, $F_{k+1}$, and $f_{k+1}$.    
Trivially, $F_{k+1}$ is elliptic with constants $\lambda$ and $\Lambda$. 
From \eqref{3.2}, $F_{k+1}(0, 0) =F(D^2P_k, 0)=0=f_{k+1}(0)$.  
Note that $F(\cdot, 0)\in \mathbb T^{\omega}(\mathcal S, \lambda, \Lambda)$.       
Let $T_{D^2P_k, 0}(M)$ be a convex or concave elliptic operator 
with constants $\lambda$, $\Lambda$ such that 
$$
|F(M,0)-T_{D^2P_k, 0}(M)| \le \omega(\|M-D^2P_k\|) \|M-D^2P_k\| 
$$
for all $M\in \mathcal S$ with $\|M-D^2P_k\| \le 1$.   
Set 
$T_{k+1}(M):=T_{D^2P_k, 0}( M+D^2P_k)$.  
Then $T_{k+1}$ is convex or concave with ellipticity constants $\lambda$, $\Lambda$,  
and satisfies  
$$
|F_{k+1}(M, 0) - T_{k+1}(M)| \le \omega(\|M\|) \|M\| \qquad\forall \ M\in \overline{\mathcal B}_1.
$$
By the induction assumption and assumptions on $f$
\begin{align*}
\|u_{k+1}\|_{L^\infty(B_1)}
&=\frac{\|u-P_k\|_{L^\infty(B_{\mu^{k+1}})}}{\mu^{2(k+1)}}\le \varepsilon_{k+1},
\\
\(\fint_{B_1}|f_{k+1}(z)|^n\,dz\)^{1/n}
&=\(\fint_{B_{\mu^{k+1}}}|f(x)|^n\,dx\)^{1/n}\le \sigma\varepsilon_{k+1}. 
\end{align*} 
From \eqref{eq:definition betaF} and \eqref{3.4}, we have
\begin{align*}
\beta_{F_{k+1}}(z) 
&=\sup_{M\in\mathcal S}
\frac{|F(M+D^2P_k,\, \mu^{k+1}z)-F(M+D^2P_k,\, 0)|} {\|M\|+1}  	\\
&\le \sup_{M\in\mathcal S}
\frac{\|M+D^2P_k\|+1}{\|M\|+1} \beta_F(\mu^{k+1}z)   
\le \( 1+C_0\sum_{j=0}^k \varepsilon_j \)\beta_F(\mu^{k+1}z). 
\end{align*}
With the selection of $\eta_3$ and the assumptions on $\beta_F$, it yields  
$$
\(\fint_{B_1}[\beta_{F_{k+1}}(z)]^n\,dz\)^{1/n}
\le  
\(1+C_0(k+1)\eta_3\)\frac{\sigma_3\varepsilon_{k+1}}{k+2}
\le \sigma\varepsilon_{k+1}.
$$    

We apply Lemma \ref{lm:Lemma 2.1} with the operator $F_{k+1}$ and $u_{k+1}$ to assert that
\begin{equation}\label{3.6}
\(\fint_{B_{1/2}}|D^2u_{k+1}(z)|^\delta\,dz\)^{1/\delta}
\le C_0\varepsilon_{k+1}  
\end{equation}
and obtain a quadratic polynomial $P(x)$ satisfying 
\begin{align}\label{3.7}
F_{k+1}(D^2P, 0)&=0, \notag  
\\
\|u_{k+1}-P\|_{L^\infty(B_\mu)}&\le \mu^{2+\alpha}\varepsilon_{k+1},
\\
|P(0)|+|DP(0)|+\|D^2P\|&\le C_0\varepsilon_{k+1}.\notag      
\end{align}
It is easy to see from \eqref{3.6} that 
\[
\(\fint_{B_{\mu^{k+1}/2}}|D^2u-D^2P_k|^\delta\,dx\)^{1/\delta}
=\(\fint_{B_{1/2}}|D^2u_{k+1}(z)|^\delta\,dz\)^{1/\delta}
\le C_0\varepsilon_{k+1}. 
\]
By the definition of $u_{k+1}$ and \eqref{3.7} 
\begin{equation}\label{3.8}
\left\| \frac{(u-P_k)(\mu^{k+1}z)}{\mu^{2(k+1)}} - P(z) 
\right\|_{L^\infty(B_\mu)}\le \mu^{2+\alpha}\varepsilon_{k+1}
\le \mu^2\varepsilon_{k+2}.
\end{equation}
Then $P_{k+1}(x)=P_k(x)+\mu^{2(k+1)}P(x/\mu^{k+1})$ is the desired polynomial. 
Indeed, \eqref{3.8} leads to
$$
\| u(x)-P_{k+1}(x) \|_{L^\infty(B_{\mu^{k+2}})} \le \mu^{2(k+2)}\varepsilon_{k+2}. 
$$
Note that $D^2P_{k+1}=D^2P_k+D^2P$.   Since 
$F_{k+1}(D^2P, 0)=F(D^2P_{k+1}, 0)$, it follows that 
$F(D^2P_{k+1}, 0)=0$.  It is also easy to obtain \eqref{3.4} with $k+1$ in place of $k$.
Thus, we have proved Lemma \ref{Lemma 3.1}. 
\end{proof}

Hereafter we consider only cubes of edges parallel to the coordinate axes. 
Let $Q_0$ be a cube.    For $g\in L^p(Q_0)$ ($p>0$), define 
$$
\|g\|_{\text{BMO}^p(Q_0)}=\sup_{Q\subset Q_0}
\inf_{c\in \mathbb R}\(\fint_{Q}|g(x)-c|^p\,dx\)^{1/p},
$$  
where the supremum is taken over all cubes $Q$ in $Q_0$.  

The following lemma is a corollary of the more general John-Str\"omberg inequality
(see \cite{St}).  Nevertheless, it can be derived directly from the John-Nirenberg 
inequality. 

\begin{lemma}\label{Lemma 3.2} 
Let $g\in L^p(Q_0)$ for some $0<p<1$.  Suppose $\|g\|_{\text{BMO}^p(Q_0)}<\infty$.      
Then $g\in \text{BMO}(Q_0)$ and 
$\|g\|_{\text{BMO}^1(Q_0)} \le C\|g\|_{\text{BMO}^p(Q_0)}$ with $C$ depending
only on $n$ and $p$. 
\end{lemma}

\begin{proof}
Let $c_1\in \mathbb R$ and $h(x)=|g(x)-c_1|^p$.   For $Q\subset Q^* \subset Q_0$, 
we have
$$
\inf_{c\in \mathbb R}\fint_{Q}|h-c|\,dx
\le \inf_{c\in \mathbb R}\fint_{Q}|h-|c|^p|\,dx
\le \inf_{c\in \mathbb R}\fint_{Q}|g-c_1-c|^p\,dx.
$$
Therefore, $\|h\|_{\text{BMO}^1(Q^*)} \le \|g\|_{\text{BMO}^p(Q^*)}^p$ 
for any cube $Q^* \subset Q_0$.  The John-Nirenberg inequality gives rise to     
$$ 
\(\fint_{Q^*}|h-h_{Q^*}|^{1/p}\,dx\)^p
\le C \|h\|_{\text{BMO}^1(Q^*)} 
\le C\|g\|_{\text{BMO}^p(Q^*)}^p. 
$$
It yields 
$$
\(\fint_{Q^*}|h|^{1/p}\,dx\)^p
\le C\|g\|_{\text{BMO}^p(Q^*)}^p + \fint_{Q^*}|h|\,dx.
$$
From the definition of $h$, we have 
$$
\fint_{Q^*}|g-c_1|\,dx
\le C\|g\|_{\text{BMO}^p(Q_0)} + C\(\fint_{Q^*}|g-c_1|^p\,dx\)^{1/p}
$$ 
for $c_1\in \mathbb R$, $Q^*\subset Q_0$.  
Taking infimum over $c_1$ and supremum over $Q^*$, we obtain 
$\|g\|_{\text{BMO}^1(Q_0)} \le C\|g\|_{\text{BMO}^p(Q_0)}$.  
The lemma now follows.        
\end{proof}

The next lemma concerns the connection between pointwise derivatives 
and weak derivatives for viscosity solutions. 

\begin{lemma}\label{Lemma 3.3} 
Let $u$ be a viscosity solution to $F(D^2u,x)=f(x)$ in $\Omega$.  
Suppose that $F$ satisfies \eqref{eq:ellipticity condition}.   
If the pointwise Hessian $D^2u$ is in $L^p_{loc}(\Omega)$ with $p\ge n$, 
then $u$ belongs to $W^{2,p}_{loc}(\Omega)$.  
\end{lemma}

\begin{proof}
In virtue of \cite{CCKS}, we obtain that the equation $F(D^2u,x)=f(x)$ is satisfied 
by $u$ both in the $L^p$-viscosity sense and for almost every $x$ in $\Omega$. 
Set $g(x)=\Delta u\in L^p_{loc}(\Omega)$. 
We claim that $u$ is an $L^p$-viscosity solution to 
\begin{equation}\label{3.9}
\Delta u=g(x) \qquad\text{in } \Omega.
\end{equation}        
To see, let $x_0\in\Omega$ and $\varphi\in W^{2,p}_{loc}(\Omega)$ such that 
$u-\varphi$ attains a strict local maximum at $x_0$.   Observe that $w=u-\varphi$ 
satisfies in the $L^p$-viscosity sense 
$$
\mathcal P^+_{\lambda/n, \Lambda}(D^2w) 
 \ge F(D^2w+D^2\varphi,x)-F(D^2\varphi,x)
 =f(x)-F(D^2\varphi, x) \qquad\text{in } \Omega.
$$
Indeed, if $x_1\in\Omega$, $\zeta\in W^{2,p}_{loc}(\Omega)$, and 
$w-\zeta=u-(\zeta+\varphi)$ has a local maximum at $x_1$, then 
\begin{align*}
ess &\limsup_{x\to x_1} [\mathcal P^+_{\lambda/n, \Lambda}(D^2\zeta)-(f(x)-F(D^2\varphi, x))]
\\
&\ge ess \limsup_{x\to x_1} [F(D^2\zeta+D^2\varphi, x)-f(x)]\ge 0.
\end{align*} 
For sufficiently small $r>0$, applying the maximum principle \cite[Proposition 3.3]{CCKS}  
for $p=n$, and \cite[Lemma 2.1]{H2} for $p>n$, it follows that the upper
contact set $\Gamma^+$ of $(w-\sup_{\partial B_r(x_0)}w)^+$ in $B_r(x_0)$
has positive measure.  
Since $D^2w=D^2u-D^2\varphi\le 0$ a.e. in
$\Gamma^+$, we get that $\Delta \varphi \ge \Delta u=g$ a.e. in $\Gamma^+$.  
It leads to  
$$
ess \limsup_{x\to x_0} [\Delta\varphi -g(x)]\ge 0. 
$$
By definition, $u$ is an $L^p$-viscosity subsolution of \eqref{3.9}.
A similar argument shows that $u$ is an $L^p$-viscosity supersolution of \eqref{3.9}.
Consequently, the claim holds.

Consider an arbitrary closed ball $\overline B\subset \Omega$.  
From the $W^{2,p}$ regularity theory of linear elliptic equations, there exists 
$v\in W^{2,p}_{loc}(B)\cap C(\overline B)$ solving the Dirichlet problem 
$\Delta v=g$ in $B$ and $v=u$ on $\partial B$.  By the uniqueness of solution, 
we have $u=v$ in $B$.  Therefore, $u\in W^{2,p}_{loc}(\Omega)$.  
This completes the proof of the lemma.         
\end{proof}

The following theorem gives the $W^{2,\text{VMO}}$ regularity of flat viscosity solutions.   

\begin{theorem}\label{Theorem 3.4}
$(W^{2,\text{VMO}} \text{ in } B_1)$   
Let $u$ be a viscosity solution to $F(D^2u,x)=f(x)$ in $B_3$.    
Suppose that $F$ satisfies \eqref{eq:ellipticity condition}, $F(0,0)=f(0)=0$, and
$F(\cdot, x)\in \mathbb T^{\omega}(\mathcal S, \lambda, \Lambda)$ 
for every $x\in \overline B_2$,  
where $\omega(t)$ is nondecreasing in $t\ge 0$, $\omega(t)>0$ for $t>0$, and
$\lim_{t\to 0^+}\omega(t)=\omega(0)=0$.  
Let $\xi$ be a positive nondecreasing 
function on $(0,\, 1]$ with $\lim_{t\to 0^+}\xi(t)=0$. 

Then there exist constants 
$\eta_4=\eta_4(n, \lambda, \Lambda, \omega)>0$ and  
$0<\sigma_4=\sigma_4(n, \lambda, \Lambda)<1$ such that if 
\begin{align*}
\varepsilon\le \eta_4, \quad  
\|u\|_{L^\infty(B_3)}\le \varepsilon, \quad 
\text{osc}_{B_3}F(0, \cdot)\le \sigma_4\varepsilon,  \quad
\xi(1)=1,    \\
\text{osc}_{\overline B_r(x_0)}f\le \sigma_4\varepsilon\xi(r),  \quad    
\(\fint_{B_r(x_0)}[\beta_F(x, x_0)]^n\,dx\)^{1/n}
\le \frac{\sigma_4\varepsilon\xi(r)}{\log(2/r)}
\end{align*}   
for $0< r\le 1$, $x_0\in \overline B_2$, 
then $u\in W^{2,\text{VMO}}(B_1)$ and 
\begin{equation}\label{3.10}
\fint_{B_R(x_1)}|D^2u-(D^2u)_{x_1,R}|\,dx \le \varepsilon\xi^*(R) 
\qquad \text{for }x_1\in \overline B_1, \ R\le \frac1{2\sqrt{n}},
\end{equation}
where $\xi^*$ is nondecreasing with $\lim_{t\to 0^+}\xi^*(t)=0$ and 
depends only on $n$, $\lambda$, $\Lambda$, and $\xi$.    
\end{theorem}

\begin{proof}
Let $x_0\in \overline B_2$.  Notice that 
$|F(0, x_0)|\le \text{osc}_{B_3}F(0, \cdot)\le \sigma_4\varepsilon$ and 
$|f(x_0)|\le 2\sigma_4\varepsilon$.  As in the proof of Theorem \ref{thm:2.4}, 
there exists $t_0\in \mathbb R$ such that $F(t_0I_n, x_0)=f(x_0)$ and 
\begin{equation}\label{3.11}
|t_0|\le \frac{|f(x_0) - F(0, x_0)|}\lambda \le \frac{3\sigma_4\varepsilon}\lambda.
\end{equation}

It is easy to see that $\widetilde u(y):= u(y+x_0)-\frac12t_0|y|^2$ satisfies in 
the viscosity sense 
$$
\widetilde F(D^2\widetilde u(y), y)=\widetilde f(y)  \qquad \text{in }B_1, 
$$ 
where
\begin{align*}
\widetilde F(M,y)&=F(M+t_0I_n,\, y+x_0) - f(x_0),  \\
\widetilde f(y)&=f(y+x_0)-f(x_0).  
\end{align*} 

Let $\eta_3$, $\sigma_3$, and $\mu$ be the constants given in Lemma \ref{Lemma 3.1}
with $\alpha=\alpha_0/2$.    
Choose $\sigma_4=\min\{\sigma_3\log\min\{2, \mu^{-1}\},\, \frac{\lambda}3\}$ 
and $\eta_4=\min\{1, \frac{\eta_3}2\}$.   

Set $\varepsilon_0^*=\xi(1)=1$ and define 
$\varepsilon_k^*=\max\{\xi(\mu^k),\, \mu^\alpha\varepsilon_{k-1}^*\}$ for $k\ge 1$. 
Also let $\varepsilon_k=2\varepsilon\varepsilon_k^*$ for $k\ge 0$.  We observe that   
\begin{align} 
&\mu^\alpha\varepsilon_k\le \varepsilon_{k+1}\le \varepsilon_k\le 2\varepsilon
\le \eta_3 \qquad\forall\ k\ge 0, 
\label{3.12}	
\\
&\lim_{k\to\infty}\varepsilon_k^*=0.
\label{3.13} 
\end{align}
Indeed, the verification of \eqref{3.12} is straightforward from the definition of 
$\varepsilon_k$ and induction.   
If $\varepsilon_\infty^*=\lim_{k\to\infty}\varepsilon_k^*>0$, we obtain 
$\varepsilon_\infty^*=\max\{0, \mu^\alpha\varepsilon_\infty^*\}
<\varepsilon_\infty^*$ from the definition 
of $\varepsilon_k^*$.  So \eqref{3.13} holds true.   

To apply Lemma \ref{Lemma 3.1}, 
we should verify that $\widetilde u$, $\widetilde F$, and $\widetilde f$ 
satisfy the assumptions of that lemma.   
Clearly, $\widetilde F$ is elliptic with constants $\lambda$, $\Lambda$, and    
$\widetilde F(0, 0) = \widetilde f(0)=0$.     
Let $M_0\in \mathcal S$ and $N_0=M_0+t_0I_n$.        
Then there exists $T_{N_0, x_0}(M)$ such that 
$$
|F(M, x_0) - T_{N_0, x_0}(M)| \le \omega(\|M-N_0\|) \|M-N_0\| 
$$
for $M-N_0\in \overline{\mathcal B}_1$, 
where $T_{N_0, x_0}$ is a convex or concave uniformly elliptic operator 
with ellipticity constants $\lambda$, $\Lambda$.  It follows that 
$$
|\widetilde F(M, 0) - \widetilde T_{M_0, 0}(M)| \le \omega(\|M-M_0\|) \|M-M_0\|
\qquad\text{for }M-M_0\in \overline{\mathcal B}_1.  
$$
Here $\widetilde T_{M_0, 0}(M)=T_{N_0, x_0}(M+t_0I_n)-f(x_0)$. 
Therefore,  
$\widetilde F(\cdot, 0)\in \mathbb T^{\omega}(\mathcal S, \lambda, \Lambda)$.  
 
From \eqref{3.11} and the choice of $\sigma_4$, we have 
$\|\widetilde u\|_{L^\infty(B_1)} \le \|u\|_{L^\infty(B_3)} + |t_0| 
\le 2\varepsilon=\varepsilon_0$.    
It is easy to see that   
$$ 
\(\fint_{B_{\mu^k}}|\widetilde f|^n\,dy\)^{1/n} 
\le \text{osc}_{\overline B_{\mu^k}(x_0)}f
\le \sigma_4\varepsilon\xi(\mu^k)
\le \sigma_3\varepsilon_k
$$
for $k\ge 0$.  
As in the proof of Theorem \ref{thm:2.4}, from \eqref{eq:definition betaF}, we get  
$$
\beta_{\widetilde F}(y) 
\le \sup_{M\in\mathcal S}
\frac{\|M+t_0I_n\|+1}{\|M\|+1} \beta_F(y+x_0, x_0)   
\le 2\beta_F(y+x_0, x_0). 
$$
Therefore for $k\ge 0$ 
$$
\(\fint_{B_{\mu^k}}[\beta_{\widetilde F}(y)]^n\,dy\)^{1/n}
\le 2\(\fint_{B_{\mu^k}(x_0)}[\beta_F(x, x_0)]^n\,dx\)^{1/n}
\le \frac{2\sigma_4\varepsilon\xi(\mu^k)}{\log(2/\mu^k)}
\le \frac{\sigma_3\varepsilon_k}{k+1}. 
$$

Applying Lemma \ref{Lemma 3.1} with $\widetilde{F}$,  $\widetilde f$, and $\widetilde u$,  
we obtain 
\begin{equation}\label{3.14}
\inf_{M\in \mathcal S} \(\fint_{B_{\mu^k/2}(x_0)}|D^2u-M|^\delta\,dx\)^{1/\delta}
=
\inf_{M\in \mathcal S} \(\fint_{B_{\mu^k/2}}
|D^2\widetilde u-M|^\delta\,dy\)^{1/\delta}
\le C\varepsilon_k  
\end{equation}
for $k\ge 0$. 
Denote by $Q_r(x_0)$ the cube centered at $x_0$ with edgelength $2r$.  Let $k_0\ge 0$.  
By a standard argument regarding mean oscillation, it follows from \eqref{3.14} that 
$$
\inf_{M\in \mathcal S} \(\fint_{Q_r(x_0)}|D^2u-M|^\delta\,dx\)^{1/\delta}
\le C\varepsilon_{k_0}
\qquad\text{for }x_0\in \overline B_2, \ 0<r\le \frac{\mu^{k_0}}{2\sqrt n}.   
$$    
This implies that 
$\|D^2u\|_{\text{BMO}^\delta(Q_{\mu^{k_0}/(2\sqrt n)}(x_1))}\le C\varepsilon_{k_0}$ 
for any $Q_{\mu^{k_0}/(2\sqrt n)}(x_1) \subset \overline B_2$.  By Lemma \ref{Lemma 3.2}, we have
$$
\inf_{M\in \mathcal S} \fint_{Q_{\mu^{k_0}/(2\sqrt n)}(x_1)}|D^2u-M|\,dx
\le C\varepsilon_{k_0}=C\varepsilon\varepsilon_{k_0}^*
\qquad\text{for }x_1\in \overline B_1, \ k_0\ge 0.  
$$
Thus, \eqref{3.10} follows and $D^2u\in \text{VMO}(B_1)$.  

According to Lemma \ref{Lemma 3.3}, $u\in W^{2, \text{VMO}}(B_1)$. 
So Theorem \ref{Theorem 3.4} is proved.  
\end{proof}

The proof of Theorem \ref{thm:Theorem 1.2} is a modification of the proof of Theorem \ref{thm:Theorem 1.1}. 
We give a sketch of it.  

\begin{proof}(of Theorem \ref{thm:Theorem 1.2})  
Let $E=\{x\in\Omega: u \text{ is twice differentiable at }x \}$.  
Given $x_0\in E$, let $P$ be the quadratic polynomial such that 
$u(x)-P(x)=o(|x-x_0|^2)$ for $x$ near $x_0$ and $F(D^2P, x_0)=f(x_0)$. 
We want to show that  
$u$ is $W^{2,\text{VMO}}$ in some neighborhood of $x_0$. 
Consider $B_{4r}(x_0)\subset \Omega$.  
As in the proof of Theorem \ref{thm:Theorem 1.1}, $v(y)= (u-P)(ry+x_0)/r^2$ is a viscosity solution to   
$$
G(D^2v, y)=g(y)  \qquad \text{in }B_3, 
$$ 
where 
\begin{align*}
G(M,y)&=F(M+D^2P,\, ry+x_0) - f(x_0),  \\
g(y)&=f(ry+x_0)-f(x_0).  
\end{align*}

We next verify that the assumptions of Theorem \ref{Theorem 3.4} hold for $v$, $G$, and $g$.       
Clearly, $G$ is elliptic with constants $\lambda$, $\Lambda$
and $G(0, 0) = g(0)=0$.  Let $M_1\in \mathcal S$ and $y_1\in \overline B_2$.  
Then set $N_1=M_1+D^2P$ and $x_1=ry_1+x_0$.    
From the assumption (H2), there exists $T_{N_1, x_1}$ that is  
convex or concave with ellipticity constants $\lambda$, $\Lambda$
and tangent to $F(\cdot, x_1)$ at $N_1$ with modulus $\omega$ 
in the sense of Definition \ref{def:tangent property operator}.  
Clearly, $\widetilde T_{M_1, y_1}(M):=T_{N_1, x_1}(M+D^2P)-f(x_0)$ is also 
tangent to $G(M, y_1)$ at $M_1$.  Therefore, 
$G(\cdot, y)\in \mathbb T^\omega(\mathcal S, \lambda, \Lambda)$ 
for any $y\in \overline B_2$.    

We will take $r$ sufficiently small
in the following estimates on $v$, $G(0, \cdot)$, $g$, and $\beta_G$. 
Thus, we have $\|v\|_{L^\infty(B_3)} \le \eta_4$ and 
$\text{osc}_{B_3}G(0, \cdot)=\text{osc}_{B_{3r}(x_0)}F(D^2P, \cdot)
\le \sigma_4\eta_4$ with $\eta_4$ and $\sigma_4$ given in Theorem \ref{Theorem 3.4}.    
For $|y_0|\le 2$, $0<\rho\le 1$, it is easy to see 
$$
\text{osc}_{\overline B_\rho(y_0)}g\le 
\sup_{x_1, x_2\in \overline B_{3r}(x_0) \\ |x_1-x_2|\le 2r\rho} |f(x_1)-f(x_2)|
=\xi_1(2r\rho)
=\sigma_4\eta_4\xi_2(\rho)
$$
with $\xi_1$ being the modulus of continuity of $f$ over $\overline B_{3r}(x_0)$.   
Note that $\beta_G(y, y_0)\le (1+\|D^2P\|)\beta_F(ry+x_0,\, ry_0+x_0)$. 
It follows that 
\begin{align*}
\(\fint_{B_\rho(y_0)}[\beta_G(y, y_0)]^n\,dy\)^{1/n}
&\le (1+\|D^2P\|)\(\fint_{B_{r\rho}(ry_0+x_0)}[\beta_F(x,\, ry_0+x_0)]^n\,dx\)^{1/n}   \\
&\le (1+\|D^2P\|)\frac{\iota(r\rho)}{\log(2/\rho)}      
=\frac{\sigma_4\eta_4\xi_3(\rho)}{\log(2/\rho)}. 
\end{align*}
Evidently, $\xi_2(1)\le 1$, $\xi_3(1)\le 1$, $\lim_{t\to 0^+}\xi_2(t)=0$, 
and $\lim_{t\to 0^+}\xi_3(t)=0$.  
 
Applying Theorem \ref{Theorem 3.4} with $v$, $G$, and $g$, 
we obtain that $v\in W^{2,\text{VMO}}(B_1)$ and consequently 
$u\in W^{2,\text{VMO}}(B_{r_{x_0}}(x_0))$ with $r_{x_0}=r$.  

Choosing $\Omega_2=\bigcup_{x_0\in E} B_{r_{x_0}}(x_0)$,  
Theorem \ref{thm:Theorem 1.2} follows.             
\end{proof}

%

\setcounter{equation}{0}
\section{$W^{3,\delta}$ estimates and Hausdorff dimensions of singular sets}

\subsection{$W^{3,\delta}$ estimates}
In this subsection, we consider $W^{3,\delta}$-estimates for viscosity solutions of 
\begin{equation}\label{4.1}
F(D^2u)=f(x) 
\end{equation}
with $f$ continuous and of class $W^{1,n}$.  
Assume that $F$ has ellipticity constants $\lambda$ and $\Lambda$.  
The following is a generalization of \cite[Theorem 5.3]{CC1}. 

\begin{lemma}\label{Lemma 4.1} 
Assume that $f$, $g\in C(\Omega)$.  
Let $u$ and $v$ be in $C(\Omega)$ such that $F(D^2u)\ge f(x)$ in $\Omega$ 
and $F(D^2v)\le g(x)$ in $\Omega$ in the viscosity sense.    
Then $u-v\in \underline S(\lambda/n, \Lambda, f-g)$ in $\Omega$.  
\end{lemma}

\begin{proof}
We follow the proof of \cite[Theorem 5.3]{CC1} with appropriate modifications.  
Let $H$ and $H_1$ be open sets such that 
$\overline H_1 \subset H \subset \overline H \subset \Omega$.  
Recall the upper $\varepsilon$-envelope of $u$ with respect to $H$ defined by
$$
u^\varepsilon(x_0)=
\sup_{x\in \overline H} [u(x)+\varepsilon-\frac1\varepsilon|x-x_0|^2]
\qquad\text{for }x_0\in H. 
$$
For small $\varepsilon>0$, we claim that $u^\varepsilon$ satisfies in the viscosity sense
\begin{equation}\label{4.2}
F(D^2u^\varepsilon(x))\ge f_\varepsilon(x)
=f(x)- \text{osc}_{B_{\sqrt{\varepsilon\text{osc}_Hu}}(x)}f
\qquad\text{in }H_1.
\end{equation}
To see, let $x_0\in H_1$ and $P(x)$ be a quadratic polynomial that touches 
$u^\varepsilon$ by above at $x_0$ in $B_{r_0}(x_0)\subset H_1$. 
There exists $x_0^*\in \overline H$ such that 
$u^\varepsilon(x_0)=u(x_0^*)+\varepsilon-\frac1\varepsilon|x_0^*-x_0|^2$.  
Since $|x_0^*-x_0|^2\le \varepsilon\text{osc}_Hu$, we may assume 
$B_{r_0}(x_0^*)\subset H$.  
As in the proof of (c) in \cite[Theorem 5.1]{CC1}, the quadratic polynomial 
$Q(x)=P(x+x_0-x_0^*)-\varepsilon+\frac1\varepsilon|x_0-x_0^*|^2$
touches $u$ by above at $x_0^*$ in $B_{r_0}(x_0^*)$.  This leads to 
$$
F(D^2P)=F(D^2Q)\ge f(x_0^*) \ge f_\varepsilon(x_0).
$$
So the claim follows.  

Similarly, $v_\varepsilon=-(-v)^\varepsilon$ (the lower $\varepsilon$-envelope 
of $v$ with respect to $H$) satisfies in the viscosity sense
\begin{equation}\label{4.3}
F(D^2v_\varepsilon(x))\le g^\varepsilon(x)
=g(x)+ \text{osc}_{B_{\sqrt{\varepsilon\text{osc}_Hu}}(x)}g
\qquad\text{in }H_1.
\end{equation}

We turn to show that $u^\varepsilon-v_\varepsilon \in 
\underline S(\lambda/n, \Lambda, f_\varepsilon-g^\varepsilon)$ in $H_1$.  
Let $x_1\in H_1$ and $\Phi$ be a quadratic polynomial such that 
$u^\varepsilon-v_\varepsilon-\Phi$ attains a strict local maximum at $x_1$. 
Consider small $\overline B_{2r}(x_1)\subset H_1$ and set 
$w=(v_\varepsilon-u^\varepsilon+\Phi)-
\inf_{\overline B_{2r}(x_1) - B_r(x_1)}(v_\varepsilon-u^\varepsilon+\Phi)$.    
According to \cite[Theorem 5.1]{CC1}, for any $x\in \overline B_r(x_1)$, there 
exists a convex quadratic polynomial $P^x$ of opening $K$ that touches $w$ by 
above at $x$ in $H$.  Here $K$ is independent of $x$.  
With the notation of 
\cite[Lemma 3.5]{CC1}, it is obvious that $P^x$ also touches $\Gamma_w$ by above 
at $x\in \overline B_r(x_1)\cap\{w=\Gamma_w\}$ in $B_{2r}(x_1)$. 
By the maximum principle \cite[Lemma 3.5]{CC1}, $B_r(x_1)\cap\{w=\Gamma_w\}$ 
(contained in the lower contact set of $w$ in $B_r(x_1)$) has positive measure.  
Note that $u^\varepsilon$ and $v_\varepsilon$ are twice differentiable a.e. in 
$H$.  
So there exists one point $x_2\in B_r(x_1)\cap\{w=\Gamma_w\}$ such that 
$D^2w(x_2) \ge 0$, that is, $D^2\Phi \ge D^2u^\varepsilon(x_2)-D^2v_\varepsilon(x_2)$.  
This together with \eqref{4.2} and \eqref{4.3} yields
\begin{align*}
\mathcal P^+_{\lambda/n, \Lambda}(D^2\Phi)
&\ge \mathcal P^+_{\lambda/n, \Lambda}(D^2u^\varepsilon(x_2)-D^2v_\varepsilon(x_2))     \\ 
&\ge F(D^2u^\varepsilon(x_2))-F(D^2v_\varepsilon(x_2))
\ge f_\varepsilon(x_2) - g^\varepsilon(x_2).
\end{align*}
Letting $r \longrightarrow 0$, we obtain 
$\mathcal P^+_{\lambda/n, \Lambda}(D^2\Phi) \ge f_\varepsilon(x_1) - g^\varepsilon(x_1)$.  

It follows from the closedness of the class $\underline S$ and uniform convergence 
of $u^\varepsilon-v_\varepsilon$ to $u-v$ 
that $u-v \in \underline S(\lambda/n, \Lambda, f-g)$ in $H_1$.  
This completes the proof of Lemma \ref{Lemma 4.1}.         
\end{proof}

\begin{corollary}\label{Corollary 4.2} 
Let $u$ be a viscosity solution to \eqref{4.1} in $\Omega$ with $f\in C(\Omega)$.  
Let $h\in\mathbb R^n$ be such that 
$\Omega_{|h|}=\left\{x\in\Omega: d(x, \partial\Omega)>|h|\right\}\neq \emptyset$.   
Then 
\[
u(x+h)-u(x)\in S\(\lambda/n, \Lambda, f(x+h)-f(x)\)
\] 
in $\Omega_{|h|}$.    
\end{corollary}

\begin{lemma}\label{Lemma 4.3} 
Let $u$ be a viscosity solution to \eqref{4.1} in $B_1$.  
Suppose that $f\in C(B_1)\cap W^{1,n}(B_1)$ and $F(0)=0$. 
Then $u\in C^{1, \alpha_1}(\overline B_{1/2})$ and   
$$
\|u\|_{C^{1, \alpha_1}(\overline B_{1/2})}
\le C(\|u\|_{L^\infty(B_1)} + \|f\|_{W^{1,n}(B_1)}), 
$$ 
where $0<\alpha_1<1$ and $C$ are universal constants.      
\end{lemma}

\begin{proof}  
The proof of the lemma is similar to that of \cite[ Corollary 5.7]{CC1}.  
For $e\in \mathbb R^n$ with $|e|=1$, $0<h<1/8$, and $0<\beta\le 1$, 
denote $\Delta_{\beta, he}u(x)=[u(x+he)-u(x)]/h^\beta$.  
In virtue of Corollary \ref{Corollary 4.2}, we have 
$$
\Delta_{\beta, he}u\in S(\lambda/n, \Lambda, \Delta_{\beta, he}f) 
\qquad\text{in }B_{7/8}.
$$
H\"older estimates (\cite[ Proposition 4.10]{CC1}) give 
\begin{equation}\label{4.4}
\|\Delta_{\beta, he}u\|_{C^{\alpha_1}(\overline B_r)}
\le C(r,s)(\|\Delta_{\beta, he}u\|_{L^\infty(B_{(r+s)/2})} 
+ \|\Delta_{\beta, he}f\|_{L^n(B_{(r+s)/2})}),  
\end{equation}
where $0<r<s\le 7/8$, $0<h<(s-r)/2$, $0<\alpha_1<1$ is universal and 
$C(r,s)$ depends only on $n$, $\lambda$, $\Lambda$, $r$, $s$.   

Set $K=\|u\|_{L^\infty(B_1)} + \|f\|_{W^{1,n}(B_1)}$.
Since $u\in S(\lambda/n, \Lambda, f)$ in $B_1$, \cite[Proposition 4.10]{CC1} 
again yields $\|u\|_{C^{\alpha_1}(\overline B_{7/8})}\le CK$.  

It follows from \eqref{4.4} with $\beta=\alpha_1$ that  
$$
\|\Delta_{\alpha_1, he}u\|_{C^{\alpha_1}(\overline B_{r_1})}
\le C(r_1)(\|u\|_{C^{\alpha_1}(\overline B_{7/8})} + \|Df\|_{L^n(B_{7/8})})\le C(r_1)K  
$$
for $r_1<7/8$ and $0<h<(7/8-r_1)/2$.  By \cite[Lemma 5.6]{CC1},  
the $C^{\alpha_1}$ estimates of $u$ and $\Delta_{\alpha_1, he}u$ 
imply that $\|u\|_{C^{2\alpha_1}(\overline B_{r_2})}\le C(r_1, r_2)K$ for $r_2<r_1$.  
Repeating this process to keep improving the H\"older exponent of $u$, we eventually 
obtain that $u$ is Lipschitz and $\|u\|_{C^{0,1}(\overline B_{3/4})}\le CK$. 
Therefore, \eqref{4.4} with $\beta=1$ gives rise to 
$$
\|\Delta_{1, he}u\|_{C^{\alpha_1}(\overline B_{1/2})}
\le C(\|u\|_{C^{0,1}(\overline B_{3/4})} + \|Df\|_{L^n(B_{3/4})}) \le CK  
$$
for $0<h<1/8$ and $|e|=1$.  Thus, Lemma \ref{Lemma 4.3} easily follows.           
\end{proof}

Let $H\subset \Omega$ be open and $t>0$.  Denote 
$$
G_t(u, H)=\left\{x_0\in H: \exists \, L \text{ affine } \text{ such that }
|(u-L)(x)|\le\frac t2|x-x_0|^2\ \ \forall\ x\in H\right\}
$$
and $A_t(u, H)=H-G_t(u, H)$.  We also define 
$$
\Theta_3(u, H)(x_0)=\inf\left\{t>0: \exists \,P\, \text{ quadratic } \text{ such that }
|(u-P)(x)|\le\frac t6|x-x_0|^3\ \ \forall\ x\in H\right\}
$$
for $x_0\in H$. We define $\Theta_3(u, H)(x_0)=\infty$ if no such quadratic polynomial 
$P$ exists.  Obviously, $\Theta_3(u, H)$ is a measurable function in $H$.  

\begin{theorem}\label{Theorem 4.4} 
$(W^{3,\delta} \text{ estimates})$
Let $u$ be a viscosity solution to $F(D^2u)=f(x)$ in $B_{8\sqrt n}$.  
Suppose that $F$ is uniformly elliptic with constants $0<\lambda \le \Lambda$ and   
$F(0)=0$.  Assume that $f\in C(B_{8\sqrt n})\cap W^{1,n}(B_{8\sqrt n})$. 
Then   
$$
|\{x\in Q_{1/2}:  \Theta_3(u, B_{8\sqrt n})(x)>t\}|
\le C(\|u\|_{L^\infty(B_{8\sqrt n})} + \|f\|_{W^{1,n}(B_{8\sqrt n})})^\delta t^{-\delta}
\qquad \forall\ t>0, 
$$  
where $\delta$ is as in \eqref{eq:Ldelta norm of Hessian} and $C$ is universal.        
\end{theorem}

\begin{proof}  
In view of Corollary \ref{Corollary 4.2}, for $0<|h|<\sqrt n$, the difference quotient of $u$ given by
$\Delta_hu(x)=[u(x+h)-u(x)]/|h|$ satisfies that  
$\Delta_hu\in S(\lambda/n, \Lambda, \Delta_hf)$ in $B_{7\sqrt n}$.  
Let $t>0$ and $K=\|u\|_{L^\infty(B_{8\sqrt n})} + \|f\|_{W^{1,n}(B_{8\sqrt n})}<\infty$.  
It follows from \cite[Lemma 7.8]{CC1} and Lemma \ref{Lemma 4.3} that   
\begin{equation}\label{4.5}
|A_t(\Delta_hu, B_{6\sqrt n})\cap Q_{1/2}| 
\le C(\|\Delta_hu\|_{L^\infty(B_{6\sqrt n})}+\|\Delta_hf\|_{L^n(B_{6\sqrt n})})^\delta t^{-\delta}  
\le CK^\delta t^{-\delta}
\end{equation}
with $\delta$ and $C$ universal.  We can further take $\delta$ as in \eqref{eq:Ldelta norm of Hessian}.   

Let $e_i$ be the $i^{\text{th}}$ standard unit vector and $h_k=\frac1k e_i$ for $k\ge 1$.   
Note that $\Delta_{h_k}u \longrightarrow D_iu$ in $C(\overline B_{6\sqrt n})$.  We claim
\begin{equation}\label{4.6}
\varlimsup_{k\to \infty} G_t(\Delta_{h_k}u, B_{6\sqrt n})\subset G_t(D_iu, B_{6\sqrt n}).
\end{equation}
Indeed, if $x_0\in \varlimsup\limits_{k\to \infty} G_t(\Delta_{h_k}u, B_{6\sqrt n})$, then there 
exists $\{p_{k_j}\}\subset \mathbb R^n$ such that 
$|\Delta_{h_{k_j}}u(x)-\Delta_{h_{k_j}}u(x_0)-p_{k_j}\cdot (x-x_0)|\le \frac{t}2|x-x_0|^2$ 
for $x\in B_{6\sqrt n}$, $j\ge 1$.  
Since $\{p_{k_j}\}$ must be bounded, there is a
convergent subsequence $\{p_{k_{j_l}}\}$.  Therefore, we obtain 
$x_0\in G_t(D_iu, B_{6\sqrt n})$.  With \eqref{4.6} and \eqref{4.5}, it is easy to see that
\begin{equation}\label{4.7}
|A_t(D_iu, B_{6\sqrt n})\cap Q_{1/2}| 
\le \varliminf_{k\to \infty} |A_t(\Delta_{h_k}u, B_{6\sqrt n})\cap Q_{1/2}| 
\le CK^\delta t^{-\delta}.
\end{equation}

Now consider $y_0\in \cap_{i=1}^n G_t(D_iu, B_{6\sqrt n}) \cap Q_{1/2}$.  
Then there exist numbers $p_{ij}$ with $1\le i$, $j\le n$ such that 
$|D_iu(x)-D_iu(y_0)-\sum_{j=1}^np_{ij}(x_j-y_{0j})| \le \frac{t}2|x-y_0|^2$ for 
$x\in B_{6\sqrt n}$.  Here $x_j$ and $y_{0j}$ are the $j^{\text{th}}$-coordinates 
of $x$ and $y_0$ respectively.    
Let $P(x)=u(y_0)+Du(y_0)\cdot(x-y_0)+\frac12\sum\limits_{i,j}p_{ij}(x_i-y_{0i})(x_j-y_{0j})$. 
For $x\in B_{6\sqrt n}$, we have 
\begin{align*}
|u(x)&-P(x)|  \\
&=\bigg |\int_0^1 \sum_{i=1}^n \big[ D_iu(y_0+\theta(x-y_0))-D_iu(y_0)
   -\sum_{j=1}^n p_{ij}\theta(x_j-y_{0j}) \big](x_i-y_{0i})\,d\theta \bigg |           \\
&\le \int_0^1 \sum_{i=1}^n\frac{t}2\theta^2|x-y_0|^2|x_i-y_{0i}| \,d\theta
   \le \frac{\sqrt{n}t}6|x-y_0|^3.
\end{align*}
By considering $P(y_0+e)\pm P(y_0-e)$ with $|e|=1$, we obtain 
$|P(y_0)|+ |DP(y_0)|+ \|D^2P\| \le C_n(K+t)$, where $C_n$ is a dimensional constant. 
Therefore, $|u(x)-P(x)|\le C_n(K+t)|x-y_0|^3$ for $x\in B_{8\sqrt n}-B_{6\sqrt n}$.  
We have shown that $\cap_{i=1}^n G_t(D_iu, B_{6\sqrt n}) \cap Q_{1/2} \subset 
\{x\in Q_{1/2}:  \Theta_3(u, B_{8\sqrt n})(x)\le C_n(K+t)\}$.  Consequently,   
\eqref{4.7} yields 
$$
|\{x\in Q_{1/2}:  \Theta_3(u, B_{8\sqrt n})(x)>C_n(K+t)\}| \le CK^\delta t^{-\delta}. 
$$
Theorem \ref{Theorem 4.4} then follows.                 
\end{proof}

\subsection{Hausdorff dimensions of singular sets}
We now use Theorems \ref{Theorem 3.4} and \ref{Theorem 4.4} to establish the Hausdorff dimension estimates 
of singular sets in Theorem \ref{thm:Theorem 1.3}.

{\bf Proof of Theorem \ref{thm:Theorem 1.3}}    
We first assume that $\overline B_{8\sqrt n}\subset \Omega$. 
Let $\eta_4$ and $\sigma_4$ be the constants given in Theorem \ref{Theorem 3.4}.    
Also let $x_0\in Q_{1/4}$ and $r\le r_0 \le 1/4$ with small $r_0$ 
to be determined later.  We claim that 
\begin{equation}\label{4.8}
B_r(x_0)\cap \{x\in B_{8\sqrt n}:  \Theta_3(u, B_{8\sqrt n})(x)\le \frac{\eta_4}{9r}\}
\neq\emptyset 
\, \ \Longrightarrow\, \  u\in W^{2,\text{VMO}}(B_r(x_0)).
\end{equation}
To see, consider $x_1\in B_r(x_0)$ such that 
$\Theta_3(u, B_{8\sqrt n})(x_1)\le \frac{\eta_4}{9r}$.  
Then there exists a quadratic 
polynomial $P$ such that $|u(x)-P(x)| \le \frac{\eta_4}{54r}|x-x_1|^3$ in $B_{8\sqrt n}$.  
Clearly, $F(D^2P)=f(x_1)$ and 
$$
|u(x)-P(x)| \le \frac{\eta_4}9 |x-x_1|^2 \qquad\text{in }B_{6r}(x_1)\subset B_{2\sqrt n}.  
$$
It is easy to check that $v(y)=(u-P)(2ry+x_1)/(2r)^2$ is a viscosity solution to 
$$
F(D^2v(y)+D^2P)-f(x_1)=f(2ry+x_1)-f(x_1) \qquad\text{in }B_3. 
$$  
As in the proof of Theorem \ref{thm:Theorem 1.2}, 
$G(M)=F(M+D^2P)-f(x_1)\in \mathbb T^\omega(\mathcal S, \lambda, \Lambda)$.  
Also $\|v\|_{L^\infty(B_3)} \le \eta_4$.  
For $|y_0|\le 2$, $0<\rho\le 1$, the function $g(y)= f(2ry+x_1)-f(x_1)$ satisfies 
$$
\text{osc}_{\overline B_\rho(y_0)}g\le 
\sup_{z_1, z_2\in \overline B_{8\sqrt n} \\ |z_1-z_2|\le 4r_0\rho} |f(z_1)-f(z_2)|
=\xi(4r_0\rho)
$$
with $\xi$ denoting the modulus of continuity of $f$ over $\overline B_{8\sqrt n}$.  
Choose $r_0$ sufficiently small so that $\xi(4r_0)\le \sigma_4\eta_4$.  
By Theorem \ref{Theorem 3.4}, we obtain $v\in W^{2,\text{VMO}}(B_1)$.  Hence,  
$u\in W^{2,\text{VMO}}(B_r(x_0))$ as desired.    

Denote 
$$
E_0=\{x\in Q_{1/4}: \exists\ r_x>0 \text{ such that }
u\in W^{2,\text{VMO}}(B_{r_x}(x))\}
$$
and $\Sigma_0=Q_{1/4}-E_0$.  Since $\Sigma_0$ is compact, for $0<r\le r_0$, 
there exist points $\{z_j\}_{j=1}^N$ in $\Sigma_0$ such that 
$\Sigma_0\subset \cup_{j=1}^N B_r(z_j)$.  
By Vitali's lemma, we can assume that 
the first $m$ balls $B_r(z_1)$, $\cdots$, $B_r(z_m)$ (where $m\le N$)  
are pairwise disjoint, and $\cup_{j=1}^N B_r(z_j) \subset \cup_{j=1}^m B_{3r}(z_j)$.  
Because $z_j\notin E_0$, \eqref{4.8} implies that 
$B_r(z_j)\subset \{x\in Q_{1/2}:  \Theta_3(u, B_{8\sqrt n})(x)>\frac{\eta_4}{9r}\}$.  
By Theorem \ref{Theorem 4.4}, we obtain 
$$
\sum_{j=1}^m |B_r(z_j)| 
\le |\{x\in Q_{1/2}:  \Theta_3(u, B_{8\sqrt n})(x)>\frac{\eta_4}{9r}\}|
\le Cr^\delta,
$$
where $C$ depends only on $n$, $\lambda$, $\Lambda$, $\omega$, $F(0)$,
$\|u\|_{L^\infty(B_{8\sqrt n})}$, and $\|f\|_{W^{1,n}(B_{8\sqrt n})}$.   
It follows immediately that 
$$
\sum_{j=1}^m \big(\text{diam }(B_{3r}(z_j))\big)^{n-\delta} \le C.
$$
Moreover, $\Sigma_0\subset \cup_{j=1}^m B_{3r}(z_j)$.  Therefore, 
$\mathcal H^{n-\delta}(\Sigma_0)\le C$, i.e., the $(n-\delta)$-dimensional Hausdorff 
measure of $\Sigma_0$ is finite.  We conclude that 
\begin{equation}\label{4.9}
\mathcal H_{dim}(\Sigma_0)\le n-\delta.
\end{equation}    

Let 
$
\Omega_2=\{x\in \Omega: \exists\ r_x>0 \text{ such that }
u\in W^{2,\text{VMO}}(B_{r_x}(x))\} 
$ 
denote the regular set of $u$.  For general $\Omega$,   
by standard rescaling and covering arguments,   
we obtain from \eqref{4.9} that $\mathcal H_{dim}(\Omega-\Omega_2)\le n-\delta$.     
Theorem \ref{thm:Theorem 1.3} is then proved.             
\qed

\setcounter{equation}{0}
\section{Nondifferentiable fully nonlinear equations satisfying (H2)}\label{sec:nondiff fully nonlinear satisfying (H2)}

The purpose of this section is to prove Theorem \ref{thm:Theorem 1.4}.

We first show that that either \hyperlink{A1}{\bf(A1)\rm} or \hyperlink{A2}{\bf(A2)\rm} implies \hyperlink{H2}{\bf(H2)\rm}.

If $F: \mathcal S \longrightarrow \Bbb R$ is of class $C^1$ and $DF$ has a modulus of continuity $\tau$, then the Taylor Theorem gives
\begin{equation}\label{(5.1)}
 |F(M)-F(M_0)-tr(DF(M_0)(M-M_0))|\le \sqrt{n}\,\tau(\|M-M_0\|)\|M-M_0\|
\end{equation}
for any $M, M_0\in \mathcal S$.  
Suppose $F(M, x)$ satisfies \hyperlink{A1}{\bf(A1)\rm}.  For $M_0\in \mathcal S$ and $x_0\in \Omega$, set 
\[
T_{M_0,x_0}(M)=F_\cup(M, x_0)+ F_C(M_0,x_0)+\text{tr}(D_MF_C(M_0, x_0)(M-M_0)).
\]
Then \eqref{(5.1)} implies that $F(\cdot, x_0)\in \Bbb T^{\sqrt{n}\tau}(\mathcal S, \lambda, 2\Lambda)$ for any $x_0\in \Omega$.  On the other hand, suppose $F(M, x)$ satisfies \hyperlink{A2}{\bf(A2)\rm}. In this case, set the tangent operator 
\[
T_{M_0,x_0}(M)=\inf_{\gamma\in\Gamma}[F_\gamma(M_0, x_0)+ \text{tr}(D_MF_\gamma(M_0, x_0)(M-M_0))].
\]
Again \eqref{(5.1)} yields that $F(\cdot, x_0)\in \Bbb T^{\sqrt{n}\tau}(\mathcal S, \lambda, \Lambda)$ for any $x_0\in \Omega$.

\subsection{Semiconvex equations}
The purpose of this section is to show the following proposition. As a consequence, \hyperlink{A3}{\bf(A3)\rm} implies \hyperlink{H2}{\bf(H2)\rm}.

\begin{proposition}\label{prop:semiconvex implies H1}
Let $F=F(M)$ be satisfying \eqref{eq:ellipticity condition}, with $F$ semi-convex, i.e., there exists a constant $K>0$ such that $F(M)+K\,|M|^2$ is convex in $\mathcal S$.
Then  
$F\in \Bbb T^{\omega}(\mathcal S, \lambda/2, 2\Lambda)$ with $\omega(t)=C\,t$ with $C=C(n, \lambda, \Lambda, K)$.
\end{proposition}
\begin{proof}
For $A,B\in \R^{n\times n}$ denote $\langle A,B \rangle=\text{tr}\(A^t\,B\)$ the Frobenius inner product and the norm $|M|=\(\text{tr}\(M^t\,M\)\)^{1/2}=\(\sum_{ij}m_{ij}^2\)^{1/2}$.
Clearly the norms $|\cdot |$ and $\|\cdot \|$ are equivalent.
Since $F:\mathcal S\to \R$ we can write 
\[
F(M)=g\(m_{11};m_{21},m_{22};m_{31},m_{32},m_{33};\cdots ;m_{n1},m_{n2},m_{n3}\,\cdots ,m_{nn}\)
\]
for $M=[m_{ij}]\in \mathcal S$ with $g$ a function of $n(n+1)/2$ variables.

Given $M\in \R^{n\times n}$ we set $\overline{F}(M)=F\(\dfrac{M+M^t}{2}\)$.

Since $F$ is semi-convex, by Rademacher's theorem $F(M)+K\,|M|^2$ is differentiable a.e. in $\mathcal S$; 
and since $|M|^2$ is differentiable everywhere, $F$ is differentiable a.e in $\mathcal S$.
If $M_1\in \mathcal S$ is a differentiability point of $F$, then $M_1$ is a differentiability point of $\overline{F}$ and the matrix of its derivatives $D\overline{F}(M_1)$ is symmetric and satisfies 
\[
\dfrac{\partial \overline{F}}{\partial m_{ii}}(M_1)=\dfrac{\partial F}{\partial m_{ii}}(M_1)\quad 1\leq i \leq n,
\]
and 
\[
\dfrac{\partial \overline{F}}{\partial m_{ij}}(M_1)=\dfrac12 \dfrac{\partial F}{\partial m_{ij}}(M_1)\quad i> j ,
\]

By \eqref{eq:ellipticity condition}, at a differentiability point $M_1$ of $F$ we have 
\[
\lambda \,t\|N\|\leq F(M_1+t\,N)-F(M_1)\leq \Lambda \,t\|N\|
\]
for all $N\geq 0$, $N\in \mathcal S$ and $t>0$.
Letting $t\to 0^+$ we get 
\begin{equation}\label{eq:differentiability of F}
\lambda \,\|N\|\leq \text{tr}\(D\overline{F}(M_1)N\)\leq \Lambda \,\|N\|
\end{equation}
for all $N\in \mathcal S$, $N\geq 0$.
Now, let $v$ be an $n$ column vector and let $N=vv^t$. Then $N\geq 0$, and $\|N\|=|v|^2$. Noticing that 
$\text{tr}\(D\overline{F}(M_1)vv^t\)=v^tD\overline{F}(M_1)v$, we obtain applying \eqref{eq:differentiability of F} that 
$\lambda \leq D\overline{F}(M_1)\leq \Lambda$.

Next fix $M_0\in \mathcal S$ and define $G(M)=\overline{F}(M)+K\,|M-M_0|^2$ for $M\in \R^{n\times n}$. Since $F$ is semi-convex, $G$ is convex in $\mathcal S$
and the set of points where $G$ is differentiable is the same as the set of points in $\mathcal S$ where $F$ is differentiable. 
Then for each $M_1\in \mathcal S$ where $G$ is differentiable, we have 
\begin{equation}\label{eq:convexity of G}
G(M)\geq G(M_1)+\text{tr}\(DG(M_1)(M-M_1)\):=L_{M_1}(M)\qquad \forall M\in \mathcal S,
\end{equation}  
where $DG(M_1)$ is the matrix whose entries are $\dfrac{\partial G}{\partial m_{ij}}(M_1)$;
and 
\[
DG(M_1)=D\overline{F}(M_1)+2\,K\(M_1-M_0\).
\]
Hence for $N\geq 0$, 
\begin{align*}
\text{tr}\(DG(M_1)N\)
&\leq 
\Lambda\,\|N\|+2\,K\,\text{tr}\(\(M_1-M_0\)N\)\\
&\leq 
\(\Lambda +2\,K\,n\,\|M_1-M_0\|\)\|N\|,
\end{align*}
and
\begin{align*}
\text{tr}\(DG(M_1)N\)
&\geq 
\(\lambda -2\,K\,n\,\|M_1-M_0\|\)\|N\|.
\end{align*}
If $\|M_1-M_0\|\leq \lambda/4Kn$ we obtain
\[
\lambda/2\,\|N\|\leq \text{tr}\(DG(M_1)N\)\leq 2\Lambda\,\|N\|.
\]
Therefore 
\[
\lambda/2\,\|N\|\leq L_{M_1}(M+N)-L_{M_1}(M)\leq 2\Lambda\,\|N\|
\]
for all $M\in \mathcal S$ and $N\geq 0$,
so $L_{M_1}$ is uniformly elliptic for each $M_1\in \mathcal B_{\lambda/4Kn}(M_0)$, $M_1$ differentiability point of $G$ (or $F$).
Since $G$ is differentiable a.e. in $\mathcal S$, the set $E$ of points of $\mathcal B_{\lambda/4Kn}(M_0)$ where $G$ is differentiable is dense in $\mathcal B_{\lambda/4Kn}(M_0)$,
and we define
\[
T_{M_0}(M)=\sup_{M_1\in E}L_{M_1}(M).
\]

From \eqref{eq:convexity of G}, $T_{M_0}(M)$ is well defined for $M\in \mathcal S$, is convex in $\mathcal S$, and is elliptic with constants $\lambda/2$ and $2\Lambda$. Notice also that  $G(M)=T_{M_0}(M)$ for $M\in E$ and by continuity $G(M)=T_{M_0}(M)$ in $\mathcal B_{\lambda/4Kn}(M_0)$.
Hence $F(M)-T_{M_0}(M)=-K\,|M-M_0|^2$ for all $M \in \mathcal B_{\lambda/4Kn}(M_0)$.
If $\lambda/4Kn<\|M-M_0\|\leq 1$, then 
\begin{align*}
|F(M)-T_{M_0}(M)|
&= 
|F(M)-F(M_0)+T_{M_0}(M_0)-T_{M_0}(M)|\\
&\leq 
3\Lambda\,\|M-M_0\| \leq \(12 n \Lambda K\)/\lambda \,\|M-M_0\|^2,
\end{align*}
since the ellipticity condition \eqref{eq:ellipticity condition} implies that $F$ and $T_{M_0}$ have Lipschitz constants $\Lambda$ and $2\Lambda$, respectively.  
Therefore,
$F\in \Bbb T^{\omega}(\mathcal S, \lambda/2, 2\Lambda)$ with $\omega(t)=12n\Lambda K/\lambda t$. 

\end{proof}

\subsection{Piecewise differentiable equations}\label{subsect:piecewise-diff-equations}
%
%
We want to show that \hyperlink{A4}{\bf(A4)\rm} implies \hyperlink{H2}{\bf(H2)\rm}.  We begin with a flat-interface lemma and then straighten the interface by a local $C^1$ change of variables.

Throughout this subsection, the norm and distance in $\mathcal S$ are those induced by the Frobenius inner product
\[
\langle A,B\rangle=\text{tr}(A^tB),\qquad |A|=\langle A,A\rangle^{1/2}.
\]
These norms are equivalent to the operator norm $\|\cdot\|$ used in the ellipticity condition.  We denote by $\mathcal B_r(M_0)$ the ball in $\mathcal S$ centered at $M_0$ with radius $r$ in the Frobenius norm.

Recall that if $F$ is differentiable at $M_0\in\mathcal S$, then there exists a unique $A\in\mathcal S$ such that
\[
F(M)=F(M_0)+\langle A,M-M_0\rangle+o(|M-M_0|)
\quad\text{as }M\to M_0.
\]
Equivalently, identifying $F$ with its symmetric extension
$\overline F(M)=F\big((M+M^t)/2\big)$, we write
$A=DF(M_0)=[D_{ij}F(M_0)]$.

Let $\nu\in\mathcal S$ be a unit vector and let
\[
\Pi_\nu=\{M\in\mathcal S:\langle M,\nu\rangle=0\},\qquad
\mathcal S_\nu^+=\{M\in\mathcal S:\langle M,\nu\rangle>0\},\qquad
\mathcal S_\nu^-=\{M\in\mathcal S:\langle M,\nu\rangle<0\}.
\]
If $M_0^*\in\Pi_\nu$ and $R>0$, set
\[
\mathcal B^*=\mathcal B_{2R}(M_0^*),\qquad
\mathcal B=\mathcal B_R(M_0^*),
\]
and write
$\mathcal B^*_\pm=\mathcal B^*\cap\mathcal S_\nu^\pm$ and
$\mathcal B_\pm=\mathcal B\cap\mathcal S_\nu^\pm$.

\begin{lemma}\label{lem:flat-interface}
Let $F=F(M)$ satisfy \eqref{eq:ellipticity condition}.  Assume that
\[
F\in C^1(\mathcal B^*\cap\mathcal S_\nu^+)\cap C^1(\mathcal B^*\cap\mathcal S_\nu^-)
\]
and that $DF$ has a modulus of continuity $\tau$ on
$\mathcal B^*\cap\mathcal S_\nu^+$ and on
$\mathcal B^*\cap\mathcal S_\nu^-$, respectively.  Then there exists
$R_0=R_0(n,\lambda,\Lambda,\tau)>0$ such that, if $R\le R_0$, then
\[
F\in\Bbb T^\omega(\overline{\mathcal B},\lambda/2,2\Lambda),
\]
where $\omega(t)=C(1+1/\tau(R))\tau(t)$ and $C=C(n,\lambda, \Lambda).$ 
\end{lemma}

\begin{proof}
Let $\{e_i:1\le i\le n^*-1\}$ be an orthonormal basis of $\Pi_\nu$, where
$n^*=n(n+1)/2$.  Since $DF$ is uniformly continuous in each of
$\mathcal B^*_+$ and $\mathcal B^*_-$, it extends continuously to
$\overline{\mathcal B^*_+}$ and $\overline{\mathcal B^*_-}$, respectively.
We write $D_{e_i}F=\langle DF,e_i\rangle$ and denote the one-sided normal derivatives in
$\overline{\mathcal B^*_+}$ and $\overline{\mathcal B^*_-}$ by
$D_\nu^+F$ and $D_\nu^-F$.

First let $M_1\in\overline{\mathcal B_+}$ and denote by
\[
\overline M_1=M_1-\langle M_1,\nu\rangle\nu
\]
the projection of $M_1$ onto $\Pi_\nu$.  Define two affine functions
\begin{align*}
L^+(M)&=F(M_1)+\sum_{i=1}^{n^*-1}D_{e_i}F(M_1)\langle M-M_1,e_i\rangle
       +D_\nu^+F(M_1)\langle M-M_1,\nu\rangle,\notag\\
L^-(M)&=F(M_1)+\sum_{i=1}^{n^*-1}D_{e_i}F(M_1)\langle M-M_1,e_i\rangle
       +D_\nu^+F(M_1)\langle \overline M_1-M_1,\nu\rangle\notag\\
&\hspace{3.1cm}+D_\nu^-F(\overline M_1)\langle M-\overline M_1,\nu\rangle.
\end{align*}
Then set
\[
T_{M_1}(M)=
\begin{cases}
L^+(M),& \text{if }\langle M,\nu\rangle\ge 0,\\
L^-(M),& \text{if }\langle M,\nu\rangle<0.
\end{cases}
\]
The two affine functions agree on $\Pi_\nu$.  Moreover, if
$D_\nu^+F(M_1)\ge D_\nu^-F(\overline M_1)$, then $L^+(M)\ge L^-(M)$ in
$\mathcal S_\nu^+$ and $L^-(M)\ge L^+(M)$ in $\mathcal S_\nu^-$.  Hence
$T_{M_1}=\max\{L^+,L^-\}$ on $\mathcal S$.  If
$D_\nu^+F(M_1)<D_\nu^-F(\overline M_1)$, then
$T_{M_1}=\min\{L^+,L^-\}$ on $\mathcal S$.  Consequently $T_{M_1}$ is either
convex or concave in $\mathcal S$.

We next prove that $T_{M_1}$ is tangent to $F$ at $M_1$.  For
$M\in\overline{\mathcal B^*_+}$, Taylor's theorem gives
\begin{equation}\label{eq:flat-tangent-plus}
\begin{aligned}
|F(M)-T_{M_1}(M)|
&=\left|F(M)-\bigl[F(M_1)+\langle D^+F(M_1),M-M_1\rangle\bigr]\right|  \\
&\le \tau(|M-M_1|)|M-M_1|.
\end{aligned}
\end{equation}
For $M\in\overline{\mathcal B^*_-}$, first observe that
\[
L^-(M)=F(\overline M_1)+E_1+\sum_{i=1}^{n^*-1}D_{e_i}F(M_1)
        \langle M-\overline M_1,e_i\rangle
        +D_\nu^-F(\overline M_1)\langle M-\overline M_1,\nu\rangle,
\]
where
\[
|E_1|\le \tau(|\overline M_1-M_1|)|\overline M_1-M_1|.
\]
Therefore
\[
L^-(M)=F(\overline M_1)+E_1+E_2+\sum_{i=1}^{n^*-1}D_{e_i}F(\overline M_1)
        \langle M-\overline M_1,e_i\rangle
        +D_\nu^-F(\overline M_1)\langle M-\overline M_1,\nu\rangle,
\]
with
\[
|E_2|\le C(n)\tau(|\overline M_1-M_1|)|M-\overline M_1|.
\]
Since $|M-\overline M_1|\le |M-M_1|$ and
$|\overline M_1-M_1|\le |M-M_1|$, we obtain
\begin{equation}\label{eq:flat-tangent-minus}
\begin{aligned}
|F(M)-T_{M_1}(M)|
&\le C(n)\tau(|M-M_1|)|M-M_1|.
\end{aligned}
\end{equation}
Together with the Lipschitz continuity of $F$ and $T_{M_1}$, estimates
\eqref{eq:flat-tangent-plus} and \eqref{eq:flat-tangent-minus} imply that
\[
|F(M)-T_{M_1}(M)|\le C(1+1/\tau(R))\tau(|M-M_1|)|M-M_1|
\qquad\text{for } |M-M_1|\le \sqrt n,
\]
where $C$ depends only on $n$, $\lambda$, and $\Lambda$.

It remains to prove the ellipticity of $T_{M_1}$.  Let $M,N\in\mathcal S$ with
$N\ge0$.  Choose $\{M_k\}\subset\mathcal B_+$ such that $M_k\to M_1$.  From
\eqref{eq:differentiability of F} with $M_k$, and the continuity of $D^+F$, we get
\[
\lambda\|N\|\le\langle D^+F(M_1),N\rangle\le\Lambda\|N\|.
\]
Similarly,
\[
\lambda\|N\|\le\langle D^-F(\overline M_1),N\rangle\le\Lambda\|N\|.
\]
Choose $R_0$ sufficiently small so that $C(n)\tau(R_0)\le \lambda/2$.  Then
\[
\left|\sum_{i=1}^{n^*-1}\bigl(D_{e_i}F(M_1)-D_{e_i}F(\overline M_1)\bigr)
\langle N,e_i\rangle\right|
\le C(n)\tau(R_0)\|N\|\le {\lambda\over2}\|N\|.
\]
Therefore
\[
{\lambda\over2}\|N\|\le L^-(M+N)-L^-(M)\le 2\Lambda\|N\|.
\]
The same bounds for $L^+$ are immediate.  Hence $T_{M_1}$ is uniformly elliptic
with constants $\lambda/2$ and $2\Lambda$.

Finally, if $M_1\in\overline{\mathcal B_-}$, the analogous construction is obtained by setting
\begin{align*}
\widetilde L^+(M)&=F(M_1)+\sum_{i=1}^{n^*-1}D_{e_i}F(M_1)\langle M-M_1,e_i\rangle
       +D_\nu^-F(M_1)\langle \overline M_1-M_1,\nu\rangle
       +D_\nu^+F(\overline M_1)\langle M-\overline M_1,\nu\rangle,\\
\widetilde L^-(M)&=F(M_1)+\sum_{i=1}^{n^*-1}D_{e_i}F(M_1)\langle M-M_1,e_i\rangle
       +D_\nu^-F(M_1)\langle M-M_1,\nu\rangle,
\end{align*}
where now $\overline M_1$ is again the projection of $M_1$ onto $\Pi_\nu$, and by taking
\[
\widetilde T_{M_1}(M)=
\begin{cases}
\widetilde L^+(M),& \text{if }\langle M,\nu\rangle\ge 0,\\
\widetilde L^-(M),& \text{if }\langle M,\nu\rangle<0.
\end{cases}
\]
The preceding argument shows that $\widetilde T_{M_1}$ is either convex or concave,
has ellipticity constants $\lambda/2$ and $2\Lambda$, and is tangent to $F$ at
$M_1$ with the same modulus.  The lemma follows.
\end{proof}

\begin{corollary}\label{cor:flat-global-interface}
Let $F=F(M)$ satisfy \eqref{eq:ellipticity condition}.  Assume that
$F\in C^1(\mathcal S^+)\cap C^1(\mathcal S^-)$, where
\[
\mathcal S^+=\{M\in\mathcal S:m_{nn}>0\},\qquad
\mathcal S^-=\{M\in\mathcal S:m_{nn}<0\},
\]
and assume that $DF$ has a modulus of continuity $\tau$ on
$\mathcal S^+$ and $\mathcal S^-$, respectively.  Then
\[
F\in\Bbb T^\omega(\mathcal S,\lambda/2,2\Lambda),
\]
where $\omega$ depends only on $n$, $\lambda$, $\Lambda$, and $\tau$.
\end{corollary}

\begin{proof}
The result follows from Lemma \ref{lem:flat-interface} when the point is near the hyperplane
$\{m_{nn}=0\}$.  Away from that hyperplane, $F$ is tangent to its affine
linearization by Taylor's theorem.  The two moduli can be combined into a single
modulus $\omega$ depending only on the stated data.
\end{proof}

Let $\Xi$ be a closed surface in $\mathcal S$, bounded or unbounded, such that
$\mathcal S\setminus\Xi$ is partitioned into two disjoint domains
$\mathcal D_\Xi^+$ and $\mathcal D_\Xi^-$.  We shall use the following uniform
$C^1$ condition for $\Xi$.

\medskip
\noindent{\bf(S)} There exist a modulus of continuity $\tau_1$ and constants
$K,d,d_0>0$ such that, for every $M_0\in\Xi$, there is a $C^1$ diffeomorphism
$\psi=\psi_{M_0}$ defined in $\mathcal B_{d_0}(M_0)$ satisfying
\begin{align*}
&\psi(\mathcal B_{d_0}(M_0)\cap\mathcal D^+_\Xi) \subset \mathcal S^+,\qquad
\psi(\mathcal B_{d_0}(M_0)\cap\mathcal D^-_\Xi) \subset \mathcal S^-,\\
&\psi(\mathcal B_{d_0}(M_0)\cap \Xi) \subset \{M\in\mathcal S: m_{nn}=0\},\\
&\psi(\mathcal B_{d_0}(M_0))\supset \mathcal B_d(\psi(M_0)),\\
&\|D\psi\|_{L^\infty(\mathcal B_{d_0}(M_0))}\le K,
  \qquad
  \|D\psi^{-1}\|_{L^\infty(\mathcal B_d(\psi(M_0)))}\le K,
\end{align*}
and $D\psi$ and $D\psi^{-1}$ have modulus of continuity $\tau_1$ in
$\mathcal B_{d_0}(M_0)$ and $\mathcal B_d(\psi(M_0))$, respectively.

We point out that if $\Xi$ is the boundary of a bounded $C^1$ domain, then
condition {\bf(S)} holds for $\Xi$.

\begin{proposition}\label{prop:curved-interface}
Let $F=F(M)$ satisfy \eqref{eq:ellipticity condition}, and let the surface
$\Xi$ satisfy {\bf(S)}.  Let $M_0\in\Xi$ and let $\psi$ straighten
$\mathcal B_{d_0}(M_0)\cap\Xi$ as above.  Assume that
\[
F\in C^1(\mathcal B_{d_0}(M_0)\cap\mathcal D_\Xi^+)
  \cap C^1(\mathcal B_{d_0}(M_0)\cap\mathcal D_\Xi^-)
\]
and that $DF$ has a modulus of continuity $\tau$ on
$\mathcal B_{d_0}(M_0)\cap\mathcal D_\Xi^+$ and on
$\mathcal B_{d_0}(M_0)\cap\mathcal D_\Xi^-$, respectively.  Then
\[
F\in\Bbb T^\omega(\overline{\mathcal B}_{R_0}(M_0),\lambda/2,2\Lambda),
\]
where $R_0$ and $\omega$ depend only on
$n$, $\lambda$, $\Lambda$, $K$, $d$, $d_0$, $\tau_1$, and $\tau$.
\end{proposition}

\begin{proof}
Define
\[
F^*(N)=F(\psi^{-1}(N)),\qquad N_0=\psi(M_0).
\]
Then
$F^*\in C^1(\mathcal B_d(N_0)\cap\mathcal S^+)
       \cap C^1(\mathcal B_d(N_0)\cap\mathcal S^-)$.
Moreover, the modulus of continuity of $DF^*$ and the $L^\infty$ bound of
$DF^*$ in the two half-balls depend only on
$n$, $\lambda$, $\Lambda$, $\tau$, $\tau_1$, and $K$; call this modulus
$\tau_2$.

Let $N_1\in\overline{\mathcal B_{d/2}(N_0)\cap\mathcal S^+}$.  Write
$N_1=(N_1',n_{1nn})$ and $\overline N_1=(N_1',0)$, where $n_{1nn}$ is the $(n,n)$-entry of $N_1$ and $N'_1$ are all non $(n,n)$-entries of $N_1$.  Following the construction in
the proof of Lemma \ref{lem:flat-interface}, with the hyperplane
$\{N\in\mathcal S:n_{nn}=0\}$, the tangent operator $T^*_{N_1}$ for $F^*$ is
given by
\[
T^*_{N_1}(N)=L^+(N)\quad\text{if } n_{nn}\ge0,
\qquad
T^*_{N_1}(N)=L^-(N)\quad\text{if } n_{nn}<0,
\]
where
\[
L^+(N)=F^*(N_1)+\langle D^+F^*(N_1),N-N_1\rangle
\]
and
\[
L^-(N)=F^*(N_1)+\langle D^+F^*(N_1),N-N_1\rangle
       +\bigl(D_{nn}^-F^*(\overline N_1)-D_{nn}^+F^*(N_1)\bigr)n_{nn}.
\]
Let $M_1=\psi^{-1}(N_1)$.  Denote by $L_\psi$ the linearization of $\psi$ at
$M_1$, namely
\[
L_\psi(M)=\left[\psi_{ij}(M_1)+\langle D\psi_{ij}(M_1),M-M_1\rangle\right]_{1\le i,j\le n}.
\]
Set
\[
T_{M_1}(M)=T^*_{N_1}(L_\psi(M)),\qquad M\in\mathcal S.
\]
Since $T^*_{N_1}$ is convex or concave in $N$ and $L_\psi$ is affine, $T_{M_1}$
is convex or concave in $M$.

Choose $R_0<\min\{d/(2K),d_0/2\}$.  If
$M_1\in \overline{\mathcal B_{R_0}(M_0)\cap\mathcal D_\Xi^+}$, then
$N_1=\psi(M_1)\in \overline{\mathcal B_{KR_0}(N_0)\cap S^+}$.  The estimates \eqref{eq:flat-tangent-plus} and \eqref{eq:flat-tangent-minus} yield
\begin{equation}\label{eq:curved-tangent-star}
|F^*(N)-T^*_{N_1}(N)|\le C\tau_2(|N-N_1|)|N-N_1|
\qquad\forall N\in \overline{\mathcal B_d(N_0)}.
\end{equation}
Since $\psi(\mathcal B_{2R_0}(M_0))\subset\mathcal B_{2KR_0}(N_0)$ and
$F(M)=F^*(\psi(M))$, it follows from \eqref{eq:curved-tangent-star} that
\begin{equation}\label{eq:curved-tangent-one}
|F(M)-T^*_{N_1}(\psi(M))|
\le C\tau_2(K|M-M_1|)|M-M_1|
\qquad\forall M\in\mathcal B_{2R_0}(M_0).
\end{equation}
By the Lipschitz continuity of $T^*_{N_1}$ and the $C^1$ control of $\psi$,
\begin{equation}\label{eq:curved-tangent-two}
|T^*_{N_1}(\psi(M))-T_{M_1}(M)|
\le C|\psi(M)-L_\psi(M)|
\le C\tau_1(|M-M_1|)|M-M_1|
\end{equation}
for $M\in\mathcal B_{d_0}(M_0)$.  From \eqref{eq:curved-tangent-one},
\eqref{eq:curved-tangent-two}, and the Lipschitz continuity of $F$ and
$T_{M_1}$, we obtain that $T_{M_1}$ is tangent to $F$ at $M_1$ with a modulus
$\omega$ depending only on the stated structural data and on $R_0$.

It remains to prove the ellipticity of $T_{M_1}$.  Since
$T^*_{N_1}=\max\{L^+,L^-\}$ or $T^*_{N_1}=\min\{L^+,L^-\}$, it is enough to
prove that $L^+\circ L_\psi$ and $L^-\circ L_\psi$ have ellipticity constants
$\lambda/2$ and $2\Lambda$.

Let $M,\Delta M\in\mathcal S$ with $\Delta M\ge0$.  Choose
$\{M_s\}_{s=2}^\infty\subset\mathcal B_{R_0}(M_0)\cap\mathcal D_\Xi^+$ such
that $M_s\to M_1$.  By the continuity of $D^+F^*$ in
$\overline{\mathcal B_d(N_0)\cap\mathcal S^+}$,
\begin{equation}\label{eq:curved-plus-chain-limit}
\begin{aligned}
&L^+(L_\psi(M+\Delta M))-L^+(L_\psi(M))\\
&=\left\langle D^+F^*(N_1),
   \left[\langle D\psi_{ij}(M_1),\Delta M\rangle\right]\right\rangle\\
&=\lim_{s\to\infty}\left\langle D^+F^*(\psi(M_s)),
   \left[\langle D\psi_{ij}(M_s),\Delta M\rangle\right]\right\rangle.
\end{aligned}
\end{equation}
By the chain rule,
\begin{equation}\label{eq:curved-plus-chain-rule}
\begin{aligned}
&\left\langle D^+F^*(\psi(M_s)),
   \left[\langle D\psi_{ij}(M_s),\Delta M\rangle\right]\right\rangle  \\
&=\sum_{k,\ell}\sum_{i,j}D^+_{ij}F^*(\psi(M_s))D_{k\ell}\psi_{ij}(M_s)\Delta m_{k\ell}\\
&=\sum_{k,\ell}D_{k\ell}(F^*\circ\psi)(M_s)\Delta m_{k\ell}
 =\langle DF(M_s),\Delta M\rangle.
\end{aligned}
\end{equation}
It follows from the ellipticity of $F$ that
\[
\lambda\|\Delta M\|\le L^+(L_\psi(M+\Delta M))-L^+(L_\psi(M))
\le\Lambda\|\Delta M\|.
\]
Thus $L^+\circ L_\psi$ has ellipticity constants $\lambda$ and $\Lambda$.

For the other affine branch, observe that for $N,\Delta N\in\mathcal S$,
\begin{align*}
L^-(N+\Delta N)-L^-(N)
&=\langle D^+F^*(N_1),\Delta N\rangle
  +\bigl(D^-_{nn}F^*(\overline N_1)-D^+_{nn}F^*(N_1)\bigr)\Delta n_{nn}\\
&=\sum_{(i,j)\ne(n,n)}\bigl(D_{ij}F^*(N_1)-D_{ij}F^*(\overline N_1)\bigr)
   \Delta n_{ij}+\langle D^-F^*(\overline N_1),\Delta N\rangle.
\end{align*}
Therefore, for $M,\Delta M\in\mathcal S$ with $\Delta M\ge0$,
\[
L^-(L_\psi(M+\Delta M))-L^-(L_\psi(M))
=E+\left\langle D^-F^*(\overline N_1),
\left[\langle D\psi_{ij}(M_1),\Delta M\rangle\right]\right\rangle,
\]
where
\begin{equation}\label{eq:curved-error-one}
|E|\le C\tau_2(KR_0)K|\Delta M|.
\end{equation}
Since $N_1,\overline N_1\in \overline{\mathcal B_{KR_0}(N_0)}$, we have
$\overline M_1=\psi^{-1}(\overline N_1)\in \overline{\mathcal B_{K^2R_0}(M_0)}\cap\Xi$.
We rewrite
\begin{align*}
&\left\langle D^-F^*(\overline N_1),
\left[\langle D\psi_{ij}(M_1),\Delta M\rangle\right]\right\rangle\\
&=\left\langle D^-F^*(\overline N_1),
\left[\langle D\psi_{ij}(\overline M_1),\Delta M\rangle\right]\right\rangle\\
&\quad+\left\langle D^-F^*(\overline N_1),
\left[\langle D\psi_{ij}(M_1)-D\psi_{ij}(\overline M_1),\Delta M\rangle\right]\right\rangle.
\end{align*}
If $K^2R_0\le d_0/2$, then
\begin{equation}\label{eq:curved-error-two}
\left|\left\langle D^-F^*(\overline N_1),
\left[\langle D\psi_{ij}(M_1)-D\psi_{ij}(\overline M_1),\Delta M\rangle\right]\right\rangle\right|
\le C\tau_1(2K^2R_0)|\Delta M|.
\end{equation}
We claim that
\begin{equation}\label{eq:curved-minus-main}
\lambda\|\Delta M\|\le
\left\langle D^-F^*(\overline N_1),
\left[\langle D\psi_{ij}(\overline M_1),\Delta M\rangle\right]\right\rangle
\le\Lambda\|\Delta M\|.
\end{equation}
To prove this, choose
$\{M_s^-\}\subset\mathcal B_{K^2R_0}(M_0)\cap\mathcal D_\Xi^-$ such that
$M_s^-\to\overline M_1$.  As in \eqref{eq:curved-plus-chain-limit} and
\eqref{eq:curved-plus-chain-rule}, the chain rule gives
\[
\left\langle D^-F^*(\psi(M_s^-)),
\left[\langle D\psi_{ij}(M_s^-),\Delta M\rangle\right]\right\rangle
=\langle DF(M_s^-),\Delta M\rangle.
\]
Letting $s\to\infty$ and using the ellipticity of $F$ proves
\eqref{eq:curved-minus-main}.

In view of \eqref{eq:curved-error-one}--\eqref{eq:curved-minus-main}, if $R_0$ is sufficiently small, then
$L^-\circ L_\psi$ has ellipticity constants $\lambda/2$ and $2\Lambda$.  Hence
$T_{M_1}$ is uniformly elliptic with constants $\lambda/2$ and $2\Lambda$.

The argument above gives the desired tangent operator for
$M_1\in \overline{\mathcal B_{R_0}(M_0)\cap\mathcal D_\Xi^+}$.  If
$M_1\in \overline{\mathcal B_{R_0}(M_0)\cap\mathcal D_\Xi^-}$, we apply the same argument
after replacing $\psi$ by the straightening map $\widehat\psi$ obtained by
composing $\psi$ with the reflection in $\mathcal S$ with respect to the
$m_{nn}$-axis; that is,
\[
\widehat\psi_{ij}=\psi_{ij}\quad\text{for }(i,j)\ne(n,n),
\qquad
\widehat\psi_{nn}=-\psi_{nn}.
\]
This completes the proof.
\end{proof}

We can now finish the verification of \hyperlink{A4}{\bf(A4)\rm} case.  Fix
$x_0\in\Omega$.  Since the surfaces $\Xi_i$ are separated by a distance greater
than $2d_0$, any ball
of radius $d_0$ meets at most one of them.  If $M_0$ lies within an
$R_0$-neighborhood of some $\Xi_i$, then Proposition \ref{prop:curved-interface}
applied to that surface gives the tangent operator property at $M_0$.  If
$M_0$ is not within this neighborhood, then $\mathcal B_{R_0}(M_0)$ is contained
in one of the domains $\mathcal D_i$, and $F(\cdot,x_0)$ is tangent at $M_0$ to
its affine linearization by Taylor's theorem and the modulus of continuity of
$D_MF$ in that domain.  Combining the two cases and using the Lipschitz
continuity of $F$ and of the corresponding tangent operators for distances
larger than $R_0$, we obtain
\[
F(\cdot,x_0)\in\Bbb T^\omega(\mathcal S,\lambda/2,2\Lambda),
\]
where $\omega$ depends only on the structural data in \hyperlink{A4}{\bf(A4)\rm}.
Therefore \hyperlink{A4}{\bf(A4)\rm} implies \hyperlink{H2}{\bf(H2)\rm}.

\setcounter{equation}{0}
\section{Isaacs equations }\label{sec:Isaacs equations}

The purpose of this section is to prove Theorems \ref{thm:Theorem 1.5}
and \ref{thm:Theorem 1.6}.

\subsection{Locally smooth enveloping and (H1)}\label{subsec:locally-smooth-enveloping}
In this subsection, we consider necessary and sufficient conditions for
locally smooth enveloping and apply them to establish Theorem
\ref{thm:Theorem 1.5}.

\begin{lemma}\label{lem:support-function-gap}
Let $F$ be a $C^1$ convex function with $DF$ bounded on $\mathcal S$.
Then, for any $R,\theta>0$, there exists $h=h(R,\theta)>0$ such that if
two support functions
\[
L_i(M)=\text{tr}(A_iM)+b_i,\qquad i=1,2,
\]
of $F$ satisfy $\|A_2-A_1\|\ge\theta$ and $L_1(N)=L_2(N)$ for some
$\|N\|\le R$, then
\[
F(N)-L_1(N)\ge h.
\]
\end{lemma}

\begin{proof}
Otherwise, there exist $L_{ik}(M)=\text{tr}(A_{ik}M)+b_{ik}$ and $N_k$,
$k\ge1$, $i=1,2$, such that
\[
\|A_{1k}-A_{2k}\|\ge\theta,\qquad
\|N_k\|\le R,\qquad L_{1k}(N_k)=L_{2k}(N_k),
\]
but
\[
0\le F(N_k)-L_{ik}(N_k)\le\frac1k,\qquad i=1,2.
\]
Since $DF$ is bounded, the sequences $\{A_{1k}\}$ and $\{A_{2k}\}$ are
bounded.  Consequently, $\{b_{1k}\}$ and $\{b_{2k}\}$ are also bounded.
Passing to a subsequence, we may assume that
\[
(A_{1k},b_{1k})\longrightarrow(A_1,b_1),\qquad
(A_{2k},b_{2k})\longrightarrow(A_2,b_2),\qquad N_k\longrightarrow N_0.
\]
Set $L_i(M)=\text{tr}(A_iM)+b_i$, $i=1,2$.  Obviously,
\[
L_i(M)=\lim_{k\to\infty}L_{ik}(M)\le F(M)
\quad\text{for }M\in\mathcal S,
\]
and
\[
F(N_0)-L_i(N_0)
=\lim_{k\to\infty}\bigl(F(N_k)-L_{ik}(N_k)\bigr)=0.
\]
Because $\|A_1-A_2\|\ge\theta$, $L_1$ and $L_2$ are two distinct support
functions of $F$ at $N_0$.  This contradicts the differentiability of $F$.
\end{proof}

Given $A\in\mathcal S$ and $b\in\R$, set
\[
L_{A,b}(M)=\text{tr}(AM)+b,\qquad M\in\mathcal S.
\]
Let $\mathcal E\subset\mathcal S\times\R$ be bounded.  Obviously,
\begin{equation}\label{eq:supremum-over-closure}
\sup_{(A,b)\in\mathcal E}L_{A,b}(M)
=\max_{(A,b)\in\overline{\mathcal E}}L_{A,b}(M).
\end{equation}
Denote
\[
\mathcal L_{\mathcal E}=\{L_{A,b}:(A,b)\in\mathcal E\}.
\]
We introduce the following locally smooth enveloping condition for
$\mathcal L_{\mathcal E}$ and $\mathcal E$:

\medskip
\begin{enumerate}
\item[\hypertarget{LSE}{\bf(LSE)\rm}]
There exists a modulus
$h=h_{\mathcal L_{\mathcal E}}:\R^+\times\R^+\longrightarrow\R^+$
such that, if $L_{A_i,b_i}\in\mathcal L_{\mathcal E}$, $i=1,2$,
\[
\|A_2-A_1\|\ge\theta>0,\qquad
L_{A_1,b_1}(M_0)=L_{A_2,b_2}(M_0)
\]
at some $\|M_0\|\le R$, then there exists
$L_{A_3,b_3}\in\mathcal L_{\mathcal E}$ satisfying
\[
L_{A_3,b_3}(M_0)-L_{A_1,b_1}(M_0)\ge h(R,\theta).
\]
\end{enumerate}
\medskip

\begin{lemma}\label{lem:LSE-closure}
Let $\mathcal E$ be bounded.  Assume that \hyperlink{LSE}{\bf(LSE)\rm}
holds for $\mathcal L_{\mathcal E}$.  Then
\[
\mathcal L_{\overline{\mathcal E}}
=\{L_{A,b}:(A,b)\in\overline{\mathcal E}\}
\]
satisfies \hyperlink{LSE}{\bf(LSE)\rm} with modulus
\[
h_{\mathcal L_{\overline{\mathcal E}}}(R,\theta)
=h_{\mathcal L_{\mathcal E}}(R+1,\theta/2).
\]
\end{lemma}

\begin{proof}
Let $L_{A_i,b_i}\in\mathcal L_{\overline{\mathcal E}}$, $i=1,2$, and
$\|M_0\|\le R$.  Assume that
\[
\|A_2-A_1\|\ge\theta,\qquad L_{A_1,b_1}(M_0)=L_{A_2,b_2}(M_0).
\]
There exist sequences $\{(A_{1k},b_{1k})\}$,
$\{(A_{2k},b_{2k})\}\subset\mathcal E$ such that
\[
(A_{ik},b_{ik})\longrightarrow(A_i,b_i),\qquad i=1,2.
\]
Clearly, $\|A_{2k}-A_{1k}\|\ge\theta/2$ for large $k$.  Set
\[
\varepsilon_k
=\bigl(L_{A_1,b_1}-L_{A_{1k},b_{1k}}\bigr)(M_0)
-\bigl(L_{A_2,b_2}-L_{A_{2k},b_{2k}}\bigr)(M_0).
\]
Then $\varepsilon_k\to0$ and
\[
L_{A_{1k},b_{1k}}(M_0)+\varepsilon_k=L_{A_{2k},b_{2k}}(M_0).
\]
Write $A_i=[a^i_{\ell m}]$ and $A_{ik}=[a^{ik}_{\ell m}]$.
There is an index $(\ell,m)$ such that
$|a^1_{\ell m}-a^2_{\ell m}|\ge\theta/n$.
Hence, $|a^{1k}_{\ell m}-a^{2k}_{\ell m}|\ge\theta/(2n)$ for large $k$.
Let $B_{\ell m}$ be the matrix whose $(\ell,m)$ and $(m,\ell)$ entries are
$1$ and whose other entries are $0$.  We then have
\[
\bigl(L_{A_{1k},b_{1k}}-L_{A_{2k},b_{2k}}\bigr)
(M_0+\widehat\varepsilon_kB_{\ell m})=0,
\]
where
\[
\widehat\varepsilon_k
=\frac{\varepsilon_k}{2(a^{1k}_{\ell m}-a^{2k}_{\ell m})}
\quad\text{if }\ell\ne m,\qquad
\widehat\varepsilon_k
=\frac{\varepsilon_k}{a^{1k}_{\ell m}-a^{2k}_{\ell m}}
\quad\text{if }\ell=m.
\]
Since \hyperlink{LSE}{\bf(LSE)\rm} holds for $\mathcal E$, for large $k$
there exists $L_{A_{3k},b_{3k}}\in\mathcal L_{\mathcal E}$ such that
\begin{equation}\label{eq:LSE-closure-gap}
\bigl(L_{A_{3k},b_{3k}}-L_{A_{1k},b_{1k}}\bigr)
(M_0+\widehat\varepsilon_kB_{\ell m})
\ge h_{\mathcal L_{\mathcal E}}(R+1,\theta/2).
\end{equation}
Since $\overline{\mathcal E}$ is compact, after passing to a subsequence
we may assume that $(A_{3k},b_{3k})\to(A_3,b_3)\in\overline{\mathcal E}$.
Letting $k\to\infty$ in \eqref{eq:LSE-closure-gap} gives
\[
L_{A_3,b_3}(M_0)-L_{A_1,b_1}(M_0)
\ge h_{\mathcal L_{\mathcal E}}(R+1,\theta/2).
\]
\end{proof}

\begin{lemma}\label{lem:LSE-differentiability}
Let $\mathcal E$ be bounded and assume that \hyperlink{LSE}{\bf(LSE)\rm}
holds for $\mathcal E$.  Then
\[
F(M)=\sup_{(A,b)\in\mathcal E}L_{A,b}(M)
\]
belongs to $C^1_{\mathrm{loc}}(\mathcal S)$.
\end{lemma}

\begin{proof}
In view of \eqref{eq:supremum-over-closure} and Lemma
\ref{lem:LSE-closure}, we may assume that $\mathcal E$ is compact.

First, we show that for each $M_0\in\mathcal S$ there is a unique
$(A(M_0),b(M_0))\in\mathcal E$ such that $L_{A(M_0),b(M_0)}$ supports
$F$ at $M_0$.  The compactness of $\mathcal E$ gives existence.
If $L_{A'_0,b'_0}\in\mathcal L_{\mathcal E}$ is another support function
of $F$ at $M_0$, then $A(M_0)\ne A'_0$.  Since
\[
L_{A(M_0),b(M_0)}(M_0)=L_{A'_0,b'_0}(M_0),
\]
\hyperlink{LSE}{\bf(LSE)\rm}, with $R=\|M_0\|+1$ and
$\theta=\|A(M_0)-A'_0\|$, gives $L_{A,b}\in\mathcal L_{\mathcal E}$
such that
\[
L_{A,b}(M_0)\ge L_{A'_0,b'_0}(M_0)+h_{\mathcal L_{\mathcal E}}(R,\theta)
>F(M_0),
\]
a contradiction.

Second, the map
\[
\mathcal N:\mathcal S\longrightarrow\mathcal E,\qquad
\mathcal N(M)=(A(M),b(M)),
\]
is continuous.  Indeed, let $M_k\to M_0$.  By compactness, any subsequence
has a further subsequence such that
\[
(A(M_{k_j}),b(M_{k_j}))\longrightarrow(A_0,b_0)\in\mathcal E.
\]
Because $L_{A(M_{k_j}),b(M_{k_j})}(M)\le F(M)$ with equality at $M_{k_j}$,
letting $j\to\infty$ shows that $L_{A_0,b_0}$ supports $F$ at $M_0$.
The uniqueness proved above yields $(A_0,b_0)=(A(M_0),b(M_0))$.

Finally, we show that $F$ is differentiable at every $M_0\in\mathcal S$ and
\begin{equation}\label{eq:LSE-tangency}
F(M)=L_{A(M_0),b(M_0)}(M)+o(\|M-M_0\|)
\quad\text{as }M\to M_0.
\end{equation}
Let $M\ne M_0$.  Then
\[
\bigl(L_{A(M),b(M)}-L_{A(M_0),b(M_0)}\bigr)(M)\ge0
\]
and
\[
\bigl(L_{A(M),b(M)}-L_{A(M_0),b(M_0)}\bigr)(M_0)\le0.
\]
If the first expression is zero, then $F(M)=L_{A(M_0),b(M_0)}(M)$ and
\eqref{eq:LSE-tangency} holds.  If the second is zero, uniqueness gives
$L_{A(M),b(M)}=L_{A(M_0),b(M_0)}$, and again \eqref{eq:LSE-tangency}
holds.  In the remaining case, there exists $0<t<1$ such that, with
$M_t=(1-t)M_0+tM$,
\[
\bigl(L_{A(M),b(M)}-L_{A(M_0),b(M_0)}\bigr)(M_t)=0.
\]
Therefore,
\begin{align*}
F(M)-L_{A(M_0),b(M_0)}(M)
&=\bigl(L_{A(M),b(M)}-L_{A(M_0),b(M_0)}\bigr)(M)\\
&\quad-\bigl(L_{A(M),b(M)}-L_{A(M_0),b(M_0)}\bigr)(M_t)\\
&=\text{tr}\bigl((A(M)-A(M_0))(M-M_t)\bigr)\\
&=o(\|M-M_0\|).
\end{align*}
The lemma follows from $DF(M)=A(M)$ and the continuity of $A(M)$.
\end{proof}

\begin{proposition}\label{prop:LSE-uniform-modulus}
Assume that $\mathcal E\subset\overline{\mathcal B}_{C_0}\times[-C_0,C_0]$
and $\mathcal E$ satisfies \hyperlink{LSE}{\bf(LSE)\rm}.  Let
\[
F=\sup_{(A,b)\in\mathcal E}L_{A,b}.
\]
Then the modulus of continuity of $DF$ in $\overline{\mathcal B}_R$
depends only on $R$, $C_0$, $n$, and the modulus $h$ in
\hyperlink{LSE}{\bf(LSE)\rm}.  More precisely, for every $R>0$ and
$\varepsilon>0$, there exists
$\widehat\delta=\widehat\delta(\varepsilon,R,C_0,n,h)>0$ such that
\[
\|DF(M_2)-DF(M_1)\|\le\varepsilon
\]
whenever $M_1,M_2\in\overline{\mathcal B}_R$ and
$\|M_2-M_1\|\le\widehat\delta$.
\end{proposition}

\begin{proof}
Suppose otherwise.  Then there exist $R_0>0$, $\varepsilon_0>0$, sequences
$\{M_{1k}\}$, $\{M_{2k}\}\subset\overline{\mathcal B}_{R_0}$, and sets
\[
\mathcal E_k\subset\overline{\mathcal B}_{C_0}\times[-C_0,C_0]
\]
satisfying \hyperlink{LSE}{\bf(LSE)\rm} with the same modulus $h$, such that
\[
\|M_{2k}-M_{1k}\|\longrightarrow0,\qquad
\|DF_k(M_{2k})-DF_k(M_{1k})\|>\varepsilon_0,
\]
where
\[
F_k=\sup_{(A,b)\in\mathcal E_k}L_{A,b}.
\]
By the Ascoli--Arzel\`a theorem, after passing to a subsequence,
$F_k\to F$ locally uniformly in $\mathcal S$.

Let $\mathcal D=\{M_j:j\ge1\}$ be a countable dense subset of $\mathcal S$.
At each $M_j$, Lemma \ref{lem:LSE-differentiability} gives a unique support
function $L_{A^k_j,b^k_j}$ for $F_k$, where
$(A^k_j,b^k_j)\in\overline{\mathcal E}_k$.
By a diagonalization argument, we may assume that
\[
(A^k_j,b^k_j)\longrightarrow(A_j,b_j)\quad\text{as }k\to\infty
\]
for every fixed $j$.  Passing to the limit in $L_{A^k_j,b^k_j}\le F_k$,
with equality at $M_j$, shows that $L_{A_j,b_j}$ supports $F$ at $M_j$.
Consequently,
\[
F(M_\ell)=\sup_{j\ge1}L_{A_j,b_j}(M_\ell)
\quad\text{for every }M_\ell\in\mathcal D.
\]
By continuity,
\[
F=\sup_{j\ge1}L_{A_j,b_j}\quad\text{on }\mathcal S.
\]

We next show that $\mathcal E_0=\{(A_j,b_j):j\ge1\}$ satisfies
\hyperlink{LSE}{\bf(LSE)\rm}.  Let $(A_j,b_j),(A_\ell,b_\ell)\in\mathcal E_0$
satisfy
\[
\|A_j-A_\ell\|\ge\theta,\qquad
L_{A_j,b_j}(M_0)=L_{A_\ell,b_\ell}(M_0)
\]
for some $\|M_0\|\le R$.  Set
\[
\varepsilon_k
=\bigl(L_{A_j,b_j}-L_{A^k_j,b^k_j}\bigr)(M_0)
-\bigl(L_{A_\ell,b_\ell}-L_{A^k_\ell,b^k_\ell}\bigr)(M_0).
\]
Then $\varepsilon_k\to0$,
\[
L_{A^k_j,b^k_j}(M_0)+\varepsilon_k=L_{A^k_\ell,b^k_\ell}(M_0),
\]
and $\|A^k_j-A^k_\ell\|\ge\theta/2$ for large $k$.
Write $A_j=[a^j_{im}]$, $A_\ell=[a^\ell_{im}]$,
$A^k_j=[a^{kj}_{im}]$, and $A^k_\ell=[a^{k\ell}_{im}]$.
For some $(i,m)$,
\[
|a^j_{im}-a^\ell_{im}|\ge\theta/n,
\]
and hence $|a^{kj}_{im}-a^{k\ell}_{im}|\ge\theta/(2n)$ for large $k$.
As in Lemma \ref{lem:LSE-closure}, there exists $\widehat\varepsilon_k\to0$
such that
\[
\bigl(L_{A^k_j,b^k_j}-L_{A^k_\ell,b^k_\ell}\bigr)
(M_0+\widehat\varepsilon_kB_{im})=0.
\]
By Lemma \ref{lem:LSE-closure}, $\overline{\mathcal E}_k$ satisfies
\hyperlink{LSE}{\bf(LSE)\rm}; hence, for large $k$, there exists
$(A^k,b^k)\in\overline{\mathcal E}_k$ such that
\[
\bigl(L_{A^k,b^k}-L_{A^k_j,b^k_j}\bigr)
(M_0+\widehat\varepsilon_kB_{im})\ge h(R+2,\theta/4).
\]
Since $F_k\ge L_{A^k,b^k}$, passing to the limit yields
\[
F(M_0)-L_{A_j,b_j}(M_0)\ge h(R+2,\theta/4).
\]
For $M_{j_0}\in\mathcal D$ sufficiently close to $M_0$,
\[
F(M_0)-L_{A_{j_0},b_{j_0}}(M_0)
=F(M_0)-F(M_{j_0})+\text{tr}\bigl(A_{j_0}(M_{j_0}-M_0)\bigr)
\]
is sufficiently small, and therefore
\[
L_{A_{j_0},b_{j_0}}(M_0)-L_{A_j,b_j}(M_0)
\ge\frac12 h(R+2,\theta/4).
\]
Thus \hyperlink{LSE}{\bf(LSE)\rm} holds for $\mathcal E_0$, and Lemma
\ref{lem:LSE-differentiability} gives $F\in C^1_{\mathrm{loc}}(\mathcal S)$.

Finally, let $L_{A_{ik},b_{ik}}$ be the support function of $F_k$ at $M_{ik}$,
where $(A_{ik},b_{ik})\in\overline{\mathcal E}_k$, $i=1,2$.
Passing to further subsequences, we may assume that
\[
M_{ik}\longrightarrow M^*,\qquad
(A_{ik},b_{ik})\longrightarrow(A_i^*,b_i^*),\qquad i=1,2.
\]
Both $L_{A_1^*,b_1^*}$ and $L_{A_2^*,b_2^*}$ support $F$ at $M^*$.
Since $F\in C^1_{\mathrm{loc}}(\mathcal S)$,
\[
A_1^*=DF(M^*)=A_2^*.
\]
This contradicts
\[
\|A_{2k}-A_{1k}\|=\|DF_k(M_{2k})-DF_k(M_{1k})\|>\varepsilon_0.
\]
\end{proof}

\begin{proof}[Proof of Theorem \ref{thm:Theorem 1.5}]
Under \hyperlink{I1}{\bf(I1)\rm}, Proposition
\ref{prop:LSE-uniform-modulus} shows that
\[
F_\gamma(M,x)
=\sup_{\beta\in\mathcal B}
\bigl[\text{tr}(A_{\beta\gamma}(x)M)-f_{\beta\gamma}(x)\bigr]
=\sup_{(A,b)\in\mathcal E_{\gamma,x}}L_{A,b}(M)
\]
belongs to $C^1_{\mathrm{loc}}(\mathcal S)$ and that
$D_MF_\gamma(\cdot,x)$ has a uniform modulus of continuity on each
$\overline{\mathcal B}_k$, $k\ge1$.  Therefore,
\[
F(M,x)=\inf_{\gamma\in\mathcal G}F_\gamma(M,x)
\]
satisfies a local version of \hyperlink{A2}{\bf(A2)\rm}, and hence
\hyperlink{H1}{\bf(H1)\rm}.

Under \hyperlink{I2}{\bf(I2)\rm}, write
\begin{align*}
F(M,x)
&=\inf_{\gamma\in\mathcal G}
\bigl[\text{tr}(A^{II}_\gamma(x)M)-f^{II}_\gamma(x)\bigr]\\
&\quad+\sup_{\beta\in\mathcal B}
\bigl[\text{tr}(A^I_\beta(x)M)-f^I_\beta(x)\bigr]\\
&=F_\cap(M,x)+F_C(M,x).
\end{align*}
Applying Proposition \ref{prop:LSE-uniform-modulus} to $F_C$ shows that
$D_MF_C(\cdot,x)$ has a uniform modulus of continuity on each
$\overline{\mathcal B}_k$.  Thus $F$ satisfies a local analogue of
\hyperlink{A1}{\bf(A1)\rm}, and \hyperlink{H1}{\bf(H1)\rm} follows.
\end{proof}

\subsection{Isaacs equations satisfying (H2)}\label{subsec:global-smooth-enveloping}
Let $\mathcal E\subset\mathcal S\times\R$ be compact.  We consider the
following global smooth enveloping condition for $\mathcal L_{\mathcal E}$
and $\mathcal E$:

\medskip
\begin{enumerate}
\item[\hypertarget{GSE}{\bf(GSE)\rm}]
There exists a nondecreasing function $\tau$ on $(0,\infty)$ with
$\lim_{t\to0^+}\tau(t)=0$ such that, if $L_1,L_2\in\mathcal L_{\mathcal E}$
support
\[
\sup_{L\in\mathcal L_{\mathcal E}}L(M)
\]
at $M_1,M_2$, respectively, and $\|M_2-M_1\|\le1$, then
\[
|L_2(M_i)-L_1(M_i)|
\le\tau(\|M_2-M_1\|)\|M_2-M_1\|,\qquad i=1,2.
\]
\end{enumerate}
\medskip

\begin{proposition}\label{prop:GSE-uniform-modulus}
Let $\mathcal E$ be compact and satisfy \hyperlink{GSE}{\bf(GSE)\rm}.
Set
\[
F(M)=\sup_{L\in\mathcal L_{\mathcal E}}L(M).
\]
Then $F\in C^1_{\mathrm{loc}}(\mathcal S)$ and, for $M_1,M_2\in\mathcal S$
with $\|M_2-M_1\|\le1$,
\begin{equation}\label{eq:GSE-gradient-modulus}
\|DF(M_2)-DF(M_1)\|\le8\tau(\|M_2-M_1\|).
\end{equation}
\end{proposition}

\begin{proof}
Given $M_0\in\mathcal S$, there exists $L_0\in\mathcal L_{\mathcal E}$
such that $L_0\le F$, with equality at $M_0$.
For $\|M-M_0\|\le1$, let $L\in\mathcal L_{\mathcal E}$ support $F$ at $M$.
By \hyperlink{GSE}{\bf(GSE)\rm},
\begin{equation}\label{eq:GSE-tangency}
|F(M)-L_0(M)|=|L(M)-L_0(M)|
\le\tau(\|M-M_0\|)\|M-M_0\|.
\end{equation}
This gives the differentiability of $F$ on $\mathcal S$.

To establish the continuity of $DF$, let $R=\|M_1-M_0\|\le1/2$ and
$\|M-M_1\|\le R$.  From \eqref{eq:GSE-tangency},
\begin{equation}\label{eq:GSE-expansion-zero}
F(M)=F(M_0)+\text{tr}\bigl(DF(M_0)(M-M_0)\bigr)+E_0,
\quad |E_0|\le\tau(2R)\,2R,
\end{equation}
and
\begin{align}
F(M)&=F(M_1)+\text{tr}\bigl(DF(M_1)(M-M_1)\bigr)+E_1\notag\\
&=F(M_0)+\text{tr}\bigl(DF(M_0)(M_1-M_0)\bigr)+E_2
+\text{tr}\bigl(DF(M_1)(M-M_1)\bigr)+E_1,\label{eq:GSE-expansion-one}
\end{align}
where $|E_1|\le\tau(R)R$ and $|E_2|\le\tau(R)R$.
Comparing \eqref{eq:GSE-expansion-zero} and \eqref{eq:GSE-expansion-one}
yields
\[
\text{tr}\bigl((DF(M_0)-DF(M_1))(M-M_1)\bigr)=E_1+E_2-E_0.
\]
Choosing a suitable $M$, we obtain
\[
\|DF(M_0)-DF(M_1)\|\le4\tau(2R)=4\tau(2\|M_0-M_1\|).
\]
For $\|M_2-M_1\|\le1$, set $M_0^*=(M_1+M_2)/2$ and apply the preceding
estimate to $DF(M_0^*)-DF(M_i)$, $i=1,2$, to obtain
\eqref{eq:GSE-gradient-modulus}.
\end{proof}

\begin{proof}[Proof of Theorem \ref{thm:Theorem 1.6}]
The proofs under \hyperlink{I3}{\bf(I3)\rm} and \hyperlink{I4}{\bf(I4)\rm}
are analogous to the corresponding parts of the proof of Theorem
\ref{thm:Theorem 1.5}, using \hyperlink{GSE}{\bf(GSE)\rm} and Proposition
\ref{prop:GSE-uniform-modulus} in place of \hyperlink{LSE}{\bf(LSE)\rm}
and Proposition \ref{prop:LSE-uniform-modulus}.

Now consider \hyperlink{I5}{\bf(I5)\rm}.  The operator $F$ is elliptic and
Lipschitz on $\mathcal S$.  Denote
\[
L_{\beta\gamma}(M)=\text{tr}(A_{\beta\gamma}M)-f_{\beta\gamma}.
\]
We claim that $F$ is piecewise linear in $\mathcal S$, associated with
\[
\{\Pi(\beta,\gamma,\beta',\gamma')\}
\quad\text{and}\quad
\{L_{\beta\gamma}:\beta\in\mathcal B,\ \gamma\in\mathcal G\}.
\]
More precisely, there exist a subfamily $\{\Pi_i\}_{i=1}^N$ of
$\{\Pi(\beta,\gamma,\beta',\gamma')\}$ and a subfamily
$\{L_j\}_{j=0}^N$ of $\{L_{\beta\gamma}\}$ such that the hyperplanes
$\{\Pi_i\}$ partition $\mathcal S$ into disjoint strips
$\{\mathcal D_j\}_{j=0}^N$, with
\[
F=L_j\quad\text{on }\mathcal D_j,\qquad L_j\ne L_{j+1}.
\]

The claim follows in two steps.  First, if $G$ is piecewise linear associated
with $\{\Pi(\beta,\gamma,\beta',\gamma')\}$ and $\{L_{\beta\gamma}\}$, then
\[
G^*=\sup\{G,L_{\beta_0\gamma_0}\}
\]
is also piecewise linear.  Indeed, if $G=L_j$ on a strip $\mathcal D_j$,
then $G^*$ is either linear on $\mathcal D_j$ or piecewise linear when the
hyperplane $L_j=L_{\beta_0\gamma_0}$ lies inside $\mathcal D_j$.
Adjacent sub-strips on which the same linear function occurs are merged.
Repeatedly applying this observation shows that
\[
F_\gamma(M)=\sup_{\beta\in\mathcal B}L_{\beta\gamma}(M)
\]
is piecewise linear and associated with the same hyperplanes and linear
functions.

Second, if $G$ and $G'$ are piecewise linear associated with these
hyperplanes and linear functions, then so is $\inf\{G,G'\}$.
Indeed, if $G=L_j$ on $\mathcal D_j$, $G'=L'_k$ on $\mathcal D'_k$, and
$\mathcal D_j\cap\mathcal D'_k\ne\emptyset$, then
\[
\inf\{G,G'\}=\inf\{L_j,L'_k\}
\]
is either linear or piecewise linear on $\mathcal D_j\cap\mathcal D'_k$.
Repeatedly applying this observation shows that
\[
F=\inf_{\gamma\in\mathcal G}F_\gamma
\]
is piecewise linear, proving the claim.

The piecewise linear representation of $F$ is unique, and
$\bigcup_{i=1}^N\Pi_i$ is the set where $F$ is nondifferentiable.
Let
\[
d_1=\min_{i\ne j}\dist(\Pi_i,\Pi_j).
\]
By Lemma \ref{lem:flat-interface},
\[
F\in\mathbb T^\omega(\mathcal S,\lambda/2,2\Lambda),\qquad\omega(t)=Ct,
\]
where $C=C(n,\lambda,\Lambda,d_1)$.  This proves the theorem.
\end{proof}

\begin{remark}\label{rem:variable-Isaacs-coefficients}
Condition \hyperlink{I5}{\bf(I5)\rm} can be modified in a straightforward
way to accommodate variable $A_{\beta\gamma}(x)$ and $f_{\beta\gamma}(x)$.
In this case, assume additionally that, for each $x\in\Omega$, the
hyperplanes $\Pi(\beta,\gamma,\beta',\gamma',x)$ given by
\[
\text{tr}(A_{\beta\gamma}(x)M)-f_{\beta\gamma}(x)
=\text{tr}(A_{\beta'\gamma'}(x)M)-f_{\beta'\gamma'}(x),
\]
on which $D_MF(\cdot,x)$ does not exist are separated by a distance not less
than a constant $d_1$.
\end{remark}


\end{document}